\documentclass[aop,MSNbibl,citesort,dvips]{arximspdf}
\usepackage{mathrsfs}
\usepackage{mathbh}
\usepackage{graphicx}

\doi{10.1214/11-AOP686}
\volume{40}
\issue{6}
\pubyear{2012}
\firstpage{2589}
\lastpage{2666}

\makeatletter

\newcommand{\lbr}{[[}
\newcommand{\rbr}{]]}

\newcommand{\eqref}[1]{(\ref{#1})}

\DeclareMathAlphabet{\mathbbl}  {U}{mt2hrb}{m}{n}
\SetMathAlphabet{\mathbbl}{bold}{U}{mt2hrb}{b}{n}

\newcommand{\mathbbm}{\mathbh}

\newtheorem{theorem}{Theorem}
\newproclaim{definition}[theorem]{Definition}
\newtheorem{proposition}[theorem]{Proposition}
\newtheorem{corollary}[theorem]{Corollary}
\newtheorem{lemma}[theorem]{Lemma}

\newcommand{\dist}{\operatorname{dist}}

\newcommand{\gw}{\mathrm{GW}_\xi}

\newcommand{\wt}{\widetilde}

\newcommand{\ov}{\overline}
\newcommand{\bF}{\mathsf{F}}
\newcommand{\bTheta}{\bolds{\Theta}}
\newcommand{\N}{\mathbb{N}}
\newcommand{\PP}{\mathbb{P}}
\newcommand{\E}{\mathbb{E}}
\newcommand{\R}{\mathbb{R}}
\newcommand{\Z}{\mathbb{Z}}
\newcommand{\T}{\mathscr{T}}
\newcommand{\M}{\mathscr{M}}
\newcommand{\Mw}{\mathscr{M}_{\mathrm{w}}}
\newcommand{\Tw}{\mathscr{T}_{\mathrm{w}}}
\newcommand{\bsigma}{\bolds{\sigma}}
\renewcommand{\P}{\mathcal{P}}
\newcommand{\sd}{{\mathcal{S}^{\downarrow}}}
\newcommand{\bT}{\mathsf{T}}
\newcommand{\bfT}{\mathbf{T}}
\newcommand{\bP}{\mathsf{P}}
\newcommand{\bQ}{\mathsf{Q}}
\newcommand{\bs}{\mathbf{s}}
\newcommand{\bX}{\mathbf{X}}
\newcommand{\bFF}{\mathbf{F}}
\newcommand{\FF}{\mathcal{F}}
\newcommand{\TT}{\mathcal{T}}
\renewcommand{\d}{\mathrm{d}}
\newcommand{\da}{\downarrow}
\newcommand{\eps}{\varepsilon}
\newcommand{\ind}{\mathbf{1}}

\makeatother

\begin{document}
\begin{frontmatter}

\title{Scaling limits of Markov branching trees with applications to
Galton--Watson and random unordered trees\thanksref{T1}}
\runtitle{Scaling limits of Markov branching trees}

\thankstext{T1}{Supported\vspace*{1pt} in part by ANR-08-BLAN-0190 and ANR-08-BLAN-0220-01.}

\begin{aug}
\author[A]{\fnms{B\'en\'edicte} \snm{Haas}\ead[label=e2]{haas@ceremade.dauphine.fr}} and
\author[B]{\fnms{Gr\'egory} \snm{Miermont}\corref{}\ead[label=e1]{Gregory.Miermont@math.u-psud.fr}}
\runauthor{B. Haas and G. Miermont}
\affiliation{Universit\'e Paris-Dauphine and Universit\'e Paris-Sud}
\address[A]{Ceremade\\
Universit\'e Paris-Dauphine\\
Place du Mar{\'e}chal de Lattre de Tassigny\\
75016 Paris\\
France\\
\printead{e2}} %adresu isvedimo komanda gale!
\address[B]{D\'epartement de Math\'ematiques d'Orsay\\
Universit\'e Paris-Sud\\
Bat. 425\\
Orsay Cedex, 91 405\\
France\\
\printead{e1}}
\end{aug}

% HISTORY:
\received{\smonth{6} \syear{2010}}
\revised{\smonth{5} \syear{2011}}

% ABSTRACT
%
\begin{abstract}
We consider a family of random trees satisfying a Markov branching
property. Roughly, this property says that the subtrees above some
given height are independent with a law that depends only on their
total size, the latter being either the number of leaves or
vertices. Such families are parameterized by sequences of
distributions on partitions of the integers that determine how the
size of a tree is distributed in its different subtrees. Under some
natural assumption on these distributions, stipulating that
``macroscopic'' splitting events are rare, we show that Markov
branching trees admit the so-called self-similar fragmentation trees
as scaling limits in the Gromov--Hausdorff--Prokhorov topology.

The main application of these results is that the scaling limit of
random uniform unordered trees is the Brownian continuum random
tree. This extends a result by Marckert--Miermont and fully proves a
conjecture by Aldous. We also recover, and occasionally extend,
results on scaling limits of consistent Markov branching models and
known convergence results of Galton--Watson trees toward the Brownian
and stable continuum random trees.
\end{abstract}

% KEYWORDS
%

\end{frontmatter}

\tableofcontents

%%%%%%%%%%%%%%%%%%%%%%%%%%
%s1 ###
\section{Introduction and main results}\label{sec1}
\label{MBTrees}
%%%%%%%%%%%%%%%%%%%%%%%%%%

The goal of this paper is to discuss the scaling limits of a model of
random trees satisfying a simple Markovian branching property that
was considered in different forms in
\cite{Devroye,Aldous96,BrDeMcLdlS,HMPW06}. Markov branching trees are
natural models of random trees defined in terms of discrete
fragmentation processes. The laws of these trees are indexed by an
integer $n$ giving the ``size'' of the tree, which leads us to
consider two distinct (but related) models, in which the sizes are,
respectively, the number of leaves and the number of vertices. We first
provide a slightly informal description of our results.

%
%f1 ###
\begin{figure}[b]

\includegraphics{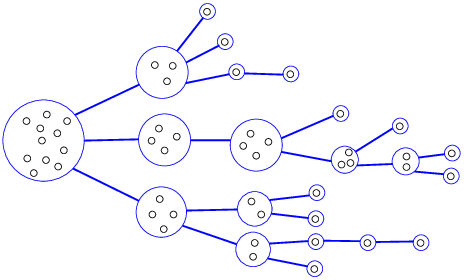}

\caption{A sample tree $T_{11}$. The first splitting arises with
probability $q_{11}(4,4,3)$.}
\label{fig:CMB}
\end{figure}

%, leaving
%a more formal approach to the rest of this introduction.

Let $q=(q_n,n\geq1)$ be a family of probability distributions,
respectively, on the set $\P_n$ of partitions of the integer $n$,
that is, of nonincreasing integer sequences with sum $n$. We assume that
$q_n$ does not assign mass $1$ to the trivial partition $(n)$
\[
q_n((n))<1 \qquad \mbox{for every } n\geq1 .
\]
In order that
this makes sense for $n=1$, we add an extra ``empty partition''
$\varnothing$ to~$\P_1$.

One constructs a random rooted tree with $n$ leaves according to the
following procedure. Start from a collection of $n$ indistinguishable
balls, and with probability $q_n(\lambda_1,\ldots,\lambda_p)$, split
the collection into $p$ sub-collections with
$\lambda_1,\ldots,\lambda_p$ balls. Note that there is a chance
$q_n((n))<1$ that the collection remains unchanged during this step of
the procedure. Then, re-iterate the splitting operation independently
for each sub-collection using this time the probability distributions
$q_{\lambda_1},\ldots,q_{\lambda_p}$. If a sub-collection consists of
a single ball, it can remain single with probability $q_1((1))$ or get
wiped out with probability $q_1(\varnothing)$. We continue the procedure
until all the balls are wiped out. There is a natural genealogy
associated with this process, which is a tree with $n$ leaves
consisting in the $n$ isolated balls just before they are wiped out,
and rooted at the initial collection of $n$ balls. See Figure
\ref{fig:CMB} for an illustration. We let $\bP^q_n$ be the law of this
tree.\vadjust{\goodbreak}

This construction can be seen as the most general form of \textit{splitting trees} of Broutin et al.~\cite{BrDeMcLdlS}, and was
referred to as trees having the so-called \textit{Markov branching
property} in~\cite{HMPW06}. There is also a variant of this
procedure that constructs a random tree with $n$ vertices rather than
$n$ leaves. This one does not need the hypothesis $q_n((n))<1$ for
$n\geq1$, and in fact we only assume $q_1((1))=1$ for consistency of
the description to follow. Informally, starting from a collection of
$n$ balls, we first remove a ball, split the $n-1$ remaining balls in
sub-collections with $\lambda_1,\ldots,\lambda_p$ balls with
probability $q_{n-1}((\lambda_1,\ldots,\lambda_p))$,\vspace*{1pt} and iterate
independently on sub-collections until no ball remains. We let
$\bQ^q_n$ be the law of the random tree associated to this procedure.

While most papers so far have been focusing on families of trees
having more structure, such as a consistency property when $n$ varies
\cite{Aldous96,HMPW06,Devroye} (with the notable exception of Broutin
et al.~\cite{BrDeMcLdlS}), the main\vspace*{1pt} goal of the present work is to
study the geometry of trees with laws $\bP^q_n$ or $\bQ^q_n$ as
$n\to\infty$ in a very general situation. The main assumption that we
make is that, as $n\to\infty$,

\begin{itemize}
\item[]``macroscopic'' splitting events of the form $n\to(ns_1,ns_2,\ldots
)\in
\P_n$ for a nonincreasing sequence $\bs=(s_1,s_2,\ldots)$ with sum
$1$ and such that $s_1<1-\eps$, for some $\eps\in(0,1)$, are rare
events, occurring with probability of order $n^{-\gamma}\nu_\eps(\d
\bs)$
for some $\gamma>0$, for some finite ``intensity'' measure $\nu_\eps$.
\end{itemize}

\noindent Note that the measures $\nu_\eps$ should satisfy a consistency
property as $\eps$ varies, and as $\eps$ goes to $0$, $\nu_\eps$
should increase to a possibly infinite measure $\nu$ on the set of
nonincreasing sequences with sum $1$. This means that splitting
events that only remove tiny parts from a large collection of balls
are allowed to remain more frequent than the order $n^{-\gamma}$.
Under this assumption, formalized in hypothesis (H) below, we
show in Theorem~\ref{sec:main-result-1} that a tree $T_n$
\label{tn}
with law
$\bP^q_n$, considered as a metric space by viewing its edges as being
real segments of lengths of order $n^{-\gamma}$, converges in
distribution toward a limiting structure $\TT_{\gamma,\nu}$, the
so-called self-similar fragmentation tree of~\cite{HM04},
\[
\frac{1}{n^\gamma}T_n\longrightarrow\TT_{\gamma,\nu} .
\]
When
$\gamma\in(0,1)$, a similar result (Theorem~\ref{sec:main-result-2})
holds when $T_n$ has distribution~$\bQ^q_n$.

The limiting tree $\TT_{\gamma,\nu}$ can be seen as the genealogical
tree of a continuous model for mass splitting, in some sense analogous
to the Markov branching property described above. The above
convergence holds in distribution in a space of measured metric
spaces, endowed with the so-called Gromov--Hausdorff--Prokhorov
topology. This result contrasts with the situation of~\cite{BrDeMcLdlS},
where it is assumed that macroscopic splitting
events occur at every step of the construction. In that case, the
height of $T_n$ is of order $\log n$, and no interesting\vadjust{\goodbreak} scaling limit
exists for the tree. A key step in our study will be to use the
results from~\cite{HaMi09}, where scaling limits of
nonincreasing
Markov chains were considered: such Markov chains are indeed obtained
by considering the successive sizes of collections containing a
particular marked ball when going up in the tree $T_n$.

This general statement allows us to recover, and sometimes improve, many
results of~\cite{HMPW06,HPW,PiWi08,CFW} dealing specifically with
Markov branching trees. It also applies to models of random trees that
are not a priori directly connected to our study. In
particular, we recover the results of Aldous~\cite{AldousCRTIII} and
Duquesne~\cite{duq02} showing that the so-called Brownian and stable
trees~\cite{AldousCRTI,LeGallLeJan,dlg02,dlg05} are universal limits
for conditioned Galton--Watson trees.

More notably, our results entail that uniform unordered trees with $n$
vertices, in which each vertex has at most $m \in[2,\infty]$
children, admit the Brownian continuum random tree as a scaling
limit. This was conjectured by Aldous~\cite{AldousCRTII} and proved in
\cite{MaMi09} in the particular case $m=2$ of a binary branching,
using completely different methods from the present paper. The
difficulty of handling such families of random trees comes from the
fact that they have no ``nice'' probabilistic representations, using,
for instance, branching processes or growth models. As a matter of
fact, uniform random unordered trees do not even have the Markov
branching property, but it turns out to be ``almost'' the case, in a
sense that will be explained below.

The rest of this section is devoted to a detailed formalization of our
results.

\subsection*{Index of notation}

Throughout the paper, we use the notation
\[
\N=\{1,2,3,\ldots\} ,\qquad \Z_+=\{0\}\cup\N,\qquad
[n]=\{1,2,\ldots,n\} ,\qquad n\in\N.
\]
The random variables appearing in this paper are either canonical or
defined on some probability space
$(\Omega,\FF,\PP)$.\vspace*{15pt}

\begin{center}
\begin{tabular}{@{}lp{308pt}@{}}
$\mathbbl{t}$ & plane tree, page \pageref{planetree}\\[2pt]
$\mathsf{t}$& unordered tree, page \pageref{unorderedtree} \\[2pt]
$\bT_n$& set of trees with $n$ vertices, page \pageref{btn} \\[2pt]
$\bT^\partial_n$& set of trees with $n$ leaves, page \pageref{btn}\\[2pt]
$p(\lambda)$& number of parts of a partition $\lambda$,
page \pageref{partsoflambda}      \\[2pt]
$\P_n$& set of partitions of $n$, page \pageref{partitionofn}     \\[2pt]
$m_j(\lambda)$& multiplicity of parts of $\lambda$ equal to $j$, page
\pageref{multiplicity} \\[2pt]
$\bP^q_n$& distributions of Markov branching trees indexed by
leaves, page~\pageref{mbtleaves} \\[2pt]
$\bQ^q_n$& distributions of Markov branching trees indexed by
vertices, page~\pageref{mbtvertices}\\[2pt]
$T_n$ & tree with distribution $\bP^q_n$ or $\bQ^q_n$, page \pageref
{tn}\\[2pt]
$d_{\mathrm{GH}}$& pointed Gromov--Hausdorff distance,
page \pageref{pointedgh}  \\[2pt]
$d_{\mathrm{GHP}}$& pointed Gromov--Hausdorff--Prokhorov distance,
page \pageref{pointedghp} \\[2pt]
$\T$& set of isometry classes of compact
rooted $\R$-trees, page \pageref{isomclasses}
\end{tabular}
\end{center}

\begin{center}
\begin{tabular}{@{}lp{308pt}@{}}
$\Tw$ & set of isometry classes of compact
rooted measured $\R$-trees, pa\-ge~\pageref{isomclasses}\\[2pt]
$\mathcal S^{\downarrow}$ & set of partitions of a unit mass, page
\pageref{partitionmass}\\[2pt]
$\TT_{\gamma,\nu}$ & $(\gamma,\nu)$-fragmentation tree, page
\pageref{fragtree}\\[2pt]
$\bT^{(m)}$ & set of trees with $n$ vertices and at most $m$
children per vertex, pa\-ge~\pageref{treesm}\\[2pt]
$\P_B$ & set of partitions of $B\subseteq\N$, page \pageref
{partitionofB}\\[2pt]
$\P$ & set of partitions with variable size,
page \pageref{partitionsvariablesizes}\\[2pt]
$\theta$ & tree with edge-lengths, page \pageref{treeedgelength}\\[2pt]
$\bTheta$ & set of trees with edge-lengths,
page \pageref{treeswithedgelengths}\\[2pt]
$\TT(\theta)$ & $\mathbb R$-tree associated to $\theta$,
page \pageref{rtreetheta}\\[2pt]
$\TT(\mathsf{t})$ & $\mathbb R$-tree associated to a tree $\mathsf
{t}$ with
edge-lengths 1, page \pageref{rtreelength1}\\[2pt]
$D^\pi_{B}$ & death time of the block $B$ in the process $\pi$, page
\pageref{dpiB}\\[2pt]
$\theta(\pi(\cdot),B)$ &tree with edge-lengths associated with a
partition-valued process, page \pageref{thetapi}\\[2pt]
$p_B(\pi)$ & exchangeable distribution on partitions of $B$
associated with $q_n$, page~\pageref{pbpi}
\end{tabular}%\
\end{center}

%s1.1 ###
\subsection{Discrete trees}\label{sec:discrete-trees}

We briefly introduce some formalism for trees. Set
$\N^0=\{\varnothing\}$, and let
\[
\mathcal{U}=\bigcup_{n\geq0}\N^n .
\]
For $u=(u_1,\ldots,u_n)\in
\mathcal{U}$, we denote by $|u| = n$ the length of $u$, also called the
height of $u$. If $u=(u_1,\ldots,u_n)$ with $n\geq1$, we let
$\mathrm{pr}(u)=(u_1,\ldots,u_{n-1})$, and for $i\geq1$, we let
$ui=(u_1,\ldots,u_n,i)$. More generally, for $u=(u_1,\ldots,u_n)$ and
$v=(v_1,\ldots,v_m)$ in $\mathcal{U}$, we let
$uv=(u_1,\ldots,u_n,v_1,\ldots,v_m)$ be their concatenation. For
$A\subset\mathcal{U}$ and $u\in\mathcal{U}$, we let $uA=\{uv\dvtx v\in
A\}$, and simply let $iA=(i)A$ for $i\in\N$. We say that $u$ is a
prefix of $v$ if $v\in u\mathcal{U}$, and write $u\preceq v$, defining
a partial order on $\mathcal{U}$.

A plane tree is a nonempty, finite subset $\mathbbl{t}\subset\mathcal
{U}$ \label{planetree}
(whose elements are called \textit{vertices}), such that:
\begin{itemize}
\item
if $u\in\mathbbl{t}$ with $|u|\geq1$, then
$\mathrm{pr}(u)\in\mathbbl{t}$;
\item
if $u\in\mathbbl{t}$, then there exists a number $c_u(\mathbbl{t})\in
\Z_+$ (the
number of children of $u$) such that $ui\in\mathbbl{t}$ if and only if
$1\leq
i\leq c_u(\mathbbl{t})$.
\end{itemize}
Let $\partial\mathbbl{t}=\{u\in\mathbbl{t}\dvtx c_u(\mathbbl{t})=0\}$ be
the set of \textit{leaves} of $\mathbbl{t}$.
If $\mathbbl{t}^{(1)},\ldots,\mathbbl{t}^{(k)}$ are plane trees, we
can define a new plane tree by
\[
\bigl\langle\mathbbl{t}^{(1)},\ldots,\mathbbl{t}^{(k)}\bigr\rangle=\{
\varnothing\}\cup\bigcup_{i=1}^k i\mathbbl{t}^{(i)} .\vadjust{\goodbreak}
\]

A plane tree has a natural graphical representation, in which every
$u\in\mathbbl{t}$ is a vertex, joined to its $c_u(\mathbbl{t})$
children by as many
edges. But $\mathbbl{t}$ carries more information than the graph, as it
has a
natural order structure. In this work, we will not be interested in
this order, and we present one way to get rid of this unwanted
structure. Let $\mathbbl{t}$ be a plane tree,
and $\bsigma=(\sigma_u,u\in
\mathbbl{t})$ be a sequence of permutations, respectively,
$\sigma_u\in\mathfrak{S}_{c_u(\mathbbl{t})}$. For $u=(u_1,\ldots
,u_n)\in\mathbbl{t}$,
let
\[
\bsigma(u)=\bigl(\sigma_{\varnothing}(u_1),
\sigma_{(u_1)}(u_2),\sigma_{(u_1,u_2)}(u_3),\ldots,
\sigma_{(u_1,\ldots,u_{n-1})}(u_n)\bigr)
\]
and
$\bsigma(\varnothing)=\varnothing$. Then the set $\bsigma
(\mathbbl{t})=\{\bsigma(u)\dvtx u\in\mathbbl{t}\}$ is a plane tree,
obtained intuitively by
shuffling the set of children of $u$ in $\mathbbl{t}$ according to the
permutation $\sigma_u$. We say that $\mathbbl{t},\mathbbl{t}'$ are
equivalent if there
exists some $\bsigma$ such that $\bsigma(\mathbbl{t})=\mathbbl{t}'$.
Equivalence
classes of plane trees will be called (rooted) \textit{unordered trees},
or simply \textit{trees} as opposed to plane trees, and denoted by
lowercase letter $\mathsf{t}$'s. \label{unorderedtree}
They are sometimes called (rooted) \textit{P\'olya trees} in the
literature~\cite{drmota09}.

Given a tree $\mathsf{t}$, we will freely adapt some notation from plane
trees when dealing with quantities that do not depend on particular
plane representatives. For instance, $\#\mathsf{t},\#\partial\mathsf
{t}$ will denote
the number of vertices and leaves of~$\mathsf{t}$, while $\varnothing,
c_\varnothing(\mathsf{t})$ will denote the root of $\mathsf{t}$ and
its degree.

We let $\bT$ be the set of trees, and for $n\geq1$, \label{btn}
\[
\bT^\partial_n=\{\mathsf{t}\in\bT\dvtx\#\partial\mathsf{t}= n\}
,\qquad
\bT_n=\{\mathsf{t}\in\bT\dvtx\#\mathsf{t}= n\}
\]
be the set of trees with $n$ leaves,
respectively, $n$ vertices. The class of $\{\varnothing\}$ is
the vertex tree $\bullet\in\bT_1=\bT^\partial_1$.

Heuristically, the information carried in a tree is its graph
structure, with a distinguished ``root'' vertex corresponding to
$\varnothing$, and considered up to root-preserving graph
isomorphisms---it is not embedded in any space, and its vertices are unlabeled.
%We will see later that it is sometimes more convenient for our
%purposes to
%re-introduce some structure in these trees by labeling the leaves,
%which will be done in an exchangeable, random way.

It is a simple exercise to see that if $\mathsf{t}^{(i)},1\leq i\leq
k$, are
trees, and $\mathbbl{t}^{(i)}$ is a choice of a plane representative of
$\mathsf{t}^{(i)}$ for each $i$, then the class of $\langle\mathbbl {t}^{(i)},1\leq
i\leq k\rangle$ does not depend on the particular choice for
$\mathbbl{t}^{(i)}$. We denote this common class by $\langle\mathsf
{t}^{(i)},1\leq
i\leq k\rangle$. Note that $j(\mathsf{t}):=\langle\mathsf{t}\rangle
$ can be seen as
the tree $\mathsf{t}$ whose root has been attached to a new root by an edge,
and similarly $j^l (\mathsf{t})$, for $l\geq0$, is the tree $\mathsf
{t}$ whose
root has been attached to a new root by a string of $l$ edges. For
instance, $j^l(\bullet)$ is the line-tree consisting of a string with
length $l$, rooted at one of its ends. Finally, for trees
$\mathsf{t}^{(1)},\ldots,\mathsf{t}^{(k)}$ and $l\geq1$ we let
\[
\bigl\langle\mathsf{t}^{(1)},\ldots,\mathsf{t}^{(k)}\bigr\rangle_l=
j^l\bigl(\bigl\langle\mathsf{t}^{(1)},\ldots,\mathsf{t}^{(k)}\bigr\rangle\bigr) ,
\]
so
$j^l(\bullet)=\langle\bullet\rangle_l$ with this notation.

%%%%%%%%%%%%%%%%%%%%%%%%%%%%%%%%%%%%
%s1.2 ###
\subsection{Markov branching trees}\label{sec:mark-branch-models}
%%%%%%%%%%%%%%%%%%%%%%%%%%%%%%%%%%%%

A partition of an integer $n\geq1$ is a sequence of integers
$\lambda=(\lambda_1,\ldots,\lambda_p)$ with $\lambda_1\geq\cdots
\geq
\lambda_p\geq1$ and $\lambda_1+\cdots+\lambda_p = n$. The number
$p=p(\lambda)$ \label{partsoflambda}
is called the number of parts of the partition
$\lambda$, and the partition is called nontrivial if $p(\lambda)\geq
2$. We let $\P_n$ \label{partitionofn}
be the set of partitions of the integer $n$. We
also add an extra element $\varnothing$ to $\P_1$, so that
$\P_1=\{(1),\varnothing\}$.

If $(c_1,c_2,\ldots)$ is a finite or infinite sequence of nonnegative
integers with finite sum and $j\geq1$, we define
%For $\lambda\in\P_n$ and $1\leq j\leq n$, we define
%
\[
m_j(c_1,c_2,\ldots)=\#\{i\dvtx c_i=j\} ,
\]
the multiplicity of terms of $c_1,c_2,\ldots$ that are equal to
$j$. In particular, if $\lambda\in\P_n$,
$m_j(\lambda)$ \label{multiplicity}
% partition $\lambda$ equal to $j$}
is the multiplicity of parts of
$\lambda$ equal to $j$.

By
convention, it is sometimes convenient to set $\lambda_i=0$ for
$i>p(\lambda)$, and to identify the sequence $\lambda$ with the
infinite sequence $(\lambda_i,i\geq1)$. Such identifications will be
implicit when needed.

%s1.2.1 ###
\subsubsection{Markov branching trees with a prescribed number of
leaves}\label{sec:mark-branch-trees}

In this paragraph, the size of a tree $\mathsf{t}\in\bT$ is going to
be the
number $\#\partial\mathsf{t}$ of its leaves.

Let $q=(q_n,n\geq1)$ be a sequence of probability distributions,
respectively, on~$\mathcal{P}_n$,
\[
q_n=\bigl(q_n(\lambda),\lambda\in\mathcal{P}_n\bigr),\qquad
\sum_{\lambda\in\P_n}q_n(\lambda)=1 ,
\]
such that
%
%e1 ###
\begin{equation}\label{eq:6}
q_n((n))<1 ,\qquad n\geq1 .
\end{equation}
Consider a family of probability
distributions $\bP^q_n,n\geq1$, \label{mbtleaves}
%indexed by leaves}
on $\bT^\partial_n$,
respectively, such that:
\begin{longlist}[(1)]
\item[(1)]
$\bP^q_1$ is the law of the line-tree $\langle\bullet\rangle_G$,
where $G$ has
a geometric distribution given by
\[
\mathbb P(G=k)=q_1(\varnothing)\bigl(1-q_1(\varnothing)\bigr)^{k}
,\qquad
k\geq0 ;
\]
\item[(2)]
for $n\geq2$, $\bP^q_n$ is the law of
\[
\bigl\langle T^{(i)},1\leq i\leq
p(\Lambda)\bigr\rangle,
\]
where
$\Lambda$ has distribution $q_n$, and conditionally on the
latter, the trees $T^{(i)},1\leq i\leq
p(\Lambda)$, are independent with distributions
$\bP^q_{\Lambda_i}$, respectively.
\end{longlist}
Alternatively, for $n\geq2$, $\bP^q_n$ is the law of $\langle
T^{(i)},1\leq i\leq p(\Lambda)
\rangle_G$ where $G$ is independent of $\Lambda$ and geometric with
\[
\mathbb P(G=k)=\bigl(1-q_n((n))\bigr)q_n((n))^k ,\qquad k\geq0 ,
\]
and
conditionally on $\Lambda$, which has law $q_n(\cdot|
\P_n\setminus\{(n)\})$, the trees $T^{(1)}, \ldots,\break T^{(p(\Lambda))}$
are independent with distributions $\bP_{\Lambda_i}$, respectively. A
simple induction argument shows that there exists a unique family
$\bP^q_n, n\geq1$, satisfying properties 1 and 2 above.\vadjust{\goodbreak}

A family of random trees $T_n, n\geq1$, with respective distributions
$\bP^q_n, n \geq1$, is called a Markov branching family. The law of
the tree $T_n$ introduced in the beginning of\vspace*{1pt}
the \hyperref[sec1]{Introduction} to
describe the genealogy of splitting collections of $n$ balls is $
\bP^q_n$.

%s1.2.2 ###
\subsubsection{Markov branching trees with a prescribed number of
vertices}\label{sec:mark-branch-trees-1}

We now consider the following variant of the above construction, in
which the size of a tree $\mathsf{t}$ is the number of its vertices. For
every $n\geq1$, let again $q_n$ be a probability distribution on
$\P_n$. We do not assume~(\ref{eq:6}), rather, we make the sole
assumption that $q_1((1))=1$. For every $n\geq1$, we construct
inductively a family of random trees $T_n$, respectively, in the set
$\bT_n$ of trees with $n$ vertices, by assuming that for
$\lambda=(\lambda_1,\ldots, \lambda_p)\in\P_{n-1}$, with probability
$q_{n-1}(\lambda)$, the $n-1$ vertices distinct from the root vertex
are dispatched in $p$ subtrees with $\lambda_1\geq\cdots\geq
\lambda_p$ vertices, and that, given these sizes, the $p$ subtrees are
independent with same distribution as
$T_{\lambda_1},\ldots,T_{\lambda_p}$, respectively.

Formally:
\begin{longlist}[(1)]
\item[(1)]
let $\bQ^q_1$ be the law of $\bullet$;
\item[(2)]
for $n\geq1$, let $\bQ^q_{n+1}$
%indexed by vertices}
\label{mbtvertices} be the law of
\[
\bigl\langle T^{(i)},1\leq i\leq
p(\Lambda)\bigr\rangle,
\]
where
$\Lambda$ has distribution $q_n$, and conditionally on the
latter, the trees $T^{(i)},1\leq i\leq
p(\Lambda)$, are independent with distributions
$\bQ^q_{\Lambda_i}$, respectively.
\end{longlist}
By induction, these two properties determine the laws $\bQ_n^q,n\geq
1$, uniquely.

The construction is very similar to the previous one, and can in fact
be seen as a special case, after a simple transformation on the tree;
see Section~\ref{proof-theor2} below.

%s1.3 ###
\subsection{Topologies on metric spaces}\label{gh-topology}

The main goal of the present work is to study scaling limits of trees
with distributions $\bP^q_n, \bQ^q_n$, as $n$ becomes large. For this
purpose, we need to consider a topological ``space of trees'' in which
such limits can be taken, and define the limiting objects.

A \textit{rooted}\setcounter{footnote}{1}\footnote{Usually such spaces are rather called
\textit{pointed}, but we prefer the term rooted which is more common
when dealing with trees.} metric space is a triple $(X,d,\rho)$,
where $(X,d)$ is a metric space and $\rho\in X$ is a distinguished
point, called the root. We say that two rooted spaces
$(X,\rho,d),(X^{\prime},\rho^{\prime},d^{\prime})$ are
isometry-equivalent if there exists a bijective isometry from $X$ onto
$X$ that sends $\rho$ to $\rho'$.

A measured, rooted metric space is a 4-tuple $(X,d,\rho,\mu)$, where
$(X,d,\rho)$ is a rooted metric space and $\mu$ is a Borel probability
measure on $X$. Two measured, rooted spaces $(X,d,\rho,\mu)$ and
$(X,d',\rho',\mu')$ are isometry-equivalent if there exists a
root-preserving, bijective isometry $\phi$ from $(X,d,\rho)$ to
$(X,d',\rho')$ such that the push-forward of $\mu$ by $\phi$ is
$\mu'$. In the sequel we will almost\vadjust{\goodbreak} always identify two
isometry-equivalent (rooted, measured) spaces, and will often use the
shorthand notation $X$ for the isometry class of a rooted space or a
measured, rooted space, in a way that should be clear from the
context. Also, if $X$ is such a space and $a>0$, then we denote by
$aX$ the space in which the distance function is multiplied by $a$.

We denote by $\M$ the set of equivalence classes of compact rooted
spaces, and by $\Mw$ the set of equivalence classes of compact
measured rooted spaces.
%These spaces are respectively endowed with the
%pointed Gromov-Hausdorff distance and the pointed
%Gromov-Hausdorff-Prokhorov distance, the definitions of which are
%recalled in Section~\ref{gh-topology}.

It is well known (this is an easy extension of the results of
\cite{EPW}) that $\M$ is a Polish space when endowed with the
so-called rooted Gromov--Hausdorff distance~$d_{\mathrm{GH}}$,
\label{pointedgh}
where by definition the distance
$d_{\mathrm{GH}}((X,d,\rho),(X^{\prime},d^{\prime},\rho^{\prime
}))$ is
equal to the infimum of the quantities
\[
\dist(\phi(\rho),\phi'(\rho'))\vee
\dist_{\mathrm{H}}(\phi(X),\phi'(X')),
\]
where $\phi,\phi'$ are isometries from $X,X^{\prime}$ into a common
metric space $(M,\dist)$, and where $\dist_{\mathrm{H}}$ is the
Hausdorff distance between compact subsets of $(M,\dist)$. It is
elementary that this distance does not depend on particular choices in
the equivalence classes of $(X,d,\rho)$ and
$(X^{\prime},d^\prime,\rho^{\prime})$. We endow $\M$ with the
associated Borel $\sigma$-algebra. Of course, $d_{\mathrm{GH}}$
satisfies a homogeneity property,
$d_{\mathrm{GH}}(aX,aX')=ad_{\mathrm{GH}}(X,X')$ for $a>0$.

We also need to define a distance on $\Mw$, that is in some sense
compatible with the Gromov--Hausdorff distance. Several complete
distances can be constructed, and we use a variation of the
Gromov--Hausdorff--Prokhorov distance used in~\cite{miermont09}. The
induced topology is the same as that introduced earlier in
\cite{EW}. The reader should bear in mind that the topology used in
the present paper involves a little extension of the two previous
references, since we are interested in rooted spaces. We let
$d_{\mathrm{GHP}}((X,d,\rho,\mu),(X',d',\rho',\mu'))$
\label{pointedghp}
be the infimum
of the quantities
\[
\dist(\phi(\rho),\phi'(\rho'))\vee
\dist_{\mathrm{H}}(\phi(X),\phi'(X'))\vee
\dist_{\mathrm{P}}(\phi_*\mu,\phi'_*\mu') ,
\]
where again
$\phi,\phi'$ are isometries from $X,X'$ into a common space
$(M,\dist)$, $\phi_*\mu,\phi'_*\mu'$ are the push-forward of
$\mu,\mu'$ by $\phi,\phi'$ and $\dist_{\mathrm{P}}$ is the Prokhorov
distance between Borel probability measures on $M$ (\cite{EK}, Chapter
3),
\[
\dist_{\mathrm{P}}(m,m')=\inf\{\eps>0\dvtx m(C)\leq
m'(C^\eps)+\eps\mbox{ for every }C\subset M\mbox{ closed}\} ,
\]
where $C^\eps=\{x\in M\dvtx\inf_{y\in C}\dist(x,y)<\eps\}$ is the
$\eps$-thickening of $C$. A simple adaptation of the results of
\cite{EW} and Section 6 in~\cite{miermont09} (in order to take into
account the particular role of the distinguished point $\rho$) shows
the following:

\begin{proposition}\label{sec:r-trees-gromov-1}
The function $d_{\mathrm{GHP}}$ is a distance on $\Mw$ that makes it
complete and separable.
\end{proposition}

This distance is called the rooted Gromov--Hausdorff--Prokhorov distance.
One must be careful that contrary to $d_{\mathrm{GH}}$, this distance
is not homogeneous: $d_{\mathrm{GHP}}(aX,aX')$ is
in general different from $ad_{\mathrm{GHP}}(X,X')$, because only the
distances, not the measures, are multiplied in $aX,aX'$.

%s1.3.1 ###
\subsubsection{Trees viewed as metric spaces}\label{sec:r-trees-gromov}

A plane tree $\mathbbl{t}$ can be naturally seen as a metric space by endowing
$\mathbbl{t}$ with the \textit{graph distance} between vertices. Namely,
\[
d_{\mathrm{gr}}(u,v)=|u|+|v|-2|u\wedge v| ,\qquad u,v\in\mathbbl {t} ,
\]
where $u\wedge v$ is the longest prefix common to $u,v$. This
coincides with the number of edges on the only simple path going from
$u$ to $v$. The space $(\mathbbl{t},d_{\mathrm{gr}})$ is naturally
rooted at
$\varnothing$. We can put two natural probability measures on $\mathbbl{t}$,
the uniform measures on the leaves or on the vertices
\[
\mu_{\partial\mathbbl{t}}=\frac{1}{\#\partial\mathbbl{t}}\sum
_{u\in
\partial\mathbbl{t}}\delta_{\{u\}} ,\qquad \mu_\mathbbl{t}=\frac
{1}{\#\mathbbl{t}}\sum_{u\in
\mathbbl{t}}\delta_{\{u\}} .
\]
If $\mathsf{t}\in\bT$ is a tree, and $\mathbbl{t},\mathbbl{t}'$ are
two plane representatives of $\mathsf{t}$, then it is elementary that the
spaces $(\mathbbl{t},d_{\mathrm{gr}},\varnothing,\mu_{\partial
\mathbbl{t}})$ and
$(\mathbbl{t}',d_{\mathrm{gr}},\varnothing,\mu_{\partial\mathbbl {t}'})$ are
isometry-equivalent rooted measured metric spaces. The same holds with
$\mu_{\mathbbl{t}},\mu_{\mathbbl{t}'}$ instead of $\mu_{\partial
\mathbbl{t}},\mu_
{\partial\mathbbl{t}'}$. We denote by
$(\mathsf{t},d_{\mathrm{gr}},\rho,\mu_{\partial\mathsf{t}})$ and
$(\mathsf{t},d_{\mathrm{gr}},\rho,\mu_\mathsf{t})$ the
corresponding elements of
$\Mw$. Conversely, it is possible to recover uniquely the discrete
tree (not a plane tree!) from the element of $\Mw$ thus defined.

%s1.3.2 ###
\subsubsection{$\R$-trees}\label{sec:r-trees}

An $\mathbb{R}$-tree is a metric space $(X,d)$ such that for every $
x$, \mbox{$y\in X$}:
\begin{longlist}[(1)]
\item[(1)] there is an isometry $\varphi_{x,y}\dvtx[0,d(x,y)]\to X$ such
that $
\varphi_{x,y}(0)=x$ and $\varphi_{x,y}(d(x,y))=y$;
\item[(2)] for every continuous, injective function $c\dvtx[0,1]\to X$ with
$c(0)=x$, \mbox{$c(1)=y$}, one has $c([0,1])=\varphi_{x,y}([0,d(x,y)])$.
\end{longlist}
In other words, any two points in $X$ are linked by a geodesic path,
which is the only simple path linking these points, up to
reparameterisation. This is a continuous analog of the
graph-theoretic definition of a tree as a connected graph with no
cycle. We denote by $\lbr x,y\rbr $ the range of $\varphi_{x,y}$.

We let $\T$ (resp., $\Tw$)
%rooted $\R$-trees}
%rooted measured $\R$-trees}
\label{isomclasses}
be the set of isometry classes of compact
rooted $\R$-trees (resp., compact, rooted measured $\R$-trees). An
important property is the following (these are easy variations on
results by~\cite{EPW,EW}):
\begin{proposition}\label{sec:spaces-metric-spaces}
The spaces $\T$ and $\Tw$ are closed subspaces of
$(\M,d_{\mathrm{GH}})$ and $(\Mw,d_{\mathrm{GHP}})$.
\end{proposition}

If $\TT\in\T$ and for $x\in\TT$, we call $d(\rho,x)$ the \textit{height}
of $x$. If $x,y\in\TT$, we say that $x$ is an ancestor of $y$ whenever
$x\in\lbr\rho,y\rbr$. We let $x\wedge y\in\TT$ be the unique element
of $\TT$ such that $\lbr \rho,x\rbr \cap\lbr\rho,y\rbr =\lbr \rho,x\wedge y\rbr $,
and call it the highest common ancestor of $x$ and $y$ in $\TT$. For
$x\in\TT$, we denote by $\TT_x$ the set of $ y\in\TT$ such that $x$
is an ancestor of $y$. The set $\TT_x$, endowed with the restriction
of the distance $d$, and rooted at $x$, is in turn a rooted
$\mathbb{R}$-tree, called the subtree of $\TT$ rooted\vadjust{\goodbreak} at $x$. If
$(\TT,d,\rho,\mu)$ is an element of $\Tw$ and $\mu(\TT_x)>0$, then
$\TT_x$ can be seen as an element of $\Tw$ by endowing it with the
measure $\mu(\cdot|\TT_x)=\mu(\cdot\cap\TT_x)/\mu(\TT_x)$.

We say that $x\in\TT$, $x\neq\rho$, in a rooted $\mathbb{R}$-tree
is a
\textit{leaf} if its removal does not disconnect $\TT$. Note that this
always excludes the root from the set of leaves, which we denote by
$\mathcal{L}(\TT)$. A \textit{branch point} is an element of $\TT$ of
the form $x\wedge y$ where $x$ is not an ancestor of $y$ nor
vice-versa. It is also characterized by the fact that the removal of a
branch point disconnects the $\mathbb{R}$-tree into three or more
components (two or more for the root, if it is a branch point). We
let $\mathcal{B}(\TT)$ be the set of branch points of $\TT$.

%s1.4 ###
\subsection{Self-similar fragmentations and associated $\R$-trees}
\label{sec:self-simil-fragm}

Self-similar fragmentation processes are continuous-time processes
that describe the dislocation of a massive object as time passes.
Introduce the set of partitions of a unit mass \label{partitionmass}
\[
\sd:=\biggl\{\bs=(s_1,s_2,\ldots)\dvtx s_1\geq s_2\geq\cdots\geq
0,\sum_{i\geq1}s_i\leq1\biggr\}.
\]
This space is endowed with the
metric $d(\bs,\bs')={\sup_{i\geq1}}|s_i-s'_i|$, which makes it
a compact space.
\begin{definition}
\label{sec:self-simil-fragm-4}
A self-similar fragmentation is a $\sd$-valued Markov process
$(\bX(t),t\geq0)$ which is continuous in probability and satisfies
the following fragmentation property. For some $a\in\R$, called the
self-similarity index, it holds that conditionally given
$\bX(t)=(s_1,s_2,\ldots)$, the process $(\bX(t+t'),t'\geq0)$ has same
distribution as the process whose value at time $t'$ is the
decreasing rearrangement of the sequences $s_i\bX^{(i)}(s_i^at'),i\geq
1$, where $(\bX^{(i)},i\geq1)$ are i.i.d. copies of $\bX$.
\end{definition}

Bertoin~\cite{BertoinSSF} and Berestycki~\cite{berest02} have shown that
the laws of self-similar fragmentation processes are characterized by
three parameters: the index $a$, a nonnegative erosion coefficient
and a dislocation measure $\nu$ on $\sd$. The idea is that every
sub-object of the initial object, with mass $x$ say, will suddenly
split into sub-sub-objects of masses $xs_1,xs_2,\ldots$ at rate
$x^a\nu(\d\bs)$, independently of the other sub-objects. Erosion
accounts for the formation of zero-mass particles that are
continuously ripped off the fragments.

For our concerns, we will consider only the special case where
the erosion phenomenon has no role and the dislocation measure does not
charge the set $\{\bs\in
\sd\dvtx\sum_is_i<1\}$. One says that $\nu$ is \textit{conservative}. This
motivates the following definition.
\begin{definition}
\label{disloc}
A \textit{dislocation measure} is a
$\sigma$-finite measure $\nu$ on $\sd$ such that
$\nu(\{(1,0,0,\ldots)\})=0$ and
%
%e2 ###
\begin{equation}
\label{eq:1}
\nu\biggl(\biggl\{\sum_{i\geq1}s_i<1\biggr\}\biggr)=0 ,\qquad
\int_\sd
(1-s_1)\nu(\d\bs)<\infty.
\end{equation}
We say that the measure is \textit{binary} when $\nu(\{s_1+s_2<1\})=0$. A
binary measure is characterized by its image $\nu(s_1\in\d x)$
through the mapping \mbox{$\bs\mapsto s_1$}.

A \textit{fragmentation pair} is a
pair $(a,\nu)$ where $a\in\R$ is called the self-similarity index,
and $\nu$ is a dislocation measure.
\end{definition}

Fragmentation pairs $(a,\nu)$ therefore characterize the distributions
of the self-similar fragmentations we are focusing on. When
$a=-\gamma<0$, small fragments tend to split faster, and it turns out
that they all disappear in finite time, a property known as \textit{formation of dust}. Using this property, it is shown in~\cite{HM04}
how to construct a \textit{fragmentation continuum random tree} encoding the
genealogy of the fragmentation processes. More precisely, a
fragmentation tree is a random element $(\TT,d,\rho,\mu)$ of $\Tw$
(often denoted $\TT$ for simplicity), such that almost surely:
\begin{longlist}[(1)]
\item[(1)]
the measure $\mu$ is supported on the set $\mathcal{ L}(\TT)$ of
leaves of $\TT$;
\item[(2)]
$\mu$ has no atom;
\item[(3)]
for every $x\in\TT\setminus\mathcal{L}(\TT)$, it holds that
$\mu(\TT_x)>0$.
\end{longlist}
Moreover, $\TT$ satisfies the following self-similarity property with
index $-\gamma$. For every $t\geq0$, let $\TT^\circ_i(t),i\geq1$, be
the connected components of the open set $\{x\in\TT\dvtx d(\rho,x)>t\}$,
and let $\TT_i(t)$ be the closure of $\TT^\circ_i(t)$ in $\TT$. It is
plain that $\TT_i(t)\setminus\TT^\circ_i(t)=\{\rho_{i,t}\}$ for some
$\rho_{i,t}\in\TT$, with $d(\rho_{i,t},\rho)=t$. The space
$(\TT_i(t),d,\rho_{i,t},\mu(\cdot|\TT_i(t)))$ is then a random element
in $\Mw$. The self-similarity property then states that for every
$t\geq0$, conditionally given $(\mu(\TT_i(s)),i\geq1)$, $s \leq t$,
the family $\{\TT_i(t),i\geq1\}$ has same distribution as
$\{\mu(\TT_{i}(t))^{\gamma}\TT^{(i)},i\geq1\}$, where $(\TT^{(i)},i
\geq1)$ are i.i.d. copies of $\TT$.

If $\TT$ is a self-similar fragmentation tree with self-similarity
index $-\gamma$, then by~\cite{HM04}, Proposition 1, the process
$((\mu(\TT_{i}(t)),i\geq1)^\da,t\geq0)$ of the nonincreasing
rearrangement of the $\mu$-masses of the trees $\TT_i(t)$, is an
$\sd$-valued self-similar fragmentation process with index
$-\gamma$. The law of this process is thus characterized by a unique
fragmentation pair $(-\gamma,\nu)$. By~\cite{HM04}, Proposition 1, the
law of $\TT$ is entirely characterized by $(-\gamma,\nu)$. In the
sequel, we will let $\TT_{\gamma,\nu}$ \label{fragtree}
be a random variable with this
law. We postpone a more constructive description of this tree to
Section~\ref{sec:cont-fragm-trees}.

It was shown in~\cite{HM04} that one can recover the celebrated
Brownian and stable continuum random trees
\cite{AldousCRTI,LeGallLeJan,dlg02} as special instances of
fragmentation trees. The parameters $\gamma$ and $\nu$ corresponding
to these
trees will be recalled when we discuss applications in Sections
\ref{sec:appl-galt-wats} and~\ref{sec:appl-polya-trees}.

%s1.5 ###
\subsection{Main results}\label{sec:main-result}

Let $(q_n(\lambda),\lambda\in\P_n),n\geq1$, satisfy~(\ref{eq:6}).
With it, we associate a finite nonnegative measure $\ov{q}_n$ on
$\sd$, defined by its integral against measurable functions $f\dvtx\sd\to
\R_+$ as
\[
\ov{q}_n(f)=\sum_{\lambda\in\P_n}q_n(\lambda)
f\biggl(\frac{\lambda}{n}\biggr).
\]
Note that in the left-hand side, we
have identified $\lambda/n$ with an element of $\sd$, in accordance
with our convention that $\lambda$ is identified with the infinite
sequence $(\lambda_i,i\geq1)$. We make the following basic
assumption:

\begin{longlist}[(H)]
\item[(H)] There exists a fragmentation pair $(-\gamma,\nu)$,
with $\gamma>0$, and a function $\ell\dvtx(0,\infty)\to(0,\infty)$ slowly
varying at $\infty$, such that we have the weak convergence of finite
nonnegative measures on $\sd$,
%
%e3 ###
\begin{equation}\label{eq:3}
n^\gamma\ell(n) (1-s_1)\ov{q}_n(\d\bs)
\mathop{\longrightarrow}^{(\mathrm{w})}_{n\to\infty}
(1-s_1)\nu(\d\bs) .
\end{equation}
\end{longlist}

\begin{theorem}\label{sec:main-result-1}
Assume
$q=(q_n(\lambda),\lambda\in\P_n),n\geq1$, satisfies assumption~\textup{(H)}.
Let $T_n$ have distribution $\bP^q_n$, and view $T_n$ as a
random element of $\Mw$ by
endowing it with
the graph distance and the uniform probability measure $\mu_{\partial
T_n}$ on $\partial
T_n$. Then we have the convergence in distribution
\[
\frac{1}{n^\gamma\ell(n)}T_n
\mathop{\longrightarrow}^{(d)}_{n\to\infty}
\TT_{\gamma,\nu}
\]
for the rooted Gromov--Hausdorff--Prokhorov topology.
\end{theorem}

There is a similar statement for the trees with laws
$\bQ^q_n$. Consider a family $(q_n(\lambda),\lambda\in\P_n),n\geq1$,
with $q_1((1))=1$.
\begin{theorem}\label{sec:main-result-2}
Assume $q=(q_n(\lambda),\lambda\in\P_n),n\geq1$, satisfies assumption~\textup{(H)}, with:
\begin{itemize}
\item
either $\gamma\in(0,1)$, or
\item
$\gamma=1$ and $\ell(n)\to0$ as $n\to\infty$.
\end{itemize}
Let $T_n$ have distribution $\bQ^q_n$. We view $T_n$ as a random
element of $\Mw$ by endowing it with the graph distance and the
uniform probability measure $\mu_{T_n}$ on $T_n$. Then we have the
convergence in distribution
\[
\frac{1}{n^\gamma\ell(n)}T_n\mathop{\longrightarrow}^{(d)}_{n\to
\infty}
\TT_{\gamma,\nu}
\]
for the rooted Gromov--Hausdorff--Prokhorov topology.
\end{theorem}

Theorem~\ref{sec:main-result-2} deals with a more restricted set of
values of values of $\gamma$ than Theorem
\ref{sec:main-result-1}. This comes from the fact that, contrary to
the set $\bT^\partial_n$ which contains trees with arbitrary height,
the set $\bT_n$ of trees with $n$ vertices has elements with height at
most $n-1$. Therefore, we cannot hope to find nontrivial limits in
Theorem~\ref{sec:main-result-2} when $\gamma>1$, or when $\gamma=1$
and $\ell(n)$ has limit $+\infty$ as $n\to\infty$. The intermediate
case where $\ell(n)$ admits finite nonzero limiting points cannot
give such a convergence with a continuum fragmentation tree in the
limit either. Indeed, the support of the height of a continuum
fragmentation tree is unbounded, whereas the heights of $T_n/n
\ell(n)$ are all bounded from above by $1/\inf_n(\ell(n))$, which is
finite under our assumption.\vadjust{\goodbreak}

Note that Theorem~\ref{sec:main-result-1} (resp., Theorem
\ref{sec:main-result-2}) implies that any fragmentation tree $\mathcal
T_{\gamma, \nu}$ is the continuous limit of a rescaled family of
discrete Markov branching trees with a prescribed number of leaves
(resp., with a prescribed number of vertices, provided $\gamma<1$),
since we have the following approximation result:\vspace*{-1pt}
\begin{proposition}
\label{Propexemple}
For every fragmentation pair $(-\gamma,\nu)$ with $\gamma>0$, there
exists a family of distributions $(q_n,n \geq1)$ satisfying
(\ref{eq:6}) and such that~(\ref{eq:3}) holds, with $\ell(x)=1$ for
every $x>0$.\vspace*{-1pt}
\end{proposition}

After some preliminaries gathered in Section~\ref{secprelim},
%Sections~\ref{sec:sett-prel} and~\ref{sec:cont-fragm-trees},
we prove Theorems~\ref{sec:main-result-1} and~\ref{sec:main-result-2}
and Proposition~\ref{Propexemple} in Section
\ref{sec:proof-theor-refs}. Before embarking in the proofs, we present
in Section~\ref{applications} some important applications of these
theorems to Galton--Watson trees, unordered random trees and particular
families of Markov branching trees studied in earlier works. Of these
applications, the first two actually involve a substantial amount of
work, so that the details are postponed to Section
\ref{sec:cond-galt-wats} and~\ref{SectionPolya}.

%%%%%%%%%%%%%%%%%%%%%%%%%%
%s2 ###
\section{Applications}\label{applications}
%%%%%%%%%%%%%%%%%%%%%%%%%%

%%%%%%%%%%%%%%%%%%%%%%%%%%%%%%
%s2.1 ###
\subsection{Galton--Watson trees}\label{sec:appl-galt-wats}
%%%%%%%%%%%%%%%%%%%%%%%%%%%%%%

A natural application is the study of Galton--Watson trees
conditioned on their total number of vertices. Let $\xi$ be a
probability measure on $\Z_+$ such that $\xi(0)>0$ and
%
%e4 ###
\begin{equation}\label{eq:9}
\sum_{k\geq0}k\xi(k)=1 .
\end{equation}
The law of the Galton--Watson tree with offspring distribution
$\xi$ is the probability measure on the set of plane trees defined by
\[
\mathrm{GW}_\xi(\{\mathbbl{t}\})=\prod_{u\in\mathbbl{t}}\xi
(c_u(\mathbbl{t})) ,
\]
for $\mathbbl{t}$ a
plane tree. That this does define a probability distribution on the
set of plane trees comes from the fact that a Galton--Watson process
with offspring distribution $\xi$ becomes a.s. extinct in finite
time, due to the criticality condition~(\ref{eq:9}). In order to fit
in the framework of this paper, we view $\mathrm{GW}_\xi$ as a
distribution on the set of discrete, rooted trees, by taking its
push-forward under the natural projection from plane trees to trees.

In order to avoid technicalities, we also assume that the support of
$\xi$ generates the additive group $\Z$. This implies that
$\mathrm{GW}_\xi(\{\#\mathbbl{t}= n\})>0$ for every $n$ large enough.
For such
$n$, we let $\gw^{(n)}=\gw(\cdot| \{\#\mathbbl{t}= n\})$, and view
it as a law on $\bT_n$.

We distinguish two different regimes.

\textit{Case} 1. The offspring distribution has finite variance
\[
\sigma^2=\sum_{k\geq0}k(k-1)\xi(k)<\infty.\vadjust{\goodbreak}
\]

\textit{Case} 2. For some $\alpha\in(1,2)$ and $c\in(0,\infty)$, it
holds that $\xi(k)\sim ck^{-\alpha-1}$ as $k\to\infty$. In
particular, $\xi$ is in the domain of attraction of a stable law of
index $\alpha$.

The \textit{Brownian dislocation measure} is the unique binary
dislocation measure such that
\[
\nu_{2}(s_1\in\d x)=\sqrt{\frac{2}{\pi x^3(1-x)^3}} \,\d x
\ind_{\{1/2\leq x< 1\}} .
\]
Otherwise said,
for every measurable $f\dvtx\sd\to\R_+$,
\[
\int_\sd\nu_2(\d\bs) f(\bs)
=\int_{1/2}^1\sqrt{\frac{2}{\pi
x^3(1-x)^3}} \,\d x f(x,1-x,0,0,\ldots) .
\]
We also define a
one-parameter family of measures in the following way. For
$\alpha\in(1,2)$, let
%Let $(T_x,x\geq0)$ be a stable subordinator with
%index $1/\alpha$, i.e. a nondecreasing L\'evy process with Laplace
%$$E[\exp(-t T_x)]=\exp(-x t^{1/\alpha}) , t,x>0 .$$ We let
%$\Delta T_{[0,1]}$ be the decreasing rearrangement of the sequence of
%jumps of $(T_x,0\leq x\leq1)$, so that $\Delta T_{[0,1]}/T_1$ is an
%element of $\sd$ with total sum $1$ (as $T$ is a.s. a pure-jump
%process). Said in another, more elementary way (which is equivalent by
%the L\'{e}vy-It\'\'{} representation of L\'{e}vy processes), the
%sequence $\Delta
%T_{[0,1]}$ is the nonincreasing rearrangement of the sequence
$\sum_{i\geq1}\delta_{\Delta_i}$ be a Poisson random measure on
$(0,\infty)$ with intensity measure
\[
\frac{1}{\alpha\Gamma(1-{1/\alpha})}\,\frac{\d
x}{x^{1+1/\alpha}}\ind_{\{x>0\}}
\]
with the atoms
$\Delta_i,i\geq1$, labeled in such a way that $\Delta_1\geq
\Delta_2\geq\cdots.$ Let $T=\sum_{i\geq1}\Delta_i$, which is finite
a.s. by standard properties of Poisson measures. In fact, $T$ follows
a stable distribution with index $1/\alpha$, with Laplace transform
\[
\mathbb{E}[\exp(-\lambda T)]=\exp(-\lambda^{1/\alpha}) ,\qquad
\lambda\geq0 .
\]
This can be seen as a stable subordinator evaluated
at time $1$, its jumps up to this time being the atoms $\Delta_i,i\geq
1$. The measure $\nu_\alpha$ is defined by its action against a
measurable function $f\dvtx\sd\to\R_+$
\[
\int_\sd\nu_\alpha(\d\bs)
f(\bs)=\frac{\alpha^2\Gamma(2-1/\alpha)}{ \Gamma(2-\alpha)}
\mathbb E\biggl[T
f\biggl(\frac{\Delta_i}{T},i\geq1\biggr)\biggr] .
\]
Because
$\mathbb E[T]=\infty$, this formula defines an infinite $\sigma
$-finite measure on
$\sd$, which turns out to satisfy~(\ref{eq:1}).
\begin{theorem}\label{sec:main-results-1}
Let $\xi$ satisfy~(\ref{eq:9}), with support that generates the
additive group~$\Z$. Let $T_n$ be a random element of $\bT_n$ with
distribution\vspace*{1pt} $\gw^{(n)}$. Consider $T_n$ as an element of $\Mw$ by
endowing it with the graph distance and the uniform probability
measure $\mu_{T_n}$ on $T_n$. Then we have, in distribution for the
Gromov--Hausdorff--Prokhorov topology:
\begin{eqnarray*}
& &\mbox{Case 1:\quad}\hspace*{015pt}
\frac{1}{\sqrt{n}}T_n
\mathop{\longrightarrow}^{(d)}_{n\to\infty}\frac{2}{\sigma}\TT
_{1/2,\nu_2} ;\\
& &\mbox{Case 2:\quad}
\frac{1}{n^{1-1/\alpha}}T_n
\mathop{\longrightarrow}^{(d)}_{n\to\infty}
\biggl(\frac{\alpha(\alpha-1)}{c\Gamma(2-\alpha)}
\biggr)^{1/\alpha}\TT_{1-1/\alpha,\nu_\alpha} .
\end{eqnarray*}
\end{theorem}

This result will be proved in Section~\ref{sec:cond-galt-wats} below,
by first showing that $\gw^{(n)}$ is of the form $\bQ^q_n$ for some
appropriate choice of $q$.

The trees $\TT_{1/2,\nu_2}$ and $\TT_{1-1/\alpha,\nu_\alpha}$
appearing in the limit are important models of continuum random trees,
called, respectively, the Brownian Continuum Random Tree and the stable
tree with index $\alpha$. The Brownian tree is somehow the archetype
in the theory of scaling limits of trees. The above theorem is very
similar to a result due to Duquesne~\cite{duq02}, but our method of
proof is totally different. While~\cite{duq02} relies on quite refined
aspects of Galton--Watson trees and their encodings by stochastic
processes, our approach requires only to have some kind of global
structure, namely the Markov branching property, and to know how mass
is distributed in one generation. We do not claim that our method is more
powerful than the one used in~\cite{duq02} (as a matter of fact, the
limit theorem of~\cite{duq02} holds in the more general case where
$\mu$ is in the domain of attraction of a totally asymmetric stable
law with index $\alpha\in(1,2]$). However, our method has some
robustness, allowing us to shift from Galton--Watson trees to other
models of trees. Our next example will try to illustrate this.

%%%%%%%%%%%%%%%%%%%%%%%%%%%%%%%%%%%
%s2.2 ###
\subsection{Uniform unordered trees}\label{sec:appl-polya-trees}
%%%%%%%%%%%%%%%%%%%%%%%%%%%%%%%%%%%

Our next application is on a different model of random trees, which is
by nature not a model of plane or labeled trees, contrary to the
previous examples. It is actually not either a Markov branching model,
but is very close from being one, as we will see.

For $2\leq m\leq\infty$, we consider the set $\bT^{(m)}_n\subset
\bT_n$
%children per vertex}
\label{treesm} of trees with $n$ vertices, in which every vertex has
at most
$m$ children. In particular, we have $\bT^{(\infty)}_n=\bT_n$. The
sets $\bT_n^{(m)}$ are harder to enumerate than sets of ordered or
labeled trees, like plane trees or Cayley trees, and there is no
closed expression for the numbers $\#\bT^{(m)}_n$. However, Otter
\cite{otter48} (see also~\cite{FlSe09}, Section VII.5) derived the
asymptotic enumeration result
%
%e5 ###
\begin{equation}\label{eq:18}
\#\bT^{(m)}_n
\mathop{\sim}_{n\to\infty} \kappa_m\frac{(\rho_m)^n}{n^{3/2}}
\end{equation}
for some $m$-dependent constants $\kappa_m>0,\rho_m>1$. This can be
achieved by studying the generating function
\[
\psi^{(m)}(x)=\sum_{n\geq1} \#\bT^{(m)}_n x^n ,
\]
which has a
square-root singularity at the point $1/\rho_m$. The behavior
(\ref{eq:18}) indicates that a uniformly chosen element of
$\bT^{(m)}_n$ should converge as $n\to\infty$, once renormalized
suitably, to the Brownian continuum random tree. We show that this is
indeed the case for any value of $m$. To state our result, let
\[
\wt{\bT}_n^{(m)}=\bigl\{\mathsf{t}\in\bT_n^{(m)}\dvtx c_\varnothing(\mathsf
{t})\leq
m-2\bigr\}.
\]
For instance, $\wt{\bT}_n^{(2)}=\varnothing$ for $n \geq
2$, while $\wt{\bT}^{(\infty)}_n=\bT^{(\infty)}_n$ for all $n$. Let
\[
\wt{\psi}^{(m)}(x)=\sum_{n\geq1}\#\wt{\bT}^{(m)}_n x^n ,
\]
and
define a finite constant $c_m$ by
\[
c_m=\frac{\sqrt{2}}{\sqrt{\pi}\kappa_m\wt{\psi}^{(m)} (1/\rho
_m)}
.
\]
Note that $\wt{\psi}^{(2)}(x)=x$ for every $x$, while
$\wt{\psi}^{(\infty)}(1/\rho_\infty)=1$
(\cite{FlSe09}, Section~VII.5). Therefore, we get
\[
c_2=\frac{\sqrt{2}\rho_2}{\sqrt{\pi}\kappa_2} ,\qquad
c_\infty=\frac{\sqrt{2}}{\sqrt{\pi}\kappa_\infty} .
\]

\begin{theorem}\label{sec:main-results-2}
Fix $m\in\{2,3,\ldots\}\cup\{\infty\}$. Let $T_n$ be uniformly
distributed in~$\bT^{(m)}_n$. We view $T_n$ as an element of $\Mw$ by
endowing it with the measure $\mu_{T_n}$, then
\[
\frac{1}{\sqrt{n}}T_n
\mathop{\longrightarrow}^{(d)}_{n\to\infty}c_m\TT_{1/2,\nu_2}
\]
for the Gromov--Hausdorff--Prokhorov topology.
\end{theorem}

The proof of this result is given in Section~\ref{SectionPolya}. We
note that this implies a similar, maybe more natural, statement for
$m$-ary trees. We say that $\mathsf{t}\in\bT$ is $m$-ary if every
vertex has
either $m$ children or no child, and we say that the vertex is
internal in the first case, that is, when it is not a leaf. Summing over
the degrees of vertices in an $m$-ary tree with $n$ internal vertices,
we obtain that $\#\mathsf{t}=mn+1$ and $\#\partial\mathsf{t}=(m-1)n+1$.

Assume now that $m<\infty$. Starting from a $m$-ary tree $\mathsf{t}$ with
$n$ internal vertices, and removing the leaves---equivalently,
keeping only the internal vertices---gives an element
$\phi(\mathsf{t})\in\bT^{(m)}_n$. The mapping $\phi$ is inverted
by attaching
$m-k$ leaves to each vertex with $k$ children, for an element of
$\bT^{(m)}_n$.
%, except for the root, to which is attached
%$m-k$ leaves when it is of degree $k$
Moreover, we leave as an easy
exercise that $d_{\mathrm{GHP}}(a\mathsf{t},a\phi(\mathsf{t}))\leq
a$ for every
$a>0$, when the trees are endowed with the uniform measures
$\mu_\mathsf{t},\mu_{\phi(\mathsf{t})}$ on vertices. Theorem
\ref{sec:main-results-2} thus implies the following:
\begin{corollary}\label{sec:unif-unord-trees}
Let $m\in\{2,3,\ldots\}$ and $T_n^{[m]}$ be a uniform $m$-ary tree
with $n$ internal vertices, endowed with the measure
$\mu_{T_n^{[m]}}$. Then
\[
\frac{1}{\sqrt{n}}T^{[m]}_n\mathop{\longrightarrow}^{(d)}_{n\to
\infty}
c_m\TT_{1/2,\nu_2} .
\]
\end{corollary}

The problem of scaling limits of random rooted unordered trees has
attracted some attention in the recent literature; see
\cite{BrFl08,MaMi09,DrGi10,drmota09}. For $m=2$,
Corollary~\ref{sec:unif-unord-trees}\vadjust{\goodbreak} readily yields the main theorem of
\cite{MaMi09}, which was derived using a completely different method,
based in a stronger way on combinatorial aspects of
$\bT^{(2)}_n$. Here, we really make use of a fragmentation property
satisfied by the uniform distributions on $\bT_n^{(m)}, n \geq1$. As
alluded to at the beginning of this section, these are not actually
laws of Markov branching trees. Nevertheless, they can be coupled with
laws of Markov branching trees in a way that the coupled trees are
close in the $d_{\mathrm{GHP}}$ metric. In the general case $m\neq2$,
Theorem~\ref{sec:main-results-2} and Corollary
\ref{sec:unif-unord-trees} are new, and were implicitly conjectured by
Aldous~\cite{AldousCRTII}. In~\cite{DrGi10}, the authors prove a
result on the scaling limit of the so-called profile of the uniform
tree for $m=\infty$, which is related to our results, although it is
not a direct consequence. Finally, we note that the problem of the
scaling limit of \textit{unrooted} unordered trees is still open,
although we expect the Brownian tree to arise again as the limiting
object.

%%%%%%%%%%%%%%%%%%%%%%%%%%%%%%%%%%%%%%%%
%s2.3 ###
\subsection{Consistent Markov branching models}\label{sec:appl-cons-models}
%%%%%%%%%%%%%%%%%%%%%%%%%%%%%%%%%%%%%%%%

Considering again in a more specific way the Markov branching models,
we stress that Theorem~\ref{sec:main-result-1} also encompasses the
results of~\cite{HMPW06}, which hold for particular families
$(q_n,n\geq1)$ satisfying a further consistency property. In this
setting, it is assumed that $q_n((n))=0$ for every $n\geq1$, so that
the trees $T_n$ do not have any vertex having only one child. The
consistency property can be formulated as follows:

\begin{itemize}
\item[]\textit{Consistency property}. Starting from $T_n$ with $n\geq2$, select
one of the leaves uniformly at random, and remove this leaf as well as
the edge that is attached to it. If this removal creates a vertex with
only one child, then remove this vertex and merge the two edges
incident to this vertex into one. Then the random tree thus
constructed has same distribution as~$T_{n-1}$.
\end{itemize}

A complete characterization of families $(q_n,n\geq1)$ giving rise to
Markov branching trees with this consistency property is given in
\cite{HMPW06}. Namely, such families can be constructed in terms of a
pair $(c,\nu)$, which is uniquely defined up to multiplication by a
common positive constant, such that $c\geq0$ is an ``erosion
coefficient'' and $\nu$ is a dislocation measure as described above
[except that $\nu(\sum s_i<1)=0$ is not required]. The cases where
$c=0$ and $\nu(\sum s_i<1)=0$ are the most interesting ones, so we
will assume henceforth that this is the case. The associated
distributions $q_n,n\geq2$, are given by the following explicit
formula: for $\lambda\in\P_n$ having $p\geq2$ parts,
%
%e6 ###
\begin{equation}\label{eq:7}
q_n(\lambda)=\frac{1}{Z_n}C_\lambda\int_{\sd}\nu(\d\bs)
\mathop{\sum_{i_1,\ldots,i_p\geq1}}_{\mathrm{distinct}}
\prod_{j=1}^ps_{i_j}^{\lambda_j} ,
\end{equation}
where
\[
C_\lambda=\frac{n!}{\prod_{i\geq1}\lambda_i !\prod_{j\geq
1}m_j(\lambda)!}
\]
is a combinatorial factor, the same that
appears in the statement of Lemma~\ref{sec:partitions-1} below, and
$Z_n$ is a normalizing constant defined by
\[
Z_n=\int_{\sd}\nu(\d\bs)\biggl(1-\sum_{i\geq1}s_i^n\biggr) .
\]
Assume further that $\nu$ satisfies the following regularity
condition:
%
%e7 ###
\begin{equation}\label{eq:8}
\nu(s_1\leq1-\eps)=\eps^{-\gamma}\ell(1/\eps) ,
\end{equation}
where $\gamma\in(0,1)$ and $\ell$ is a function that is slowly
varying at $\infty$.
\begin{theorem}\label{sec:main-results}
If $\nu$ is a dislocation measure satisfying~(\ref{eq:1}) and
(\ref{eq:8}), and if $(q_n,n\geq1)$ is the consistent family of
probability measures defined by~(\ref{eq:7}), then the Markov
branching trees $T_n$, viewed as random measured $\R$-trees by
endowing the sets of their leaves with the uniform probability
measures, satisfies
\[
\frac{1}{\Gamma(1-\gamma)n^\gamma
\ell(n)}T_n\mathop{\longrightarrow}^{(d)}_{n\to\infty}
\TT_{\gamma,\nu} ,
\]
for the Gromov--Hausdorff--Prokhorov topology.
\end{theorem}

This theorem is in some sense more powerful than~\cite{HMPW06},
Theorem 2, because the latter result needed one extra technical
hypothesis that is discarded here. Moreover, our result holds for the
Gromov--Hausdorff--Prokhorov topology, which is stronger than the
Gromov--Hausdorff topology considered in~\cite{HMPW06}. However, the
setting of~\cite{HMPW06} also provided a natural coupling of the trees
$T_n,n\geq1$, and $\TT_{\gamma,\nu}$ on the same probability
space, for which the convergence in Theorem~\ref{sec:main-results} can
be strengthened to a convergence in probability. This coupling is not
provided in our case.
\begin{pf*}{Proof of Theorem~\ref{sec:main-results}}
Let $\mathbf s \in\mathcal S^{\downarrow}$ be
such that $\sum_{i \geq1} s_i=1$. Let $K_1,\ldots, K_n$ be
i.i.d. random variables in $\N$ such that $\PP(K_1=i)=s_i$ for every \mbox{$i
\geq1$}. Call $ \Lambda^{(i)}(n)$ the number of variables $K_j$ equal to
$i$, and let $\Lambda^{(\mathbf
s)}(n)=(\Lambda^{(1)}(n),\Lambda^{(2)}(n),\ldots
)^{\downarrow}$,
where $\mathbf{x}^\downarrow$ denotes the decreasing rearrangement of
the nonnegative sequence $\mathbf{x}=(x_1,x_2,\ldots)$ with finite
sum.
It is not hard to
see that the probability distributions $q_n$ defined by~(\ref{eq:7})
are also given, for
$\lambda\neq(n)$, by
\[
q_n(\lambda)=\frac{1}{Z_n}\int_{\mathcal S^{\downarrow}} \mathbb
P\bigl(\Lambda^{(\mathbf s)}(n)= \lambda\bigr) \nu(\mathrm d \mathbf s ).
\]
See, for example, the forthcoming Lemma~\ref{sec:partitions-1} in
Section~\ref{sec:discrete-model}.
The normalizing constant $Z_n$ is regularly varying, according to
the assumption of regular variation~(\ref{eq:8}). Indeed, by
Karamata's Tauberian theorem (see~\cite{BGT},\vadjust{\goodbreak} Theorem~1.7.1'), we have that
\begin{eqnarray*}
Z_n&=&\int_{\mathcal S^{\downarrow}}\biggl(1-\sum_{i \geq
1}s_i^n\biggr)\nu(\mathrm d \mathbf s)
\mathop{\sim}_{n\to\infty}
\int_{\mathcal S^{\downarrow}}(1-s_1^n)\nu(\mathrm d
\mathbf s) \\
&\displaystyle \mathop{\sim}_{n\to\infty}&
\Gamma(1-\gamma)\nu(s_1 \leq
1-1/n)\\
&=&\Gamma(1-\gamma)n^{\gamma}\ell(n).
\end{eqnarray*}
Now, to get a convergence
of the form~(\ref{eq:3}), note that for all continuous functions
$f\dvtx\mathcal S^{\downarrow} \rightarrow\mathbb R_+$,
%(by Fubini's theorem since $f$ is bounded on the compact space
%$\mathcal S^{\downarrow}$)
%
\begin{eqnarray*}
Z_n \sum_{\lambda\in\mathcal P_n} q_n(\lambda) \biggl(1-\frac
{\lambda_1}
{n}\biggr) f\biggl( \frac{\lambda}{n}\biggr) &=& \int_{\mathcal
S^{\downarrow}} \nu(\mathrm d \mathbf s) \mathbb E\biggl[
\biggl(1-\frac{\Lambda^{(\mathbf s)}_1(n)}{n}
\biggr)f\biggl(\frac{\Lambda^{(\mathbf s)}(n)}{n} \biggr) \biggr]
\\
&\displaystyle \mathop{\rightarrow}_{n\to\infty}& \int_{\mathcal
S^{\downarrow}}\nu(\mathrm d \mathbf s) (1-s_1)f(\mathbf s) ,
\end{eqnarray*}
which follows by dominated convergence, since $f$ is bounded (say by
$K$) on the compact space $\mathcal S^{\downarrow}$ and
\begin{eqnarray*}
\mathbb
E\biggl[ \biggl(1-\frac{\Lambda^{(\mathbf s)}_1(n)}{n} \biggr)
f\biggl(\frac{\Lambda^{(\mathbf s)}(n)} {n} \biggr) \biggr] &\leq&
K\mathbb E\biggl[1-\frac{\Lambda^{(\mathbf s)}_1(n)}{n} \biggr] \leq
K\mathbb E\biggl[1-\frac{\Lambda^{(1)}(n)}{n} \biggr] \\
&=&K(1-s_1) .
\end{eqnarray*}
We conclude by applying Theorem~\ref{sec:main-result-1}.
\end{pf*}

%%%%%%%%%%%%%%%%%%%%%%%%%%%%%%%%%%%%%%%%%%%%%%%
%s2.4 ###
\subsection{\texorpdfstring{Further nonconsistent cases: $(\alpha,\theta)$-trees}
{Further nonconsistent cases: (alpha, theta)-trees}}\label{sec:some-natural-non}
%%%%%%%%%%%%%%%%%%%%%%%%%%%%%%%%%%%%%%%%%%%%%%%

One final application concerns a family of binary labeled trees
introduced by Pitman and Winkel~\cite{PiWi08} and built inductively
according to a growth rule depending on two parameters $ \alpha\in
(0,1)$ and $\theta\geq0$. Roughly, at each step, given that the tree
$T_n^{\alpha,\theta, \mathrm{lab}}$ with $n$ leaves branches at the
branch point adjacent to the root into two subtrees with $k\geq1$
leaves for the subtree containing the smallest label in
$T_n^{\alpha,\theta, \mathrm{lab}}$ and $n-k \geq1$ leaves for the
other one, a weight $ \alpha$ is assigned to the root edge and weights
$k-\alpha$ and $n-k-1+\theta$ are assigned, respectively, to the trees
with sizes $k$, $n-k$. Then choose either the root edge or one of the
two subtrees according with probabilities proportional to these
weights. If a subtree with two or more leaves is selected, apply this
weighting procedure inductively to this subtree until the root edge or
a subtree with a single leaf is selected. If a subtree with single
leaf is selected, insert a new edge and leaf at the unique edge of
this subtree. Similarly, if the root edge is selected, add a new edge
and leaf to this root edge. This gives the tree
$T_{n+1}^{\alpha,\theta, \mathrm{lab}}$. We then denote by
$T_n^{\alpha,\theta}$ the tree $T_n^{\alpha,\theta, \mathrm{lab}}$
without labels, $n \geq1$.

Pitman and Winkel show that the family $(T_n^{\alpha,\theta}, n\geq
1)$ is not consistent in general (\cite{PiWi08}, Proposition 1),
except when $\theta=1-\alpha$ or $\theta=2- \alpha$, and has the\vadjust{\goodbreak}
Markov branching property (\cite{PiWi08}, Proposition 11) with the
following probabilities $q_n$:
\begin{enumerate}
\item[$\bullet$] $q_n((k,n-k,0,\ldots) ) =
q_{\alpha,\theta}(n-1,k)+q_{\alpha,\theta}(n-1,n-k), $ for
$n-k<k \leq n-1$;
\item[$\bullet$] $ q_n (n/2,n/2 )=
q_{\alpha,\theta}(n-1,n/2) $,
\end{enumerate}
where
\[
q_{\alpha,\theta}(n,k)= \pmatrix{n\cr k} \frac{ \alpha(n-k)
+\theta k }{n}
\frac{\Gamma(k-\alpha)\Gamma(n-k+\theta)}{\Gamma(1-\alpha)\Gamma
(n+\theta)},\qquad
1 \leq k \leq n.
\]
Now consider the binary measure $\nu_{\alpha,\theta}$ defined on
$\mathcal S^{\downarrow}$ by $\nu_{\alpha,\theta}(s_1+s_2<1)=0$ and
$\nu_{\alpha,\theta}(s_1 \in
\mathrm dx)=f_{\alpha,\theta}(x) \,\mathrm dx$ where $f_
{\alpha,\theta}$ is defined on $[1/2,1)$ by
\begin{eqnarray*}
f_{\alpha,\theta}(x)&=&\frac{1}{\Gamma(1-\alpha)}\bigl( \bigl(
\alpha
(1-x ) + \theta x \bigr) x^{-
\alpha-1}(1-x)^{\theta-1}\\
&&\hphantom{\frac{1}{\Gamma(1-\alpha)}\bigl(}
{} + \bigl( \alpha x + \theta(1-x)
\bigr) (1-x)^{-\alpha-1}x^{\theta-1}\bigr).
\end{eqnarray*}

\begin{theorem}
Endow $T_n^{\alpha,\theta}$ with the uniform probability measure on
$\partial T_n^{\alpha,\theta}$. Then,
\[
\frac{1}{n^{\alpha}}T_n^{\alpha,\theta} \mathop{\longrightarrow
}^{(d)}_{n\to\infty}
\mathcal
T_{\alpha,\nu_{\alpha,\theta}}
\]
for the rooted
Gromov--Hausdorff--Prokhorov topology.
\end{theorem}

This result reinforces Proposition 2 of~\cite{PiWi08} which states
the a.s. convergence of~$T_n^{\alpha,\theta} $, in a certain
finite-dimensional sense, to a continuum fragmentation tree with
parameters $\alpha,
\nu_{\alpha,\theta}$.
In view of Theorem~\ref{sec:main-result-1},
it suffices to check that hypothesis~(H) holds, which in the
present case states that
for any $f\dvtx\sd\to\R$ continuous with $|f(\bs)|\leq(1-s_1)$,
\begin{eqnarray*}
&&n^{\alpha}\sum_{k=\lceil n/2 \rceil}^{n-1} f
\biggl(\frac{k}{n},\frac{n-k}{n},0,\ldots\biggr) q_n
\bigl((k,n-k,0,\ldots)
\bigr) \\
&&\qquad\rightarrow\int_{1/2}^1 f(x,1-x,0,\ldots)
f_{\alpha,\theta}(x)\,\mathrm dx .
\end{eqnarray*}
To prove this, we use that
$\int_0^1 x^{a-1}(1-x)^{b-1} \,\mathrm dx = \Gamma(a) \Gamma
(b)/\Gamma(a+b)$ and rewrite $q_{\alpha,\theta}(n-1,k)$ as
\begin{eqnarray*}
q_{\alpha,\theta}(n-1,k)&=&\pmatrix{n-1\cr k} \frac{\alpha
(n-1-k)
+ \theta k}{n-1}
\frac{\Gamma(n-1+\theta-\alpha)}{\Gamma(1-\alpha)\Gamma
(n-1+\theta)} \\
&&{}\times\int_0^1 x^{k-\alpha-1}(1-x)^{n-k+\theta-2} \,\mathrm
dx.
\end{eqnarray*}
Then set for $x \in[0,1]$,
\[
F(x):=f(x, 1-x,0,\ldots) \mathbf
1_{\{x > 1/2\}} +f(1-x, x,0,\ldots) \mathbf1_{\{x \leq1/2\}},\vadjust{\goodbreak}
\]
and
note that $F(0)=0$ and $|F(x) |\leq(1-x) \wedge x$, for every $x \in
[0,1]$. We have
\begin{eqnarray*}
&&\sum_{k=\lceil n/2 \rceil}^{n-1} f
\biggl(\frac{k}{n},\frac{n-k}{n},0,\ldots\biggr) q_n
\bigl((k,n-k,0,\ldots)
\bigr)\\[-2pt]
&&\qquad= \sum_{k=0}^{n-1} F \biggl(\frac{k}{n}\biggr)
q_{\alpha,\theta}(n-1,k) \\[-2pt]
&&\qquad=
\frac{\Gamma(n-1+\theta-\alpha)}{\Gamma(1-\alpha)\Gamma
(n-1+\theta)}\\[-2pt]
&&\qquad\quad{}\times
\int_0^1 \sum_ {k=0}^{n-1} \pmatrix{n-1\cr k} \frac{ \alpha
(n-1-k)
+ \theta k}{n-1} F \biggl(\frac{k}{n} \biggr)
x^{k-\alpha-1}(1-x)^{n-k+\theta-2} \,\mathrm dx \\[-2pt]
&&\qquad=
\frac{\Gamma(n-1+\theta-\alpha)}{\Gamma(1-\alpha)\Gamma
(n-1+\theta)}\\[-2pt]
&&\qquad\quad{}\times
\int_0^1 \mathbb E \biggl[ \biggl( \alpha
\biggl(1-\frac{B_{n-1}^{(x)}}{n-1} \biggr)+ \theta\frac{ B_{n-1}^{(x)}}{n-1}
\biggr) F \biggl(\frac{B_{n-1}^{(x)}}{n} \biggr)\biggr]
x^{-\alpha-1}(1-x)^{\theta-1} \,\mathrm dx,
\end{eqnarray*}
where $B_{n-1}^{(x)}$ denotes a binomial random variable with
parameters $(n-1, x)$. We can assume that $B_{n-1}^{(x)}/n \rightarrow
x$ a.s. on the probability space $(\Omega,\FF,\PP)$, and since $F$ is
continuous and bounded on $[0,1]$, we have
\begin{eqnarray*}
&&\mathbb E \biggl[ \biggl( \alpha\biggl(1-\frac{B_{n-1}^{(x)}}{n-1}
\biggr)+\theta\frac{B_{n-1}^{(x)}}{n-1} \biggr) F
\biggl(\frac{B_{n-1}^{(x)}}{n} \biggr)\biggr]\\[-2pt]
&&\qquad \rightarrow\bigl(
\alpha
(1-x ) + \theta x \bigr) F(x) \qquad\mbox{for every
}x \in[0,1].
\end{eqnarray*}
Moreover,
\begin{eqnarray*}
\label{eqexp}
&&\mathbb E \biggl[ \biggl( \alpha\biggl(1-\frac{B_{n-1}^{(x)}}{n-1}
\biggr)+\theta\frac{B_{n-1}^{(x)}}{n-1} \biggr) F
\biggl(\frac{B_{n-1}^{(x)}}{n} \biggr)\biggr]\\[-2pt]
&&\qquad\leq\biggl((\alpha+ \theta
) \mathbb E \biggl[ \frac{B_{n-1}^{(x)}}{n} \biggr] \biggr)
\wedge\biggl(\alpha\mathbb E \biggl[
1-\frac{B_{n-1}^{(x)}}{n-1}\biggr] + \theta\mathbb E \biggl[\frac
{B_{n-1}^{(x)}}{n-1} \biggr] \biggr) \\[-2pt]
&&\qquad\leq\bigl((\alpha+
\theta)x \bigr) \wedge\bigl( \alpha(1-x) + \theta x\bigr).
\end{eqnarray*}
This is enough to conclude by dominated convergence
that
\begin{eqnarray*}
&&\int_0^1 \mathbb E \biggl[ \biggl(
\alpha
\biggl(1-\frac{B_{n-1}^{(x)}}{n-1} \biggr)+\theta\frac{B_{n-1}^{(x)}}
{n-1} \biggr) F \biggl(\frac{B_{n-1}^{(x)}}{n} \biggr)\biggr]
x^{-\alpha-1}(1-x)^{\theta-1} \d x\\[-2pt]
&&\qquad
\mathop{\longrightarrow}_{n\to\infty}
\int_0^1 \bigl( \alpha(1-x ) +
\theta x \bigr) F (x)x^{-\alpha-1}(1-x)^{\theta-1} \,\mathrm dx\\[-2pt]
&&\qquad=
\Gamma(1-\alpha) \int_{1/2}^1 f(x,1-x,\ldots) f_{\alpha,
\theta}(x)\,\mathrm dx.
\end{eqnarray*}
Last, Stirling's formula implies that
\[
\frac{\Gamma(n-1+\theta-\alpha)}{\Gamma(n-1+\theta)}
\mathop{\sim}_{n\to\infty}
n^{- \alpha}
\]
as wanted.\vspace*{-2pt}

%%%%%%%%%%%%%%%%%%%%%%%%%%%%%%%%%%%%%%%%%%%%%%
%s3 ###
\section{Preliminaries on self-similar fragmentations and trees}\label{secprelim}\vspace*{-2pt}
%%%%%%%%%%%%%%%%%%%%%%%%%%%%%%%%%%%%%%%%%%%%%

%%%%%%%%%%%%%%%%%%%%%%%%%%%%%%%%%%%%%%
%s3.1 ###
\subsection{Partition-valued self-similar fragmentations}
\label{sec:sett-prel}
%%%%%%%%%%%%%%%%%%%%%%%%%%%%%%%%%%%%%%

In this section, we recall the aspects of the theory of self-similar
fragmentations that will be needed to prove Theorems
\ref{sec:main-result-1} and~\ref{sec:main-result-2}. We refer the
reader to~\cite{BertoinBook} for more details.\looseness=-1\vspace*{-2pt}

%s3.1.1 ###
\subsubsection{Partitions of sets of integers}\label{sec:partitions}

Let $B\subset\N$ be a possibly infinite, nonempty subset of the
integers, and $\pi=\{\pi_1,\pi_2,\ldots\}$ be a partition of $B$. The
(nonempty) sets $\pi_1,\pi_2,\ldots$ are called the \textit{blocks} of
$\pi$, we denote their number by~$b(\pi)$. In order to remove the
ambiguity in the labeling of the blocks, we will use, unless otherwise
specified, the convention that $\pi_i,i\geq1$, is defined inductively
as follows: $\pi_i$ is the block of $\pi$ that contains the least
integer of the set
\[
B\Bigm\backslash\bigcup_{j=1}^{i-1}\pi_j ,
\]
if the latter is not empty.
For $i\in B$, we also let $\pi_{(i)}$ be the block of $\pi$ that
contains~$i$.

We let $\P_B$ \label{partitionofB}
be the set of partitions of $B$. This forms a partially
ordered set, where we let $\pi\preceq\pi'$ if the blocks of $\pi'$
are all included in blocks of $\pi$ (we also say that $\pi'$ is finer
than $\pi$). The minimal element is $\mathbbm{O}_B=\{B\}$, and the
maximal element is $\mathbb{I}_B=\{\{i\}\dvtx i\in B\}$.

If $B'\subseteq B$ is nonempty, the restriction of $\pi$ to $B'$,
denoted by $\pi|_{B'}$ or $B'\cap\pi$ with a slight abuse of
notation, is the element of $\P_{B'}$ whose blocks are the nonempty
elements of $\{B'\cap\pi_1,B'\cap\pi_2,\ldots\}$.

If $B\subset\N$ is finite, with say $n$ elements, then any partition
$\pi\in\P_B$ with $b$ blocks induces an element $\lambda(\pi)\in
\P_n$ with $b$ parts, given by the nonincreasing rearrangement of the
sequence $(\#\pi_1,\ldots,\#\pi_b)$.

A subset $B\subset\N$ is said to admit an asymptotic frequency if the
limit
\[
\lim_{n\to\infty}\frac{\#(B\cap[n])}{n}
\]
exists. It is
then denoted by $|B|$. It is a well-known fact, due to Kingman, that
if $\pi$ is a random partition of $\N$ with distribution invariant
under the action of permutations (simply called exchangeable
partition), then a.s. every block of $\pi$ admits an asymptotic
frequency.
%The law of $\pi$ is then given by the paintbox construction
%of next section, for some probability measure $\nu$.
We then let
$|\pi|^\da\in\sd$ be the nonincreasing rearrangement of the sequence
$(|\pi_i|,i\geq1)$. The exchangeable partition $\pi$ is called \textit{proper} if $\sum_{i=1}^{b(\pi)}|\pi_i|=1$, which is equivalent to
the fact that $\pi$ has a.s. no singleton blocks.\vadjust{\goodbreak}

%s3.1.2 ###
\subsubsection{Paintbox construction}\label{sec:paintb-constr}

Let $\nu$ be a dislocation measure, as defined in Definition \ref
{disloc}. We
construct a $\sigma$-finite measure on $\P_\N$ by the so-called \textit{paintbox construction}. Namely, for every $\bs\in\sd$ with
$\sum_{i\geq1}s_i=1$, consider an i.i.d. sequence $(K_i,i\geq1)$
such that
\[
\mathbb P(K_1=k)=s_k ,\qquad k\geq1 .
\]
Then the partition
$\pi$ such that $i,j$ are in the same block of $\pi$ if and only if
$K_i=K_j$ is exchangeable. We denote by $\rho_\bs(\d\pi)$ its
law. Note that $\rho_\bs(\d\pi)$-a.s., it holds that $|\pi|^\da
=\bs$,
and $\pi$ is a.s. proper under $\rho_\bs$. The measure
\[
\kappa_\nu(\d\pi):=\int_{\sd}\nu(\mathrm d \bs)\rho_\bs(\d
\pi)
\]
is a
$\sigma$-finite measure on $\P_\N$, invariant under the action of
permutations. From the integrability condition~(\ref{eq:1}) on $\nu$,
it is easy to check that for $k\geq2$,~if
\[
A_k=\bigl\{\pi\in\P_\N\dvtx\pi|_{[k]}\neq\{[k]\}\bigr\}
\]
is the set of partitions whose trace on $[k]$ has at least two blocks,
then
%
%e8 ###
\begin{equation}
\label{eq:2}
\kappa_\nu(A_k)=\int_{\sd}\nu(\d\bs)\biggl(1-\sum_{i\geq
1}s_i^k\biggr)<\infty
\end{equation}
for every $k\geq2$, since $1-\sum_{i\geq1}s_i^k\leq1-s_1^k\leq
k(1-s_1)$.

%s3.1.3 ###
\subsubsection{Exchangeable partitions of finite and infinite
sets}\label{sec:exch-part-finite}

In this section, we establish some elementary results concerning
exchangeable partitions of $[n]$ or $\N$. The set of partitions with
variable size, namely \label{partitionsvariablesizes}
\[
\P=\P_\N\cup\bigcup_{n\geq1}\P_{[n]}
\]
is endowed with the
distance
\[
d_\P(\pi,\pi')=\exp\bigl(-\sup\bigl\{k\geq1\dvtx\pi|_{[k]}=\pi'|_{[k]}\bigr\}\bigr) .
\]
In the sequel, convergence in distribution for partitions will be
understood with respect to the separable and complete space
$(\P,d_\P)$. We will use the falling factorial notation
\[
(x)_n=x(x-1)\cdots(x-n+1)=\frac{\Gamma(x+1)}{\Gamma(x-n+1)}
\]
for
$x$ a real number and $n \in\N$, $n <x+1$. When $x \in\mathbb N$,
we extend the notation to all $n \in\mathbb N$, by setting $(x)_n=0$
for $n \geq x+1$.
\begin{lemma}\label{sec:exch-part-finite-1}
Let $\pi$ be an exchangeable partition of $[n]$, and let $k\leq
n$. Then for every $B\subset[k]$ with $l$ elements such that $1\in
B$,
\[
\mathbb P\bigl([k]\cap\pi_{(1)}=B |
\#\pi_{(1)} \bigr)=\frac{(\#\pi_{(1)}-1)_{l-1}(n-\#\pi
_{(1)})_{k-l}}{(n-1)_{k-1}}
.\vadjust{\goodbreak}
\]
\end{lemma}
\begin{pf}
By exchangeability, the probability under
consideration depends on $B$ only through its cardinality, and
this equal to $\mathbb P(i_2,\ldots,i_l\in
\pi_{(1)},j_1,\ldots,\break j_{k-l}\notin\pi_{(1)}|\#\pi_{(1)})$
for any
pairwise disjoint $i_2,\ldots,i_l,j_1,\ldots,j_{k-l}\in
\{2,3,\ldots,n\}$ [note that there are
${n-1\choose l-1}{n-l\choose k-l}$ such choices]. Consequently,
\begin{eqnarray*}
&&
\mathbb P\bigl([k]\cap\pi_{(1)}=B |
\#\pi_{(1)}\bigr)\\
&&\qquad=\frac{\mathbb E[\sum_{i_2,\ldots
,i_l,j_1,\ldots,j_{k-l}}
\ind_{\{i_2,\ldots,i_l\in
\pi_{(1)}\}}\ind_{\{j_1,\ldots,j_{k-l}\notin\pi_{(1)}\}}
|
\#\pi_{(1)}]}{{n-1\choose l-1}{n-l\choose k-l}}\\
&&\qquad=\frac{{\#
\pi_{(1)}-1\choose
l-1}{n-\#\pi_{(1)}\choose k-l}}{{n-1\choose l-1}{n-l\choose
k-l}},
\end{eqnarray*}
where the sum in the expectation is over indices considered
above. This yields the result.
\end{pf}
\begin{lemma}\label{sec:tightness-1}
Let $(\pi^{(n)},n\geq1)$ be a sequence of random exchangeable
partitions, respectively, in $\P_{[n]}$. We assume that $\pi^{(n)}$
converges in distribution to $\pi$. Then $\pi$ is exchangeable and
\[
\frac{\#\pi^{(n)}_{(i)}}{n}
\mathop{\longrightarrow}^{(d)}_{n\to\infty}\bigl|\pi_{(i)}\bigr| ,
\]
the
latter convergences holding jointly for $i\geq1$, and jointly with
the convergence $\pi^{(n)}\to\pi$.
\end{lemma}
\begin{pf}
The fact that $\pi$ is invariant under the action of permutations of
$\N$ with finite support is inherited from the exchangeability of
$\pi^{(n)}$, and one concludes that $\pi$ is exchangeable
\cite{aldous85}.

The random variables $(\#\pi^{(n)}_{(i)}/n,i\geq1)$ take values in
$[0,1]$, so their joint distribution is tight, and up to extraction,
we may assume that they converge in distribution to a random vector
$(x_{(i)},i\geq1)$, jointly with the convergence
$\pi^{(n)}\to\pi$. We want to show that a.s. $x_{(i)}=|\pi_{(i)}|$,
which will characterize the limiting distribution.
%By the Skorokhod representation theorem, we
%may assume that (still along some extraction), we $\pi^n\to\pi$ and
%$\#\pi^n_{(i)}/n\to x_{(i)},i\geq1$ almost-surely as $n\to\infty$.

%Let $\FF$ be the $\sigma$-algebra generated by the random variables
%$\#\pi^n_{(i)},n,i\geq1$ and $|\pi_{(i)}|,i\geq1$.
For $k\geq l\geq1$ fixed, by summing the formula of Lemma
\ref{sec:exch-part-finite-1} over all $B\subset[k]$ containing $i$,
with $l$ elements, we get
\begin{eqnarray*}
&&
\mathbb P\bigl(\#\bigl([k]\cap\pi_{(i)}^{(n)}\bigr)=l | \#\pi
^{(n)}_{(i)}\bigr)\\
&&\qquad=
\pmatrix{k-1\cr l-1}\frac{(\#\pi^{(n)}_{(i)}-1)_{l-1}(n-\#\pi
^{(n)}_{(i)})_{k-l}}{
(n-1)_{k-1}}\\
&&\qquad\mathop{\longrightarrow}_{n\to\infty}
\pmatrix{k-1\cr l-1}x_{(i)}^{l-1}\bigl(1-x_{(i)}\bigr)^{k-l} ,
\end{eqnarray*}
which entails that, conditionally on $x_{(i)}$, $\#([k]\cap
\pi_{(i)})-1$ follows a binomial distribution with parameters
$(k-1,x_{(i)})$. Therefore,
\[
\bigl|\pi_{(i)}\bigr|=\lim_{k\to\infty}\frac{\#([k]\cap
\pi_{(i)})}{k}=x_{(i)} \qquad\mbox{a.s.}
\]
by the law of large numbers.
\end{pf}
\begin{lemma}\label{sec:self-simil-fragm-5}
Let $(\pi^{(i)},1\leq i\leq r)$ be a sequence of random elements of
$\P_\N$, which is exchangeable in the sense that
$(\sigma\pi^{(i)},1\leq i\leq r)$ has the same distribution as
$(\pi^{(i)},1\leq i\leq r)$, for every permutation $\sigma$ of $\N$.
Then for
every $k\geq2$,
\[
\mathbb P\bigl(2,3,\ldots,k\in\pi^{(1)}_{(1)} | \bigl|\pi
^{(i)}_{(j)}\bigr|,1\leq
i\leq r,j\geq1\bigr)=\bigl|\pi_{(1)}^{(1)}\bigr|^{k-1} .
\]
\end{lemma}
\begin{pf}
Let $n\geq k$, and set
$\pi^{(i,n)}=\pi^{(i)}|_{[n]}$, so that $(\pi^{(i,n)},1\leq i\leq r)$
is a random sequence of $\P_{[n]}$ that is exchangeable. Then, by
using the same argument as in the proof of Lemma
\ref{sec:exch-part-finite-1}, it holds that
\[
\mathbb P\bigl(2,3,\ldots,k\in\pi_{(1)}^{(1,n)} |
\#\pi^{(i,n)}_{(j)},1\leq i\leq r,1\leq j\leq n\bigr)=\frac{
(\#\pi^{(1,n)}_{(1)}-1)_{k-1}}{(n-1)_{k-1}} .
\]
Using\vspace*{1pt} Lemma
\ref{sec:tightness-1}, and the fact that $(\pi^{(i,n)},1\leq i\leq r)$
converges in distribution to $(\pi^{(i)},1\leq i\leq r)$ as $n\to
\infty$, it is then elementary to get the result by taking limits.
\end{pf}

%s3.1.4 ###
\subsubsection{Poisson construction of homogeneous
fragmentations}\label{sec:poiss-constr-homog}

We now recall a useful construction of homogeneous fragmentations
using Poisson point processes. We again fix a dislocation measure
$\nu$.

Consider a Poisson random measure $\mathcal{N}(\d t\,\d\pi\,\d i)$ on the
set $\R_+\times\P_\N\times\N$, with intensity measure $\d
t\otimes\kappa_\nu(\d\pi)\otimes\#_\N(\d i)$, where $\#_\N$ is the
counting measure on $\N$. We use a Poisson process notation
$(\pi^0_t,i^0_t)_{t\geq0}$ for the atoms of $\mathcal{N}$: for
$t\geq
0$, if $(t,\pi,i)$ is an atom of $\mathcal{N}$, then we let
$(\pi^0_t,i^0_t)=(\pi,i)$, and if there is no atom of $\mathcal{N}$ of
the form $(t,\pi,i)$, then we set $\pi^0_t=\mathbbm{O}_\N$ and
$i^0_t=0$ by convention. One constructs a process $(\Pi^0(t),t\geq0)$
by letting $\Pi^0(0)=\mathbbm{O}_\N$, and given that $\Pi
^0(s),0\leq
s<t$, has been defined, we let $\Pi^0(t)$ be the element of $\P_\N$
obtained\vspace*{1pt} from
$\Pi^0(t-)$ by leaving its blocks unchanged, except the $i^0_t$th
block $\Pi^0_{i^0_t}(t-)$,
which\vspace*{-1pt} is intersected with $\pi^0_t$.
%such that for $i\geq1$
%$$\Pi^0_i(t)=\{\begin{array}{ll}\Pi^0_i(t-) & \mbox{ if }i\neq
%i^0_t\\
%
%$$
Of course,\vspace*{1pt} this construction is only informal, since the times $t$ of
occurrence of an atom of $\mathcal{N}$ are everywhere dense in
$\R_+$. However, using~(\ref{eq:2}), it is possible to perform a
similar construction for partitions restricted to $[k]$, and check
that these constructions are consistent as $k$ varies (\cite{BertoinBook},
Section 3.1.1). The process $(\Pi^0(t),t\geq0)$ is called a
partition-valued homogeneous fragmentation with dislocation measure
$\kappa_\nu$.

Note in particular that
the block $\Pi_{(1)}^0(t)$ that contains $1$ at time $t$, is given by
%
%e9 ###
\begin{equation}\label{eq:23}
\Pi_{(1)}^0(t)=\mathop{\bigcap_{0<s\leq t}}_{i^0_s=1} (\pi
^0_s)_{(1)} ,
\end{equation}
and that the restriction of $\mathcal{N}$ to $\R_+\times\P_\N
\times
\{1\}$ is a Poisson measure with intensity $\d t\otimes
\kappa_\nu(\d\pi)$.

For $k\geq2$, let $D_k^0=\inf\{t\geq0\dvtx\Pi^0(t)\in A_k\}$ be the
first time when the restriction of $\Pi^0(t)$ to $[k]$ has at least
two blocks. By the previous construction, it is immediate to see that
$D_k^0$ has an exponential distribution with parameter
$\kappa_\nu(A_k)$: it is the first time $t$ such that $i^0_t=1$ and
$\pi^0_t\in A_k$. Moreover, by standard properties of Poisson random
measures, conditionally on $D_k^0=s$, the random variables $\pi^0_s$
and $(\pi^0_t,i^0_t)_{0\leq t<s}$
are independent, and the law of $\pi^0_s$ equals
$\kappa_\nu(\cdot| A_k)=\kappa_\nu(\cdot\cap
A_k)/\kappa_\nu(A_k)$, while $(\pi^0_t,i^0_t)_{0\leq t<s}$ has the same
distribution as the initial process conditioned on $\{(\pi
^0_t,i_t^0)\notin
A_k \times\{1\},0\leq t<s\}=\{D^0_k\geq s\}$, which has probability
$e^{-s\kappa_\nu(A_k)}$. It is also equivalent to condition on
$\{D^0_k>s\}$, since $\mathbb P(D^0_k=s)=0$. The next statement sums up this
discussion. By definition, we let $X(t\wedge
s-)=X(t)\ind_{\{t<s\}}+X(s-)\ind_{\{t\geq s\}}$ for $X$
c\`adl\`ag.
\begin{lemma}
\label{sec:continuous-model-1}
Let $F,f$ be nonnegative measurable functions. Then
\begin{eqnarray*}
&&\mathbb E\bigl[F\bigl(\Pi^0(t\wedge D^0_k-),t\geq0\bigr)f(\pi
^0_{D_k^0})\bigr]\\
&&\qquad=\kappa_\nu(f|A_k) \int_0^\infty
\kappa_\nu(A_k)\,\d s \mathbb E\bigl[F\bigl(\Pi^0(t\wedge
s),t\geq0\bigr) \ind_{\{ D_k^0> s\}}\bigr] .
\end{eqnarray*}
Otherwise\vspace*{-1pt} said, $\pi^0_{D^0_k}$ and $(\Pi^0(t\wedge
D^0_k-),t\geq0)$, are independent with respective laws
$\kappa_\nu(\cdot| A_k)$, and the law\vspace*{1pt} of $(\Pi^0(t\wedge
\mathbbm{e}),t\geq0)$ conditioned not to split $[k]$, where $\mathbbm
{e}$ is an exponential random variable, independent of $\Pi^0$, and
with parameter $\kappa_\nu(A_k)$.
\end{lemma}

%s3.1.5 ###
\subsubsection{Self-similar fragmentations}\label{sec:self-simil-fragm-1}

From a homogeneous fragmentation $\Pi^0$ constructed as above, one can
associate a one-parameter family of $\P_\N$-valued processes by a
time-changing method. Let $a\in\R$. For every $i\geq1$ we let
$(\tau_{(i)}^a(t),t\geq0)$ be defined as the right-continuous inverse
of the nondecreasing process
\[
\int_0^t\bigl|\Pi^0_{(i)}(u)\bigr|^{-a}\,\d u,\qquad t\geq0 .
\]
For
$t\geq0$, let $\Pi(t)$ be the random partition of $\N$ whose blocks
are given by $\Pi_{(i)}^0(\tau_{(i)}^a(t)),i\geq1$. One\vspace*{-1pt} can check
that this definition is consistent, namely, that for every $j\in
\Pi_{(i)}^0(\tau_{(i)}^a(t))$, one has
$\Pi_{(i)}^0(\tau_{(i)}^a(t))=\Pi_{(j)}^0(\tau_{(j)}^a(t))$.\vadjust{\goodbreak}%\vspace*{2pt}

The process $(\Pi(t),t\geq0)$ is called the self-similar
fragmentation with index $a$ and dislocation measure $\nu$
(\cite{BertoinBook}, Chapter 3.3). We now assume that $a=-\gamma<0$ is
fixed once and for all. Let $D_k=\inf\{t\geq0\dvtx\Pi(t)\in A_k\}$.
\begin{proposition}\label{sec:self-simil-fragm-2}
Conditionally given $\sigma\{\Pi_{(i)}(t\wedge D_k)\dvtx t\geq0,1\leq
i\leq k\}$ and letting $\pi=\Pi(D_k)$, the random\vspace*{2pt} variable
$(\Pi_i(t+D_k),t\geq0)_{1\leq i\leq b([k]\cap\pi)}$ has the same
distribution as $(\pi_i\cap\Pi^{(i)}(|\pi_i|^at),t\geq0)_{1\leq
i\leq b([k]\cap\pi)}$, where $(\Pi^{(i)},i\geq1)$ are
i.i.d. copies of $\Pi$.
\end{proposition}
\begin{pf}
For every $i\geq1$ we let $L_i=\inf\{t\geq
0\dvtx\Pi_{(i)}(t)\cap[k]\neq[k]\}$. Then $L=(L_i,i\geq1)$ is a
so-called \textit{stopping line}, that is, for every $i\geq1$, $L_i$ is a
stopping time with respect to the natural filtration of $\Pi_{(i)}$,
while $L_i=L_j$ for every $j\in\Pi_{(i)}(L_i)$. We let $\Pi(L)$ be
the partition whose blocks are $\Pi_{(i)}(L_i),i\geq1$---by
definition of a stopping line, two such blocks are either equal or
disjoint. Note that $t+L=(t+L_i,i\geq1)$ is also a stopping line, as
well as $t\wedge L=(t\wedge L_i,i\geq1)$.

From the so-called \textit{extended branching property} (\cite{BertoinBook},
Lemma 3.14), we obtain that conditionally given
$\sigma\{\Pi(t\wedge L),t\geq0\}$, the process $(\Pi(t+L),t\geq0)$
has same distribution as
\[
\bigl(\bigl\{\pi_i\cap\Pi^{(i)}(|\pi_i|^at),i\geq1\bigr\},t\geq0\bigr) ,
\]
where
$\pi=\Pi(L)$ and $(\Pi^{(i)},i\geq1)$ are i.i.d. copies of
$\Pi$. The result is then a specialization of this, when looking only
at the blocks of $\Pi$ that contain $1,2,\ldots,k$.
\end{pf}

%The process $(\Pi^a(t),t\geq0)$ is an exchangeable, c\`adl\`ag,
%$\P_\N$-valued process such that $\Pi(0)=\mathbbm{O}_\N$, and such
%that a.s., for every $t\geq0$, the blocks of $\Pi^a(t)$ admit an
%asymptotic frequency. Moreover, for every $t\geq0$, conditionally
%given $\Pi^a(t)=\pi$, $(\Pi^a(t+s),s\geq0)$ has the same distribution
%as the process $(\Pi'(s),s\geq0)$, where $\Pi'(s)$ is the partition
%of $\P_\N$ whose blocks are given by
%$$\pi_i\cap\widetilde{\Pi}^{(i)}(|\pi_i|^as) , i\geq1 .$$
%$\widetilde{\Pi}^{(i)},i\geq1$ being independent copies of
%$(\Pi^a(t),t\geq0)$.

It will be of key importance to characterize the joint distribution of
$D_k$, $(\Pi_{(i)}(D_k),1\leq i\leq k)$. This can be obtained as
a consequence of Lemma~\ref{sec:continuous-model-1}. Recall the
construction of $\Pi$ from $\Pi^0$, let
$\tau_{(i)}=\tau_{(i)}^a$ and define
$\pi_t=\pi^0_{\tau_{(1)}(t)}$. The\vspace*{-1pt} latter is equal to
$\pi^0_{\tau_{(i)}(t)}$ for every $i\in[k]$ and $t\leq D_k$.
\begin{proposition}\label{sec:poiss-constr-homog-1}
Let $F,f$ be nonnegative, measurable functions. Then
\begin{eqnarray*}
&&\mathbb E\bigl[F\bigl(\bigl|\Pi_{(1)}(t\wedge D_k-)\bigr|,t\geq0\bigr)f(\pi
_{D_k})\bigr]\\[-1.5pt]
&&\qquad=\kappa_\nu(f|A_k)
\int_0^\infty\d u
\kappa_\nu(A_k)\mathbb E \bigl[\bigl|\Pi_{(1)}(u)\bigr|^{k-1+a} \mathbf1_{\{
|\Pi_{(1)}(u)|>0\}}\\[-1.5pt]
&&\hspace*{158.5pt}{}\times F\bigl(\bigl|\Pi_{(1)}(t\wedge
u)\bigr|,t\geq0\bigr)\bigr].
\end{eqnarray*}
\end{proposition}
\begin{pf}
By definition, $D_k$ (resp., $D^0_k$) is
the first time when $[k]\cap\Pi(t)\neq[k]$ (resp., $[k]\cap
\Pi^0(t)\neq[k]$). It follows that $D^0_k=\tau_{(1)}(D_k)$,
and that the process
\[
\Pi_{(1)}(t\wedge D_k-)=\Pi^0_{(1)}\bigl(\tau_{(1)}(t\wedge
D_k-)\bigr)=\Pi^0_{(1)}\bigl(\tau_{(1)}(t)\wedge D^0_k-\bigr),\qquad
t\geq0,
\]
is
measurable with respect to $\sigma\{\Pi^0_{(1)}(t\wedge D^0_k-),t\geq
0\}$. Lemma~\ref{sec:continuous-model-1} implies that
\begin{eqnarray*}
&&\mathbb E\bigl[F\bigl(\bigl|\Pi_{(1)}(t\wedge D_k-)\bigr|,t\geq
0\bigr)f(\pi_{D_k})\bigr]
\\
&&\qquad=\mathbb E\bigl[F\bigl(\bigl|\Pi_{(1)}^0\bigl(\tau_{(1)}(t)\wedge
D^0_k-\bigr)\bigr|,t\geq
0\bigr)f(\pi^0_{D^0_k})\bigr] \\
&&\qquad=\kappa_\nu(f|A_k) \int_0^\infty\d
s \kappa_\nu(A_k)\mathbb E\bigl[F\bigl(\bigl|\Pi_{(1)}^{0}\bigl(\tau
_{(1)}(t)\wedge
s\bigr)\bigr|,t\geq0\bigr)\ind_{\{
D_k^0>s\}}\bigr]\\
&&\qquad=\kappa_\nu(f|A_k)\mathbb E\biggl[\int_0^\infty
\d
u
\kappa_\nu(A_k)\bigl|\Pi_{(1)}(u)\bigr|^a\\
&&\qquad\quad\hspace*{72.7pt}{}\times F\bigl(\bigl|\Pi_{(1)}^0\bigl(\tau_{(1)}(t)\wedge
\tau_{(1)}(u)\bigr)\bigr|,t\geq0\bigr)\ind_{\{
D_k>u\}}\biggr]\\
&&\qquad=\kappa_\nu(f|A_k)\mathbb E\biggl[\int_0^\infty
\d
u \kappa_\nu(A_k)\bigl|\Pi_{(1)}(u)\bigr|^aF\bigl(\bigl|\Pi_{(1)}(t\wedge u)\bigr|,t\geq
0\bigr)\ind_{\{ D_k>u\}}\biggr] ,
\end{eqnarray*}
where in the third equality, we used successively Fubini's theorem
and the change of variables $s=\tau_{(1)}(u)$, so that $\d
s=|\Pi_{(1)}(u)|^{a}\,\d u$. We conclude by using the fact that
%
%e10 ###
\begin{equation}\label{eq:24}
\mathbb P\bigl(D_k>u | \bigl|\Pi_{(1)}(t)\bigr|,0\leq t\leq u\bigr)=
\bigl|\Pi_{(1)}(u)\bigr|^{k-1} ,
\end{equation}
which can be argued as follows. Let $0\leq
t_1<t_2<\cdots<t_r=u$ be fixed times, then by applying Lemma
\ref{sec:self-simil-fragm-5} to the sequence $(\Pi(t_i),1\leq i\leq
r)$, we obtain that
\[
\mathbb P\bigl(D_k>u | \bigl|\Pi_{(1)}(t_i)\bigr|,1\leq i\leq r\bigr)=
\bigl|\Pi_{(1)}(u)\bigr|^{k-1} .
\]
This yields~(\ref{eq:24}) by a monotone
class argument, using the fact that $\sigma\{|\Pi_{(1)}(t)|,\break0\leq
t\leq u\}$ is generated by finite cylinder events.
\end{pf}

The last important property of self-similar fragmentations is
that the process $(|\Pi_{(1)}(t)|,t\geq0)$ is a Markov process, which
can be described as follows~\cite{BertoinBook}. Let $(\xi_t,t\geq0)$
be a
subordinator with Laplace transform
\[
\mathbb E[\exp(-r\xi_t)]=\exp\biggl(-t\int_0^\infty
\biggl(1-\sum_{i\geq
1}s_i^{r+1}\biggr)\nu(\d\bs)\biggr) .
\]
Then $(|\Pi^0_{(1)}(t)|,t\geq0)$ has the same distribution as
$(\exp(-\xi_t),t\geq0)$, and consequently, the process $(|\Pi
_{(1)}(t)|,t\geq
0)$ is a so-called self-similar Markov process:
\begin{proposition}[(Corollary 3.1 of~\cite{BertoinBook})]
\label{sec:self-simil-fragm-3}
The process $(|\Pi_{(1)}(t)|,t\geq
0)$ has same distribution as $\exp(-\xi_{\tau(t)},t\geq0)$, where
$\tau$ is the right-continuous inverse of the process
$(\int_0^u\exp(a\xi_s)\,\d s,u\geq0)$.
\end{proposition}

%Conditionally on $[k]\cap\Pi(D_k)=\pi'$, a partition with say $b$
%blocks, and on $(\Pi_{i}(D_{[k]}),1\leq i\leq b)=(\pi_1,\ldots,\pi_b)$,
%it also holds that the processes $(\pi_i\cap\Pi(D_{[k]}+t),t\geq0)$
%for
%$1\leq i\leq b$ are independent, with the same distribution as the
%images of $(\Pi(|\pi_i|^{-\gamma}t),t\geq0)$ under (any) bijections
%$\N\to\pi_i,1\leq i\leq b$.

%%%%%%%%%%%%%%%%%%%%%%%%%%%%%%%%%%%
%s3.2 ###
\subsection{Continuum fragmentation trees}\label{sec:cont-fragm-trees}
%%%%%%%%%%%%%%%%%%%%%%%%%%%%%%%%%%%

This section is devoted to a more detailed description of the limiting
self-similar fragmentation tree $\TT_{\gamma,\nu}$~\cite{HM04}. In
particular, we will need a new decomposition result of reduced trees
at the first branchpoint (Proposition~\ref{sec:cont-fragm-trees-2}).

%s3.2.1 ###
\subsubsection{Trees with edge-lengths and $\R$-trees}
\label{tree-edge-lengths}

We saw in Section~\ref{sec:r-trees-gromov} how to turn a tree into a
(finite) measured metric space. It is also easy to ``turn discrete
trees into $\R$-trees,'' viewing the edges as real segments of length
$1$.

More generally, a plane tree \textit{with edge-lengths} is a pair
$\theta=(\mathbbl{t},(\ell_u,u\in\mathbbl{t}))$
\label{treeedgelength}
where $\ell_u\geq0$ for every $u\in
\mathbbl{t}$, and a \textit{tree with edge-lengths} is obtained by ``forgetting
the ordering'' in a way that is adapted from the discussion of Section
\ref{sec:discrete-trees} in a straightforward way. Namely, the plane
trees with edge-lengths $(\mathbbl{t},(\ell_u,u\in\mathbbl{t}))$ and
$(\mathbbl{t}',(\ell'_u,u\in\mathbbl{t}'))$ are equivalent if there
exist permutations
$\bsigma=(\sigma_u,u\in\mathbbl{t})$ such that $\bsigma\mathbbl {t}=\mathbbl{t}'$ and
$\ell'_{\bsigma(u)}=\ell_u$, for every $u\in\mathbbl{t}$. We let
$\bTheta$
\label{treeswithedgelengths}
be
the set of trees with edge-lengths, that is, of equivalence classes of
plane trees with edge-lengths. There is a natural concatenation
transformation, similar to $\langle\cdot\rangle$, for elements of
$\bTheta$. Namely, if $\theta^{(i)}=(\mathbbl{t}^{(i)},(\ell
^{(i)}_u,u\in
\mathbbl{t})),1\leq i\leq k$, is a sequence of plane trees with
edge-lengths and
$\ell\geq0$, let
\[
\bigl\langle\theta^{(i)},1\leq i\leq k\bigr\rangle_\ell=\bigl(\mathbbl{t},(\ell
_u,u\in
\mathbbl{t})\bigr)
\]
be defined by
\[
\mathbbl{t}=\bigl\langle\mathbbl{t}^{(i)},1\leq i\leq k\bigr\rangle
\]
and
\[
\ell_{\varnothing}=\ell,\qquad \ell_{iu}=\ell^{(i)}_u ,\qquad
1\leq i\leq k, u\in\mathbbl{t}^{(i)} .
\]
If we replace each $\theta^{(i)}$
by another equivalent plane tree with edge-lengths, then the resulting
concatenation is equivalent to the first one, so that this operation
is well defined for elements of $\bTheta$.

Let $\theta\in\bTheta$, and consider a plane representative
$(\mathbbl{t},(\ell_u,u\in\mathbbl{t}))$. We construct an $\R$-tree
$\TT$ by imagining
that the edge from $\mathrm{pr}(u)$ to $u$ has length $\ell_u$. Note
that this intuitively involves a new edge with length
$\ell_\varnothing$ pointing from the root $\rho$ of the resulting
$\R$-tree to $\varnothing$ (this is sometimes called \textit{planting}).
Formally, $\TT$ is the isometry-equivalence class of a subset of
$\R^\mathbbl{t}$ endowed with the $l^1$-norm $\|(x_u,u\in\mathbbl {t})\|_1=\sum_{u\in
\mathbbl{t}}|x_u|$, defined as the union of segments
\[
\bigcup_{u\in\mathbbl{t}}\biggl[\sum_{v\prec u}\ell_ve_v,\sum
_{v\prec
u}\ell_ve_v+\ell_u e_u \biggr] ,
\]
where $(e_u,u\in\mathbbl{t})$ is the
canonical basis of $\R^\mathbbl{t}$ and $v\prec u$ means that $v$ is a strict
ancestor of $u$ in $\mathbbl{t}$. This $\R$-tree is naturally rooted
at $0\in
\R^\mathbbl{t}$. Of course, its isometry
class does not depend on the choice of the plane representative of
$\theta$, and can be written $\TT(\theta)$
\label{rtreetheta}\vadjust{\goodbreak}
unambiguously. Note that
there is a natural ``embedding'' mapping $\iota:\mathsf{t}\to\TT
(\theta)$
inherited from
%
%e11 ###
\begin{equation}\label{eq:10}
\iota_0\dvtx\mathbbl{t}\to\TT,\qquad \iota_0(u)=\sum_{v\preceq
u}\ell_ve_v ,
\end{equation}
and the latter is an isometry if $\theta$ is endowed with the
(pseudo-)metric $d_\theta$ on its vertices, defined by
\[
d_\theta(u,v)=\sum_{w\preceq u\ \mathrm{xor}\ w\preceq v}\ell_w
,
\]
where xor denotes ``exclusive or.''

Conversely, it is an elementary exercise to see that any rooted
$\R$-tree $\TT$ with a finite number of leaves can be written in the
form $\TT=\TT(\theta)$ for some $\theta\in\bTheta$, which is in fact
unique. In the sequel, we will often identify the tree $\theta\in
\bTheta$ with the $\R$-tree $\TT(\theta)$. For instance, this
justifies the notation $\langle\TT^{(1)},\ldots,\TT^{(r)}\rangle
_\ell$
for $\R$-trees $\TT^{(1)},\ldots,\TT^{(r)}$\vspace*{1pt} with finitely many leaves
and for $\ell\geq0$, which stands for the $\R$-tree in which the
roots of $\TT^{(1)},\ldots,\TT^{(r)}$ have been identified, and
attached to a segment of length $\ell$ to a new root.

With a discrete tree $\mathsf{t}$, we canonically associate the tree with
edge-lengths $\theta$ in which all lengths are equal to $1$, and the
rooted $\R$-tree $\TT(\mathsf{t})=\TT(\theta)$.
%edge-lengths 1}
\label{rtreelength1}
In this case,
$d_\theta=d_{\mathrm{gr}}$ is the graph distance. Using the isometry
$\iota\dvtx\mathsf{t}\mapsto\TT(\mathsf{t})$, we get the following
statement, left as an
exercise to the reader.
\begin{proposition}\label{sec:trees-with-edge}
Viewing $\mathsf{t}\in\bT$ as the element
$(\mathsf{t},d_{\mathrm{gr}},\rho,\mu_{\partial\mathsf{t}})$ of
$\Mw$ as in Section
\ref{sec:r-trees-gromov}, and endowing $\TT(\mathsf{t})$ with the uniform
probability distribution on $\mathcal{L}(\TT(\mathsf{t}))$, it holds that
\[
d_{\mathrm{GHP}}(a\mathsf{t},a\TT(\mathsf{t}))\leq a ,\qquad
a>0 .
\]
\end{proposition}

Due to this statement, in order to prove that the Markov branching
tree $T_n$ with law $\bP_n^q$ converges after rescaling toward
$\TT_{\gamma,\nu}$, it suffices to show the same statement for the
$\R$-tree $\TT(T_n)$. We will often make the identification of $T_n$
with $\TT(T_n)$.

%s3.2.2 ###
\subsubsection{Partition-valued processes and
$\R$-trees}\label{sec:part-valu-proc}

Let $(\pi(t),t\geq0)$ be a process with values in $\mathcal P_C$,
$C\subset\N$ finite or infinite, which is nondecreasing and indexed
either by $t\in\Z_+$ or $t\in\R_+$, in which case we also assume
that $\pi(\cdot)$ is right-continuous. We assume that there exists
some $t_0>0$ such that $\pi(t_0)=\mathbb{I}_C$. Let $B\subseteq C$ be
finite. If $B=\{i\}$, we let
\[
D^\pi_{\{i\}}=\inf\bigl\{t\geq0\dvtx\{i\}\in\pi(t)\bigr\}
\]
be the first time where
$i$ is isolated in a singleton block, and for $\#B\geq2$, let \label{dpiB}
\[
D^\pi_{B}=\inf\{t\geq0\dvtx B\cap\pi(t)\neq B\} .
\]
We can build a tree with edge-lengths (and labeled leaves)
$\theta(\pi(\cdot),B)$ by the following inductive procedure:\vadjust{\goodbreak}
\begin{longlist}[(1)]
\item[(1)] if $B=\{i\}$ we let $\theta(\pi(\cdot), B)$
\label{thetapi}
be the tree $\bullet$ with
length $D^\pi_{\{i\}};$
\item[(2)] if $\#B\geq2$, we let
\[
\theta(\pi(\cdot),B)=\bigl\langle\theta\bigl( \pi(D^\pi_{B}+\cdot),B\cap
\pi_i(D^\pi_{B}) \bigr),1\leq i\leq b\bigr\rangle_{D^\pi_{B}} ,
\]
where $b$
is the number of blocks of $\pi(D^\pi_{B})$ that intersect $B$, and
which are denoted by $\pi_1(D^\pi_{B}),\ldots,\pi_b(D^\pi_{B})$.
\end{longlist}
Note that the previous labeling convention for blocks may not agree
with our usual convention of labeling with respect to order of least
element.

If $(\pi(t),t\in\Z_+)$ is indexed by nonnegative integers, and
satisfies $\pi(0)=\mathbbm{O}_C$, there is a similar construction with
trees rather than trees with edge-lengths. Namely, we let
$\mathsf{t}_{\pi(\cdot)}$ be defined by:
\begin{longlist}[(1)]
\item[(1)] $\mathsf{t}_{\pi(\cdot)}=\bullet$ if $\#C=1$;
\item[(2)] $\mathsf{t}_{\pi(\cdot)}=\langle\mathsf{t}_{\pi
_i(1)\cap\pi(\cdot+1)},1\leq
i\leq b\rangle$ otherwise, where $b$ is the number of blocks of
$\pi(1)$, denoted by $\pi_1(1),\ldots,\pi_b(1)$.
\end{longlist}
It is then easy to see that, with the notation of the previous
section,
%
%e12 ###
\begin{equation}\label{eq:12}
\TT\bigl(\mathsf{t}_{\pi(\cdot)}\bigr)=\TT(\theta(\pi(\cdot),C)) ,
\end{equation}
and one can view $\theta(\pi(\cdot),B)$ as the subtree of
$\mathsf{t}_{\pi(\cdot)}$ spanned by the root and the leaves with
labels in $B$.

%s3.2.3 ###
\subsubsection{Continuum fragmentation trees}\label{sec:continuous-model}

Let $(\Pi(t),t\geq0)$ be the self-similar fragmentation process with
index $-\gamma<0$ and dislocation measure $\nu$. The formation of dust
property alluded to in Section~\ref{sec:self-simil-fragm} amounts to
the fact that almost surely, there exists some time $t_0>0$ such that
$\Pi(t)=\mathbb{I}_\N$ for every $t\geq t_0$. Consequently, the
construction of the previous paragraph applies with $C=\N$, and allows
us to construct a family of $\R$-trees
\[
\mathcal{R}_{B}=\theta(\Pi(\cdot),B)
\]
indexed by finite subsets
$B\subset\N$. Recall that a tree $\theta\in\bTheta$ has been
identified with $\TT(\theta)\in\Tw$. These $\R$-trees have finitely
many leaves that are naturally indexed by elements of $B$. Moreover,
they satisfy an obvious consistency property, meaning that taking the
subtree spanned by the root and the leaves indexed by $B'\subset B$
yields an $\R$-tree with same law as $\mathcal{R}_{B'}$. This is the
key to the definition of the fragmentation tree
$\TT_{\gamma,\nu}$.
\begin{proposition}[(\cite{HM04})]\label{sec:cont-fragm-trees-1}
Conditionally given $\TT_{\gamma,\nu}=(\TT,d,\rho,\mu)$, let
$L_1$, $L_2,\ldots$ be an i.i.d. sequence of leaves of $\TT$ distributed
according to $\mu$. Then for every finite $B\subset\N$, the \textit{reduced subtree}
\[
\mathcal{R}(\TT_{\gamma,\nu},B)=\bigcup_{i\in B}\lbr \rho,L_i\rbr
\]
has same distribution as $\mathcal{R}_B$.\vadjust{\goodbreak}

Moreover, the law of $\TT_{\gamma,\nu}$ is the only one having this
property, among distributions on $\Tw$ supported on the set of
elements satisfying properties (1), (2) and (3) in Section
\ref{sec:self-simil-fragm}.
%
%that charge only the set of
%$\{(\TT,d,\rho,\mu)\in\Tw: \forall x\in\TT,
%x\notin\mathcal{L}(\TT)\implies\mu(\TT_x)>0\mbox{ and
%}\mu(\mathcal{L}(\TT))=1\}$.
\end{proposition}

As an easy consequence, we have the following ``converse
construction'' of fragmentations from $\TT_{\gamma,\nu}$. With the
notation of the proposition, for every $t\geq0$, let $\Pi(t)$ be the
partition of $\N$ such that $i,j$ are in the same block of $\Pi(t)$ if
and only if $d(\rho,L_i\wedge L_j)>t$. Then $(\Pi(t),t\geq0)$ is a
self-similar fragmentation process with dislocation measure $\nu$ and
index $-\gamma$.

Also, note that the \textit{reduced trees}
$\mathcal{R}(\TT_{\gamma,\nu},B)$ rooted at $\rho$ and endowed with
the empirical measure
\[
\mu_B=\frac{1}{\#B}\sum_{i\in B}\delta_{L_i}
\]
converge in
distribution as $\#B\to\infty$ in $\Tw$ toward
$(\TT,d,\rho,\mu)$. In fact, the convergence holds a.s. if $B=[k]$
with $k\to\infty$: this is a simple exercise using the fact that
$\{L_i,i\geq1\}$ is a.s. dense in $\mathcal{L}(\TT)$ [by property
(3) in the definition of
$\TT_{\gamma,\nu}$], and the weak convergence of $\mu_{[k]}$ to
$\mu$
as $k\to\infty$.

The following statement gives a decomposition of the reduced tree
$\mathcal{R}(\TT,[k])$ at its first branchpoint above the root. Recall
the notation $D_k=\inf\{t\geq0\dvtx\Pi(t)\in A_k\}$.
\begin{proposition}\label{sec:cont-fragm-trees-2}
Let $k\geq2$ and $\pi=\Pi(D_k),\pi'=\pi|_{[k]}, b=b(\pi')$. Then
conditionally on $\{\pi,D_k\}$, the reduced tree
$\mathcal{R}(\TT_{\gamma,\nu},[k])$ has same distribution as
\[
\TT\bigl(\bigl\langle|\pi_i|^\gamma\mathcal{R}\bigl(\TT^{(i)},\pi
'_i\bigr),1\leq
i\leq b\bigr\rangle_{D_k}\bigr) ,
\]
where the $\TT^{(i)}$ are i.i.d. with same distribution as
$\TT_{\gamma,\nu}$, independent of $\sigma\{\pi,D_k\}$.

Moreover, for every $i\in\N$, the tree
$\mathcal{R}(\TT_{\gamma,\nu},\{i\})$ has the same distribution as the
$\R$-tree associated with the tree $(\varnothing,D_1)\in\bTheta$,
that is, a real segment with length $D_1=\inf\{t\geq0\dvtx\{1\}\in
\Pi(t)\}$.
\end{proposition}
\begin{pf}
The second statement is just a matter of
definitions, so we only need to prove the first one. By Proposition
\ref{sec:self-simil-fragm-2}, the process $\Pi(D_k+\cdot)$, in
restriction to the blocks containing at least one element in $[k]$,
has same distribution as the partitions-valued process whose blocks
are those of $\pi_i\cap\Pi^{(i)}(|\pi_i|^{-\gamma}\cdot),1\leq
i\leq
b$, for i.i.d. copies $\Pi^{(i)},i\geq1$, of $\Pi$, independent on
$\pi,D_k$. Therefore, one gets from the definition of $\mathcal{R}_B$
that
\[
\mathcal{R}_{[k]}\stackrel{(\mathrm{d})}{=}\TT\bigl\langle\theta\bigl(
\Pi^{(i)}(|\pi_i|^{-\gamma}\cdot),\pi'_i\bigr),1\leq i\leq
b\bigr\rangle_{D_k} ,
\]
from which the result follows immediately.
\end{pf}

Note that Proposition~\ref{sec:poiss-constr-homog-1} gives the joint
distribution of $D_k,|\pi_i|,1\leq i\leq b,\pi'$, as a special case,
while Proposition~\ref{sec:self-simil-fragm-3} characterizes the law\vadjust{\goodbreak}
of $D_1$, since it is the first time where the process
$(|\Pi_{(1)}(t)|,t\geq0)$ attains $0$. This, together with the
previous proposition, allows us to characterize entirely the laws of the
reduced trees of $\TT_{\gamma,\nu}$, hence the law of
$\TT_{\gamma,\nu}$ itself.

%s3.2.4 ###
\subsubsection{Markov branching trees as discrete
fragmentation trees}\label{sec:discrete-model}

Recall the informal description of Markov branching trees $\bP^q_n$ in
the \hyperref[sec1]{Introduction}, relying on collections of balls in urns. Rather than
collections of indistinguishable balls that split randomly, it is
convenient to consider instead a collection of balls that are
distinguished by a random, exchangeable labeling. This is achieved by
replacing partitions of integers by partitions of sets. We start with
a preliminary lemma.
\begin{lemma}\label{sec:partitions-1}
Let $n\geq1$ be fixed, as well as a partition $\lambda\in\P_n$ with
$p=p(\lambda)$ parts.

\begin{longlist}
\item
There are
\[
C_\lambda=\frac{n!}{\prod_{i=1}^{p}\lambda_i!\prod_{j=1}^n
m_j(\lambda)!}
\]
partitions $\pi\in\P_{[n]}$ such that
$\lambda(\pi)=\lambda$.

\item If $1\leq k\leq n$, and $\pi'\in\P_{[k]}$ has $b$
blocks, then for every $i_1,\ldots,i_b\in\{1,2,\ldots,p\}$ pairwise
distinct, there are
\[
C_\lambda^{\pi'}(i_1,\ldots,i_b)=C_\lambda
\frac{1}{(n)_{k}}\prod_{j=1}^b(\lambda_{i_j})_{\#\pi'_j}\prod
_{l\geq
1}\frac{m_l(\lambda)!}{m_l(\lambda_i,i\notin\{i_1,\ldots,i_b\})!}
\]
partitions $\pi\in\P_{[n]}$ such that $\lambda(\pi)=\lambda$,
$\pi|_{[k]}=\pi'$ and $\#\pi_j=\lambda_{i_j},1\leq j\leq b$.
\end{longlist}
\end{lemma}
\begin{pf}
Let $p$ be the\vspace*{1pt} number of parts of
$\lambda$. Then there are $p!/\prod_{j=1}^{n}m_j(\lambda)!$ sequences
$(c_1,\ldots,c_p)$ whose nonincreasing rearrangement is
$\lambda$. With any such sequence, we can associate
\[
\frac{n!}{\prod_{i=1}^pc_i!}=\frac{n!}{\prod_{i=1}^p\lambda_i!}
\]
sequences of the form $(B_1,\ldots,B_p)$ such that
$\{B_1,\ldots,B_p\}$ is a partition of $[n]$ (beware that the labeling
of the blocks $B_i$ will differ, in general, from labeling convention
described above for the blocks of a partition), with $\#B_i=c_i,1\leq
i\leq p$. Finally, exactly $p!$ sequences of the form
$(B_1,\ldots,B_p)$ induce the same partition
$\{B_1,\ldots,B_p\}$. Putting things together easily yields the formula
for~$C_\lambda$.

For the second formula, if $\lambda\in\P_n,\pi'\in\P_{[k]}$ and
$i_1,\ldots,i_b$ are given with $b=b(\pi')$, then any partition
$\pi\in\P_{[n]}$ with $\lambda(\pi)=\lambda$ and $\pi|
_{[k]}=\pi'$
must have $\pi_i|_{[k]}=(\pi|_{[k]})_i=\pi'_i$, for $1\leq i\leq b$,
the first equality coming from our choice of the labeling of blocks of
partitions. The restriction of $\pi$ to $[k]$ is thus entirely
determined. The blocks $\pi'_1,\ldots,\pi'_b$ should be completed with,
respectively, $\lambda_{i_1}-\#\pi'_1,\ldots,\lambda_{i_b}-\#\pi'_b$
elements of $[n]\setminus[k]$ to form the blocks
$\pi_1,\ldots,\pi_b$, while the remaining subset of $[n]\setminus[k]$
should be partitioned in such a way that the block sizes are given by
the sequence $(\lambda_i,i\notin\{i_1,\ldots,i_b\})$. There are
\[
\frac{(n-k)!}{\prod_{j=1}^b(\lambda_{i_j}-\#\pi'_j)!\prod_{i\notin
\{i_1,\ldots,i_b\}}\lambda_i!\prod_{l\geq
1}m_l(\lambda_{i},i\notin\{i_1,\ldots,i_b\})! }
\]
such partitions,
and this can be rewritten as $C^{\pi'}_\lambda(i_1,\ldots,i_b)$.
\end{pf}

Going back to Markov branching trees, let $B\subset\N$ have $n\geq2$
elements. Let $q=(q_n,n\geq1)$ satisfy~(\ref{eq:6}), and also assume
that $q_1(\varnothing)=1$. For every $\pi\in\P_B$, set \label{pbpi}
%
%e13 ###
\begin{equation}\label{eq:4}
p_B(\pi)=\frac{ q_n(\lambda(\pi))}{C_{\lambda(\pi)}} ,
\end{equation}
where $C_\lambda$ is the constant appearing in Lemma
\ref{sec:partitions-1}. Given the partition of $n$ that it induces
(which has distribution $q_n$), a $p_B$-distributed partition is thus
uniform among possible choices of partitions of $B$. In particular, a
random partition with distribution $p_B$ is exchangeable; that is, its
law is invariant under the action of permutations of $B$. By
convention, the law $p_B$, if $B=\{i\}$ is a singleton, is the Dirac
mass at the partition $\{\{i\}\}$.

For every $\pi\in\P_B$ with blocks $\pi_1,\pi_2,\ldots,\pi_k$ say,
consider random partitions $\widetilde{\pi}^i,1\leq i\leq k$, of
$\pi_1,\ldots,\pi_k$, respectively, chosen independently with
respective distributions $p_{\pi_1},\ldots,p_{\pi_k}$. We let
$Q(\pi,\cdot)$ be the law of the partition of $B$ made up of the
collection of all blocks of $\widetilde\pi^i$, $1 \leq i \leq k$.
%$\bigcup_{1\leq i\leq k}\widetilde{\pi}^i\in\P_B$.
Then $Q$ is the transition kernel of a
Markov chain on $\P_B$ (for any finite $B\subset\N$), that ends at the
state $\mathbb{I}_B$. It is easily seen that this Markov chain is
exchangeable as a process. Moreover, the chain started from the state
$\{B\}$, with $\#B = n$ has same distribution as the image of the chain
started from $[n]$ under the action of any bijection $[n]\to B$.

\begin{figure}

\includegraphics{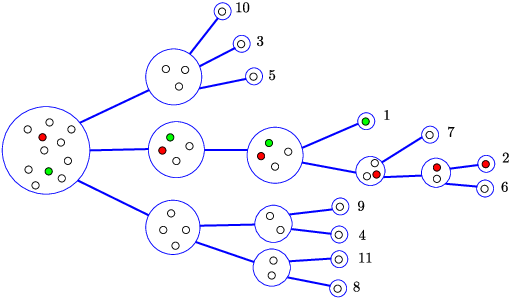}

\caption{A sample tree $T_{[n]}$ for $n=11$, with the labeled
leaves. The process $\Pi^{(11)}$ can be easily deduced: for
instance,
$\Pi^{(11)}(1)=\{\{1,2,6,7\},\{3,5,10\},\{4,8,9,11\}\}$. As opposed to
Figure \protect\ref{fig:CMB}, leaves are all connected to vertices
with at
least $2$ children, because of the requirement $q_1(\varnothing)=1$.}
\label{fig2}
\end{figure}

For finite $C\subset\N$, we let $(\Pi^C(r),r\geq0)$ be the chain
with transition matrix $Q$ and started from
$\Pi^C(0)=\mathbbm{O}_C$. Plainly, $\Pi^C$ is nondecreasing and
attains $\mathbb{I}_C$ in finite time a.s., so the construction of
Section~\ref{sec:part-valu-proc} applies and yields a family
$\theta(\Pi^C(\cdot),B)\in\bTheta,B\subseteq C$, as well as a tree
$T_C:=\mathsf{t}_{\Pi^C(\cdot)}$ (see Figure~\ref{fig2} for an example).
By construction, given that $\Pi^{C}(1)$
has blocks $\pi_1,\ldots,\pi_b$, the trees $\mathsf{t}_{\pi_i\cap
\Pi^C(\cdot+1)},1\leq i\leq b$, are independent\vspace*{1pt} with same
%
%f2 ###
distribution as $T_{\pi_i},1\leq i\leq b$, respectively. Since\vspace*{1pt} the law
of the nonincreasing rearrangement of $\#\pi_i,1\leq i\leq b$, is
$q_{\#C}$, we readily obtain the following statement.\footnote{There is
one subtlety in this statement, which is in the case $C=\{i\}$ for
some $i\in\N$. Indeed, by construction we have $T_C=\bullet$ a.s.,
and this is the only place where we have to require that
$q_1(\varnothing)=1$.}
\begin{lemma}\label{sec:mark-branch-trees-2}
The tree $T_C$ has law $\bP^q_{\#C}$.
\end{lemma}

In fact, the leaves of the tree $T_C$ are naturally labeled by
elements of $C$. We will use this in the sequel, without further
formalizing the notion of trees with labeled leaves.

We will also use the shorthand notation $T_C^B$ for the \textit{reduced
tree} $\theta(\Pi^C(\cdot),B)$. Using the above description, and
applying the Markov property for $\Pi^C$ at time $D^{\Pi^C}_B$ and the
particular form of the Markov kernel $Q$, we immediately obtain the
following, in the particular case $B=[k],C=[n]$.
\begin{proposition}\label{sec:discrete-model-1}
Let $2\leq k\leq n$. %, and let $\pi'\in\P_{[k]}$ have at least two
%blocks.
Then, conditionally on $D^{\Pi^{[n]}}_{[k]}=\ell$ and
$\Pi^{[n]}(D^{\Pi^{[n]}}_{[k]})=\pi$, with $\pi|_{[k]}=\pi'$, it holds
that $T_{[n]}^{[k]}$ has same distribution as
\[
\bigl\langle\theta^{(i)},1\leq i\leq b(\pi')\bigr\rangle_{\ell} ,
\]
where
$\theta^{(i)},1\leq i\leq b$, are independent with respective laws that
of $T_{\pi_i}^{\pi'_i},1\leq i\leq b(\pi')$.
\end{proposition}
\begin{pf}
The only subtle point is that $[k]\cap\pi_i=\pi'_i,1\leq i\leq b$,
since the labeling of the blocks of $\pi,\pi'$ could differ. But
since these partitions are, respectively, of $[n]$ and $[k]$, this
cannot be the case.
\end{pf}

%%%%%%%%%%%%%%%%%%%%%%%%%%%%%%%%%%%%%%%%%%%%%%%%%%%%%
%s4 ###
\section{\texorpdfstring{Proofs of Theorems \protect\ref{sec:main-result-1} and \protect\ref{sec:main-result-2}}
{Proofs of Theorems 5 and 6}}\label{sec:proof-theor-refs}
%%%%%%%%%%%%%%%%%%%%%%%%%%%%%%%%%%%%%%%%%%%%%%%%%%%%%

Let $q=(q_n,n\geq1)$ be a sequence of laws on~$\P_n$, respectively,
that satisfies~(\ref{eq:6}) and (H), for some fragmentation pair
$(-\gamma,\nu)$ and some slowly varying function $\ell$. In order to
lighten notation, we let $a_n = n^\gamma\ell(n)$.

Consider a sequence of trees $(T_n,n\geq1)$, where $T_n$ has
distribution $\bP_n^q$, $n\geq1$.
As we noticed in the \hyperref[sec1]{Introduction}, it is easy to pass from the
situation where $q_1(\varnothing)=1$ to the general situation, by adding
independent linear strings with
geometric$ (q_1(\varnothing))$-distributed lengths to the $n$ leaves of
$T_n$. Since geometric distributions have exponential tails, the
longest of these $n$ strings will have a length at most $C\log n$ with
probability going to $1$ as $n\to\infty$, for some $C>0$. If we let
$T_n^1$ be the tree for which $q_1(\varnothing)>0$ and $T_n^2$ the one
for which $q_1(\varnothing)=1$, coupled in the way depicted above, we
easily get, for any $\gamma>0$,
\[
\mathbb P\bigl(d_{\mathrm{GHP}}(a_n^{-1}T_n^1,a_n^{-1}T_n^2)\leq
Ca_n^{-1}\log n\bigr)\mathop{\longrightarrow}_{n\to\infty} 1 .
\]
Thus,
we can deduce the convergence in distribution of $a_n^{-1}T_n^1$ to
$\TT_{\gamma,\nu}$ from that of $a_n^{-1}T_n^2$. Therefore,\vspace*{1pt} from now
on and until the end of the present section, we make the following
hypothesis, which will allow us to apply Lemma
\ref{sec:mark-branch-trees-2}:
\begin{longlist}[(H$'$)]
\item[(H$'$)] The sequence $(q_n,n\geq1)$ satisfies
(H) and $q_1(\varnothing)=1$.
\end{longlist}
%
%Said otherwise, this also means that if a vertex of $T_n$ has only one
%child, then this child is not a leaf.

%s4.1 ###
\subsection{Preliminary convergence lemmas}\label{sec:prel-conv-lemm}

We now establish a couple of intermediate convergence results for the
discrete model. Recall that the sequence of distributions $q_n,n\geq
2$, on $\P_n$, respectively, induce distributions $p_B$ on $\P_B$ for
finite $B$ by formula~(\ref{eq:4}). By convention we set
$p_n=p_{[n]}$.
\begin{lemma}\label{sec:prel-conv-lemm-1}
Let $k\geq2$, and let $\pi'$ be an element in $\P_{[k]}$ with $b$
blocks, $b\geq2$. Let
%f:\sd\to\R,
$g\dvtx(0,\infty)^b\to\R$ be a continuous function with compact support.
Then, under assumption \textup{(H$'$)},
\[
a_np_n\biggl(
g\biggl(\frac{\#\pi_1}{n},\ldots,\frac{\#\pi_b}{n}\biggr)\ind_{\{
\pi|_{[k]}=\pi'\}}
\biggr)\mathop{\longrightarrow}_{n\to\infty}
\int_{\P_\N}\kappa_\nu(\d\pi)
g(|\pi_1|,\ldots,|\pi_b|)\ind_{\{\pi|_{[k]}=\pi'\}} ,
\]
where
$\kappa_\nu$ is the paintbox construction associated with $\nu$. Note
that on the event $\{\pi|_{[k]}=\pi'\}$, the quantities $\#\pi_i/n$
and $|\pi_i|$ for $1\leq i\leq b$ that appear above are
a.e. nonzero, respectively, under $p_n$ and $\kappa_\nu$.
\end{lemma}
\begin{pf}
For simplicity, we let
\[
B_n=p_n\biggl(
g\biggl(\frac{\#\pi_1}{n},\ldots,\frac{\#\pi_b}{n}\biggr)\ind_{\{
\pi|_{[k]}=\pi'\}}
\biggr) .
\]
Let $\lambda\in\P_n$, and let $i_1,\ldots,i_b\in\N$ be pairwise
distinct. This induces a sequence
$(\lambda_{i_1},\ldots,\lambda_{i_b})$. Note that there are exactly
\[
\prod_{l\geq1}\frac{m_l(\lambda)!}{m_l(\lambda_i,i\notin\{
i_1,\ldots,i_b\})!}
\]
choices of such pairwise distinct indices $i'_1,\ldots,i'_b$ such that
$(\lambda_{i'_1},\ldots,\lambda_{i'_b})=(\lambda_{i_1},\ldots
,\lambda_{i_b})$.
Hence, by the definition of $q_n$ and Lemma~\ref{sec:partitions-1},
\begin{eqnarray*}
B_n&=&\sum_{\lambda\in\P_n}q_n(\lambda)
\mathop{\sum_{i_1,\ldots,i_b\geq1}}_{\mathrm{pairwise}\
\mathrm{distinct}}\prod_{l\geq1}\frac{m_l(\lambda_i,i\notin\{
i_1,\ldots,i_b\})!}{m_l(\lambda)!} g\biggl(\frac{\lambda
_{i_1}}{n},\ldots,\frac{\lambda_{i_b}}{n}\biggr)\\
&&\hspace*{117.8pt}{}\times\frac{C_
\lambda^{\pi'}(i_1,\ldots,i_b)}{C_\lambda}
\\ &=&\sum_{\lambda\in\P_n}q_n(\lambda)
\frac{1}{(n)_k}\mathop{\sum_{i_1,\ldots,i_b\geq1}}_{
\mathrm{pairwise}\ \mathrm{distinct}}g\biggl(\frac{\lambda_{i_1}}{n},\ldots,\frac
{\lambda_{i_b}}{n}\biggr)
\prod_{j=1}^b(\lambda_{i_j})_{\#\pi'_j} .
\end{eqnarray*}
Now, the function
\[
h(\bs)=
\mathop{\sum_{i_1,\ldots,i_b\geq1}}_{\mathrm{pairwise}\
\mathrm{distinct}}g(s_{i_1},\ldots,s_{i_b})
\prod_{j=1}^bs_{i_j}^{\#\pi'_j} ,\qquad \bs\in\sd,
\]
is continuous
and bounded, because $g$ is compactly supported in $(0,\infty)^b$, so
that the sum is really a finite sum. Moreover,
%
%e14 ###
\begin{equation}\label{eq:5}\qquad
h(\bs)%\leq h_0(\bs)
\leq K\mathop{\sum_{0\leq
k_1,k_2,\ldots< k}}_{k_1+k_2+\cdots=k}\frac{k!}{\prod_{j\geq
1}k_j!}\prod_{j\geq1}s_j^{k_j}=K\biggl(1-\sum_{j\geq
1}s_j^k\biggr)\leq kK(1-s_1) ,
\end{equation}
where $K$ is an upper-bound of $|g|$, and for every $\lambda\in\P_n$,
it is easily checked that for large $n$, if $\eps>0$ is such that
$g(x_1,\ldots,x_b)=0$ as soon as \mbox{$\min_{1\leq i\leq b} x_i\leq\eps$},
\begin{eqnarray*}
\biggl(1-\frac{k}{\eps n}\biggr)^kh(\lambda/n)&\leq&
\frac{1}{(n)_k}\mathop{\sum_{i_1,\ldots,i_b\geq1}}_{\mathrm{
pairwise}\ \mathrm{distinct}}g\biggl(\frac{\lambda_{i_1}}{n},\ldots,
\frac{\lambda_{i_b}}{n}\biggr)
\prod_{j=1}^b(\lambda_{i_j})_{\#\pi'_j}\\
&\leq&\biggl(\frac{n}{n-k}\biggr)^kh(\lambda/n) .
\end{eqnarray*}
Letting $n\to\infty$
and applying \textup{(H$'$)}, which is validated by~(\ref{eq:5}),
\[
\lim_{n\to\infty}a_nB_n=
\int_\sd\nu(\d\bs)h(\bs) = \int_{\P_\N}\kappa_\nu(\d\pi)
g(|\pi_1|,\ldots,|\pi_b|)\ind_{\{\pi|_{[k]}=\pi'\}} ,
\]
the latter
equality being a simple consequence of the paintbox construction of
Section~\ref{sec:paintb-constr}.
\end{pf}

Now, we associate with $(q_n,n\geq1)$ a family of process
$(\Pi^B(r),r\geq0)$ with values in $\P_B$, as in Section
\ref{sec:discrete-model}. We let $\Pi^n=\Pi^{[n]}$ for simplicity, and
set
\[
D_k^n=D_{[k]}^{\Pi^n}=\inf\bigl\{r\geq0\dvtx[k]\cap\Pi^n(r)\neq\{[k]\}\bigr\}
\]
for $2\leq k\leq n$, and $D_1^n=D^{\Pi^n}_{\{1\}}=\inf\{r\geq
0\dvtx\{1\}\in\Pi^n(r)\}$. Also, for $r\geq0$ we let
\[
X_n(r)=\#\Pi^n_{(1)}(r) .
\]

\begin{lemma}\label{sec:prel-conv-lemm-2}
Let $n,k\in\N$ be fixed, with $n\geq k\geq2$, and let $\pi'\in\P_{[k]}$
have $b \geq2$ blocks. Let $F,f$ be measurable nonnegative functions. Then
\begin{eqnarray*}
&&\mathbb E\bigl[F\bigl(X_n\bigl(\cdot\wedge(D_{k}^n-1)
\bigr)\bigr)f\bigl(\#\Pi^n_i(D_{k}^n),1\leq i\leq b\bigr)\ind_{\{[k]\cap
\Pi^n(D_{k}^n)=\pi'\}}\bigr]\\
&&\qquad=\sum_{r'>
0}\mathbb E\biggl[\frac{(X_n(r'-1)-1)_{k-1}}{(n-1)_{k-1}}\\
&&\qquad\hphantom{=\sum_{r'>
0}\mathbb E\biggl[}
{}\times F
\bigl(X_n\bigl(\cdot\wedge
(r'-1)\bigr)\bigr) p_{X_n(r'-1)}\bigl(f(\#\pi_i,1\leq i\leq
b)\ind_{\{\pi|_{[k]}=\pi'\}}\bigr)\biggr].
\end{eqnarray*}
\end{lemma}
\begin{pf}
We first consider an expression of a more
general form. For nonnegative functions $G,g$, we have, using the
Markov property at time $r'-1$ in the second step,
\begin{eqnarray*}
\hspace*{-4pt}&&\mathbb E\bigl[G\bigl(\Pi^n\bigl(\cdot\wedge(D_{k}^n-1)\bigr)\bigr)
g\bigl(\Pi^n_{(1)}(D_{k}^n-1)\cap\Pi^n(D_{k}^n)\bigr)\bigr]\\
\hspace*{-4pt}&&\qquad=\sum_{r'>
0}\mathbb E\bigl[G\bigl(\Pi^n\bigl(\cdot\wedge(r'-1)\bigr)\bigr)\\
\hspace*{-4pt}&&\qquad\hphantom{=\sum_{r'>
0}\mathbb E\bigl[}{}
\times\ind_{\{
[k]\subset
\Pi^n_{(1)}(r'-1)\}}g\bigl(\Pi^n_{(1)}(r'-1)\cap\Pi^n(r')\bigr)\ind
_{\{[k]\cap\Pi^n(r')\neq
\{[k]\}\}}\bigr]\\
\hspace*{-4pt}&&\qquad=\sum_{r'>0}\mathbb E\bigl[G\bigl(\Pi^n\bigl(\cdot
\wedge(r'-1)\bigr)\bigr)\ind_{\{[k]\subset
\Pi^n_{(1)}(r'-1)\}}p_{\Pi^n_{(1)}(r'-1)}\bigl(g(\pi)\ind_{\{[k]\cap
\pi\neq\{[k]\}\}}\bigr)\bigr].
\end{eqnarray*}
Specializing this formula to $G$ depending only on $X_n$ and
$g(\pi)=f(\#\pi_1,\ldots,\break\#\pi_b)\ind_{\{\pi|_{[k]}=\pi'\}}$, and
using obvious exchangeability properties, we obtain
\begin{eqnarray*}
&&\mathbb E\bigl[F\bigl(X_n\bigl(\cdot\wedge(D_{k}^n-1)\bigr)\bigr)
f\bigl(\#\Pi^n_i(D_{k}^n),1\leq i\leq b\bigr)\ind_{\{[k]\cap
\Pi^n(D_{k}^n)=\pi'\}}\bigr]\\
&&\qquad= \sum_{r'>
0}\mathbb E\bigl[F\bigl(X_n\bigl(\cdot\wedge(r'-1)\bigr)\bigr)\\
&&\qquad\hphantom{= \sum_{r'>
0}\mathbb E\bigl[}{}
\times p_{X_n(r'-1)}\bigl(f(\#\pi_i,1\leq i\leq b)\ind_{\{ \pi|_{[k]}=\pi'
\}}\bigr)\ind_{\{[k]\subset\Pi^n_{(1)}(r'-1)\}}\bigr].
\end{eqnarray*}
All the terms in the expectation depend on $(X_n(r),0\leq r\leq
r'-1)$, except the last one which is a function of
$\Pi^{n}_{(1)}(r'-1)$. But by Lemma~\ref{sec:exch-part-finite-1} (in
fact, the variant used in the proof of Lemma
\ref{sec:self-simil-fragm-5}),
\[
\mathbb P\bigl([k]\subset\Pi^n_{(1)}(r'-1) | \bigl(X_n(r),0\leq
r\leq
r'-1\bigr)\bigr) =\frac{(X_n(r'-1)-1)_{k-1}}{(n-1)_{k-1}}
\]
giving the
result.
\end{pf}

In the sequel, $\Pi(\cdot)$ will denote a continuous-time self-similar
fragmentation with characteristic pair $(-\gamma,\nu)$, and
$D_k,k\geq1$, will be defined as in Section~\ref{sec:self-simil-fragm-1}.
\begin{lemma}\label{sec:prel-conv-lemm-3}
Under assumption \textup{(H$'$)}, it holds that
\[
\biggl(\frac{X_n(\lfloor a_n t\rfloor)}{n}, t\geq
0\biggr)\mathop{\longrightarrow}^{(d)}_{n\to\infty}
\bigl(\bigl|\Pi_{(1)}(t)\bigr|,t\geq0\bigr) ,
\]
in distribution for the Skorokhod
topology, jointly with the convergence
\[
\frac{1}{a_n}D^n_{1}
\mathop{\longrightarrow}^{(d)}_{n\to\infty}D_1 .
\]
\end{lemma}
\begin{pf}
For $n> k\geq1$, let
$p_{n,k}=\mathbb P(X_n(1)=k)$. Note that the process $X_n$ is a nonincreasing
Markov chain started from $n$, with probability transitions $p_{i,j}, 1
\leq j \leq i$. Then by a simple exchangeability argument,
\[
p_{n,k}= \sum_{\pi\in\P_{[n]}} p_n(\pi)m_k(\pi)\frac{k}{n}
=\sum_{\lambda\in\P_{n}} q_n(\lambda)m_k(\lambda)\frac{k}{n},\qquad
1 \leq k \leq n,
\]
where
$m_k(\pi)=m_k(\lambda(\pi))$ is the number of blocks of $\pi$ with
size $k$. Consider the associated generating function for $x \geq0$,
\begin{eqnarray*}
F_n(x)&=&\sum_{k=1}^{n}
\biggl(\frac{k}{n}\biggr)^xp_{n,k}=\sum_{\lambda\in\P_{n}}
q_n(\lambda)\sum_{k=1}^{n}
m_k(\lambda)\biggl(\frac{k}{n}\biggr)^{x+1}\\
&=&\sum_{\lambda\in\P_{n}}
q_n(\lambda) \sum_{i\geq1}\biggl(\frac{\lambda_i}{n}\biggr)^{x+1} .
\end{eqnarray*}
Hence, $1-F_n(x)= \ov{q}_n(f)$,
where $f(\bs)=1-\sum_{i\geq1}s_i^{x+1}$. Note that $f\dvtx\sd\to\R$ is
continuous whenever $x>0$. Indeed, the norm
$\|\bs\|_{x+1}=(\sum_{i\geq1}s_i^{x+1})^{1/(x+1)}$ of any $\bs\in
\sd$ is finite, and satisfies for every $\bs,\bs'\in\sd$,
\begin{eqnarray*}
\bigl|
\|\bs\|_{x+1}-\|\bs'\|_{x+1}\bigr|&\leq& \|\bs-\bs'\|_{x+1}\\
&\leq& \sup_{i\geq1}|s_i-s'_i|^{x/(1+x)}\biggl(\sum_{i\geq
1}s_i+\sum_{i\geq1}s'_i\biggr)^{1/(x+1)}\\
&\leq& 2^{1/(x+1)}d(\bs,\bs')^{x/(x+1)} .
\end{eqnarray*}
Thus we may apply (H$'$), and obtain
\[
a_n\bigl(1-F_n(x)\bigr)\mathop{\longrightarrow}_{n\to\infty}
\int_\sd\biggl(1-\sum_{i\geq1} s_i^{x+1}\biggr) \nu(\d\bs) .
\]
This is exactly
what we need to use~\cite{HaMi09}, Theorem 1, stating that\break
$(n^{-1}X_n(\lfloor a_n t\rfloor)$, $t\geq0)$ converges in
distribution to the self-similar Markov process
$\exp(-\xi_{\tau(\cdot)})$, as defined around Proposition~\ref{sec:self-simil-fragm-3}. Moreover, this convergence holds jointly
with the convergence of absorption times at $1$, so
$a_n^{-1}D_{1}^n$ converges to the absorption time
at $0$ of $\exp(-\xi_{\tau(\cdot)})$. By Proposition~\ref{sec:self-simil-fragm-3}, the process $\exp(-\xi_{\tau(\cdot)})$
has same distribution as $(|\Pi_{(1)}(t)|,t\geq0)$, which reaches $0$
for the first time at time $D_1$. Hence the result.
\end{pf}

Finally, the combination of the last two lemmas gives the last of our
preliminary ingredients.
\begin{lemma}\label{sec:prel-conv-lemm-4}
The following joint convergence in distribution holds:
\begin{eqnarray*}
&&\biggl(\frac{D_k^n}{a_n},[k]\cap\Pi^n(D_k^n),
\biggl(\frac{\#\Pi^n_{(i)}(D_k^n)}{n},i\in
[k]\biggr)\biggr)\\
&&\qquad\mathop{\longrightarrow}^{(d)}_{n\to\infty}
\bigl(D_k,[k]\cap\Pi(D_k),\bigl(\bigl|\Pi_{(i)}(D_k)\bigr|,i\in
[k]\bigr)\bigr) .
\end{eqnarray*}
\end{lemma}
\begin{pf}
Let $\pi'\in\P_k$ have $b\geq2$ blocks, and
$f,g\dvtx(0,\infty)\to\R,h\dvtx(0,\infty)^b\to\R$ be continuous functions with
compact support. Then by Lemma~\ref{sec:prel-conv-lemm-2},
\begin{eqnarray*}
&&\mathbb E\biggl[f\biggl(\frac{
D_k^n}{a_n}\biggr)g\biggl(\frac{X_n(D_k^n-1)}{n}\biggr)
h\biggl(\frac{\#\Pi^n_i(D_k^n)}{X_n(D_k^n-1)},1\leq i\leq
b\biggr)\ind_{\{[k]\cap\Pi^n(D_k^n)=\pi'\}}\biggr]\\
&&\qquad=\sum
_{r'>0}f\biggl(\frac{r'}{a_n}\biggr)
\mathbb E\biggl[\frac{(X_n(r'-1)-1)_{k-1}}{(n-1)_{k-1}}g\biggl(\frac
{X_n(r'-1)}{n}\biggr)\\
&&\qquad\hphantom{=\sum
_{r'>0}f\biggl(\frac{r'}{a_n}\biggr)
\mathbb E\biggl[}
{}\times p_{X_n(r'-1)}\biggl(h\biggl(\frac{\#\pi_i}{X_n(r'-1)},1\leq i\leq
b\biggr)\ind_{\{\pi|_{[k]}=\pi'\}}\biggr)\biggr]\\
&&\qquad=\frac{1}{a_n}
\sum_{r'>0}f\biggl(\frac{r'}{a_n}\biggr) \mathbb E
\biggl[\Phi\bigl(n,X_n(r'-1)\bigr)
g\biggl(\frac{X_n(r'-1)}{n}\biggr)\Psi\bigl(X_n(r'-1)\bigr)\biggr]\\
&&\qquad=\int
_{1/a_n}^\infty
f\biggl(\frac{\lfloor a_nu\rfloor}{a_n} \biggr)\,\d u
\mathbb E\biggl[\Phi\bigl(n,X_n(\lfloor a_nu\rfloor-1)\bigr)\\
&&\qquad\hphantom{=\int
_{1/a_n}^\infty
f\biggl(\frac{\lfloor a_nu\rfloor}{a_n} \biggr)\,\d u
\mathbb E\biggl[}
\times g\biggl( \frac
{X_n(\lfloor a_n
u\rfloor-1)}{n}\biggr) \Psi\bigl(X_n(\lfloor a_n u
\rfloor-1)\bigr)\biggr],
\end{eqnarray*}
where
\[
\Phi(n,x)=\frac{(x-1)_{k-1}}{(n-1)_{k-1}}\frac{a_n}{a_{x}}
\mathop{\longrightarrow}_{(n,x/n)\to(\infty,c)}
c^{k-1-\gamma}
\]
and
\[
\Psi(m)=a_m p_m\biggl(h\biggl(\frac{\#\pi_i}{m},1\leq i\leq
b\biggr)\ind_{\{\pi|_{[k]}=\pi'\}}\biggr) .
\]
Note\vspace*{1pt} that the Potter's
bounds for regularly varying functions (\cite{BGT}, Theorem 1.5.6) imply
that $\Phi(n,x) \leq C( \frac{x}{n})^{k-1-\gamma-1}$ for
all $n \geq x \geq A$ for some finite positive constants $C,A$. In
particular\vspace*{1pt} there exists some $n_0$ such that $\sup_{n \geq n_0, 0<x
\leq n}\Phi(n,x)\times g(x/n)<\infty$ (since $g$ is null in a neighborhood
of 0). The joint use of Lemmas~\ref{sec:prel-conv-lemm-1} and
\ref{sec:prel-conv-lemm-3} entails by dominated convergence that the
expectation term in the integral converges to (note that the
quantities $|\pi_i|,1\leq i\leq b$, are all a.e. positive on
$\{\pi|_{[k]}=\pi'\}$ under $\kappa_\nu$)
\[
\mathbb E\biggl[\bigl|\Pi_{(1)}(u)\bigr|^{k-1-\gamma}g\bigl(\bigl|\Pi_{(1)}(u)\bigr|\bigr)
\int_{\P_\N}\kappa_\nu(\d\pi)h(|\pi_i|,1\leq i\leq
b)\ind_{\{\pi|_{[k]}=\pi'\}}\biggr] ,
\]
and since $f,g,h$ are
compactly supported, the whole integral converges to
\begin{eqnarray*}
&&\int_0^\infty f(u)\,\d u
\mathbb
E\biggl[\bigl|\Pi_{(1)}(u)\bigr|^{k-1-\gamma}g\bigl(\bigl|\Pi_{(1)}(u)\bigr|\bigr)\\
&&\qquad\hspace*{44.8pt}{}\times\int_{\P_\N}\kappa_\nu(\d\pi)h(|\pi_i|,1\leq i\leq
b)\ind_{\{\pi|_{[k]}=\pi'\}}\biggr],
\end{eqnarray*}
which, by Proposition
\ref{sec:poiss-constr-homog-1}, equals
\[
\mathbb E\biggl[f(D_k)g\bigl(\bigl|\Pi_{(1)}(D_k-)\bigr|\bigr)h\biggl(\frac{|\Pi_i(D_k)|}{
|\Pi_{(1)}(D_k-)|},1\leq i\leq b\biggr)\ind_{\{[k] \cap
\Pi(D_k)=\pi'\}}\biggr] .
\]
It is now easy to conclude, since
$|\Pi_i(D_k)|>0$ almost-surely.
\end{pf}

%s4.2 ###
\subsection{Convergence of finite-dimensional
marginals}\label{sec:conv-finite-dimens}

The first step in the proof of Theorem~\ref{sec:main-result-1} is the
following result on reduced trees $T_C^B$ of Section~\ref{sec:discrete-model}.
\begin{proposition}\label{sec:conv-finite-dimens-1}
Let $B\subset\N$ be finite. Under assumption \textup{(H$'$)}, we have the
following convergence in distribution in $\Tw$:
\[
\frac{1}{a_n}T_{[n]}^B
\mathop{\longrightarrow}^{(d)}_{n\to\infty}
\mathcal{R}(\TT_{\gamma,\nu},B) .
\]
\end{proposition}
\begin{pf}
We use an induction argument on the cardinality
of $B$. For $B=\{i\}$, one can assume by exchangeability (as soon as
$n\geq i$) that $B=\{1\}$, and in this case, the reduced tree is
$T_{[n]}^B=(\varnothing,D_1^n)$, while
$\mathcal{R}(\TT_{\gamma,\nu},\{1\})=(\varnothing,D_1)$ by
Proposition~\ref{sec:cont-fragm-trees-2}. By the second part of Lemma
\ref{sec:prel-conv-lemm-3}, under (H$'$), it holds that
\[
\frac{D^n_1}{a_n}\mathop{\longrightarrow}^{(\mathrm d)}_{n\to\infty}D_1 .
\]
This initializes the
induction. Now, assume that Proposition
\ref{sec:conv-finite-dimens-1} has been proved for every set $B$
with cardinality at most $k-1$, for some $k\geq2$. We want to
show that the same is true of any set of cardinality $k$, and by
exchangeability, we may assume that $B=[k]$.\vadjust{\goodbreak}

We now recall, using Proposition~\ref{sec:discrete-model-1}, that
conditionally on $D_k^n=\ell$, $[k]\cap\Pi^n(D_k^n)=\pi'$
having $b\geq2$ blocks and on $\Pi^n_i(D_k^n)=\pi_i,1\leq i\leq
b$, with respective cardinality $\#\pi_i = n_i$, the tree $T_{[n]}^{[k]}$
has same distribution as
\[
\bigl\langle\theta^{(i)},1\leq i\leq b\bigr\rangle_{\ell} ,
\]
where
$\theta^{(i)}$ has same distribution as $T^{\pi'_i}_{\pi_i}$, and
these trees are independent.

The joint\vspace*{1pt} distribution of $D_k^n,[k]\cap
\Pi^n(D_k^n),(\#\Pi^n_{(i)}(D_k^n),1\leq i\leq k)$ is
specified by Lemma~\ref{sec:prel-conv-lemm-2}, and its scaling limit
by Lemma~\ref{sec:prel-conv-lemm-4}. We obtain by the induction
hypothesis that jointly with the above convergence, conditionally on $
[k]\cap\Pi^n(D_k^n)=\pi'$,
\begin{eqnarray*}
\frac{1}{a_n}\theta^{(i)}&=&
\frac{a_{n_i}}{a_n}
\frac{1}{a_{n_i}}
\theta^{(i)}\\
&\displaystyle \mathop{\longrightarrow}^{(\mathrm{d})}_{n\to\infty}& |\Pi
_i(D_k)|^\gamma\TT^{(i)}
,\qquad 1\leq i\leq b ,
\end{eqnarray*}
where the $\TT^{(i)}$ are independent with same laws as
$\mathcal{R}(\TT_{\gamma,\nu},\pi'_i)$, respectively. Finally,
$a_n^{-1}T^{[k]}_{[n]}$ converges to
\[
\bigl\langle|\Pi_{i}(D_k)|^{\gamma}\TT^{(i)},1\leq i\leq
b\bigr\rangle_{D_k} ,
\]
and the $\R$-tree associated with this tree has
same distribution as $\mathcal{R}(\TT_{\gamma,\nu},[k])$ by
Proposition~\ref{sec:cont-fragm-trees-2}.
\end{pf}

%%%%%%%%%%%%%%%%%%%%%%%%%%%%%%%%%%%%%%%%
%s4.3 ###
\subsection{Tightness in the Gromov--Hausdorff topology}\label{sec:tightness}
%%%%%%%%%%%%%%%%%%%%%%%%%%%%%%%%%%%%%%%%

We now want to improve the convergence of Proposition
\ref{sec:conv-finite-dimens-1} into a convergence of nonreduced trees
for the Gromov--Hausdorff topology. Namely
\begin{proposition}\label{sec:tightn-grom-hausd-1}
Under hypothesis \textup{(H$'$)}, we have the convergence in distribution
\[
\frac{1}{a_n}T_n\mathop{\longrightarrow}^{(d)}_{n\to\infty}
\TT_{\gamma,\nu}
\]
in $\T$, for the Gromov--Hausdorff topology.
\end{proposition}

This will be proved by first showing a couple of intermediate lemmas.
\begin{lemma}\label{sec:tightness-2}
Under assumption \textup{(H$'$)}, we have the convergence in
distribution
\[
\bigl(\Pi^n(\lfloor a_nt\rfloor),t\geq0\bigr)
\mathop{\longrightarrow}^{(d)}_{n\to\infty}\bigl(\Pi(t),t\geq
0\bigr)
\]
jointly with
\[
\biggl(\frac{\#\Pi^n_{(i)}(\lfloor a_nt\rfloor)}{n},t\geq
0\biggr)\mathop{\longrightarrow}^{(d)}_{n\to\infty}
\bigl(\bigl|\Pi_{(i)}(t )\bigr|,t\geq0\bigr)
\]
for every $i\geq1$, all these
convergences holding jointly.\vadjust{\goodbreak}
\end{lemma}
\begin{pf}
The fact that $([k]\cap\Pi^n(\lfloor a_n
t\rfloor),t\geq0)$ converges in the Skorokhod space to $([k]\cap
\Pi(t),t\geq0)$ for every $k\geq1$ is obtained by using an inductive
argument similar to that used in the proof of Proposition
\ref{sec:conv-finite-dimens-1}. We only sketch the argument. The
statement is trivial for $k=1$, so we can assume that $k\geq2$. The
process $[k]\cap\Pi^n(\lfloor a_n\cdot\rfloor)$ remains constant
equal to $[k]$ up to time $a_n^{-1}D^n_k$, and jumps to the state
$\pi'=[k]\cap\pi$, $\pi=\Pi^n(D^n_k)$. By Lemma
\ref{sec:prel-conv-lemm-4}, $a_n^{-1}D^n_k\to D_k$ as $n\to\infty$,
and the latter has a diffuse law by Proposition
\ref{sec:poiss-constr-homog-1}.

After time $a_n^{-1}D^n_k$, given $\pi$, the restrictions $\pi'_i\cap
\Pi^n(\lfloor a_n\cdot\rfloor+D^n_k)$ have same distribution as
$\pi'_i\cap\Pi^{\pi_i}(\lfloor a_n\cdot\rfloor)$ and are
independent. By the induction hypothesis and exchangeability, still
conditionally on $\pi$, this converges to $\pi'_i\cap
\Pi^{(i)}(\cdot)$, where $\Pi^{(i)},i\geq1$, are i.i.d. copies of
$\Pi$. Moreover, since the jump times have diffuse laws, two such
copies never jump at the same time, from which one concludes that
given $\pi$, the process $(\pi'_i\cap\Pi^n(\lfloor
a_nt\rfloor+D^n_k),1\leq i\leq b(\pi'),t\geq0)$ converges in the
Skorokhod space to $(\pi'_i\cap\Pi^{(i)}(t),1\leq i\leq b(\pi
'),t\geq
0)$. This concludes the inductive step by gluing on the initial
constancy interval of the process, with length $a_n^{-1}D^n_k$.

The convergence of $\Pi^n(\lfloor a_n\cdot\rfloor)$ in the Skorokhod
space follows, because $d_{\mathcal{P}}([k]\cap\pi,\pi)\leq e^{-k}$
for every $\pi\in\P_\N$. This shows that $[k]\cap\Pi^n(\lfloor
a_n\cdot\rfloor)$ remains uniformly close to $\Pi^n(\lfloor
a_n\cdot\rfloor)$.

Next, by Lemma~\ref{sec:tightness-1}, it follows that, jointly with
this convergence, for every $i\geq1$, the one-dimensional marginals
of $(n^{-1}\#\Pi^n_{(i)}(\lfloor a_n t\rfloor),t\geq0)$ converge\vspace*{1pt} in
distribution to those of $(|\Pi_{(i)}(t)|,t\geq0)$, at least for
times which are not fixed discontinuity times of the limiting
process---the set of such points is always countable, and it turns out that
there are none in the present case. The convergence of
finite-dimensional marginals (outside of possible fixed
discontinuities) is obtained in a similar way, using a straightforward
generalization of Lemma~\ref{sec:tightness-1} to the case of a
sequence $((\pi^{(n,1)},\ldots,\pi^{(n,k)}),n\geq1)$ of jointly
exchangeable random partitions, respectively, of $[n]$, that converges
to a limiting $k$-tuple of random partitions of $\N$. This
generalization is left to the reader.

Since we also
know that the laws of the processes \mbox{$(n^{-1}\#\Pi^n_{(i)}(\lfloor a_n
t\rfloor),t\geq0)$} are tight when $n$ varies, by Lemma
\ref{sec:prel-conv-lemm-3} (these processes all have same distribution
as $(n^{-1}X_n(\lfloor a_n t\rfloor),t\geq0)$ by exchangeability),
this allows us to conclude.
\end{pf}

For $k+1\leq
i\leq n$, let
\[
S_i^n=\inf\bigl\{r\geq0\dvtx[k]\cap\Pi^n_{(i)}(r)=\varnothing\bigr\},
\]
the
first time when the ball indexed $i$ is separated from the $k$ first
balls. The random variables $S^n_i,k+1\leq i\leq n$, have same
distribution by exchangeability. The strong Markov property at the
stopping time $S_i^n$ also shows that conditionally on
$\Pi^n_{(i)}(S_i^n)=B$, the process $(B\cap\Pi^n(S_i^n+r),r\geq0)$
has same distribution as $\Pi^B$. Conditionally on\vadjust{\goodbreak}
$\Pi^n_{(i)}(S_i^n)=B$, the tree $\mathsf{t}_{B\cap
\Pi^n(S_i^n+\cdot)}$ has thus the same distribution as $T_B$ and can be
seen as a subtree of $T_{[n]}$, characterized by the fact that this
subtree contains the leaf labeled $i$, does not contain any of the
leaves labeled by an element of $[k]$ and is the maximal subtree of
$T_{[n]}$ with this property. In particular, the Gromov--Hausdorff
distance between $T_{[n]}^{[k]}$ and $T_{[n]}$ is at most
\[
d_{\mathrm{GH}}\bigl(T_{[n]}^{[k]},T_{[n]}\bigr)\leq\max_{k+1\leq i\leq
n}\mathrm{ht} \bigl(\mathsf{t}_{\Pi^n_{(i)}(S^n_i)\cap\Pi
^n(S_i^n+\cdot)}
\bigr) ,
\]
where $\mathrm{ht}(\mathsf{t})$, called the height of $\mathsf{t}$, is
the maximal height of a vertex in $\mathsf{t}$.

Note\vspace*{1pt} that if $j\in\Pi^n_{(i)}(S^n_i)$, then $S^n_j=S^n_i$. Therefore,
the blocks $\Pi^n_{(i)}(S^n_i), k+1\leq i\leq n$, are either disjoint
or equal. Moreover, the partition $\pi$ of $[n]\setminus[k]$ with
these blocks is clearly exchangeable. By putting the previous
observations together,
we obtain by first conditioning on $\pi$, and for every $\eta>0$,
%
%e15 ###
\begin{equation}\label{eq:13}
\mathbb P \bigl(d_{\mathrm{GH}}\bigl(T_{[n]}^{[k]},T_{[n]}\bigr)\geq\eta a_n
\bigr) \leq
\mathbb E\biggl[\sum_{i\geq1}\bP^q_{\#\pi_i}(\mathrm{ht}\geq\eta
a_n)\biggr]
,\qquad \eta>0 .
\end{equation}
Here and later, we adopt the convention that quantities involving
$\#\pi_i$ are always equal to $0$ when $\#\pi_i=0$. At this point, we
need the following uniform estimate for the height of a
$\bP^q_n$-distributed tree, which is the key lemma of this section.
\begin{lemma}\label{sec:tightn-grom-hausd}
Assume \textup{(H$'$)}. Then for all $p>0$, there exists a finite
constant $C_p$ such that
\[
\bP_n^q(\mathrm{ht} \geq x a_n ) \leq
\frac{C_p}{x^p}\qquad \forall x>0, \forall n \geq1.
\]
\end{lemma}

Before giving the proof of this statement, we end the proof of
Proposition~\ref{sec:tightn-grom-hausd-1}. Using Lemma
\ref{sec:tightn-grom-hausd} for $p=2/\gamma$ and~(\ref{eq:13}), we
obtain
\[
\mathbb P \bigl(d_{\mathrm{GH}}\bigl(T_{[n]}^{[k]},T_{[n]}\bigr)\geq\eta
a_n\bigr) \leq C_{2/\gamma} \eta^{-2/\gamma} \mathbb
E\biggl[\sum_{i\geq1} \frac{a_{\#\pi_i}^{2/\gamma}}{a_n^{2/\gamma
}}\biggr] .
\]
By the exchangeability of the
partition $\pi$ of $[n]\setminus[k]$, note that for every measurable
function $f$,
\[
\mathbb E\bigl[f\bigl(\#\pi_{(k+1)}\bigr)\bigr]=\frac{1}{n-k}\mathbb
E\Biggl[\sum_{i=k+1}^nf\bigl(\#\pi_{(i)}\bigr)\Biggr]= \mathbb E\biggl[\sum
_{i\geq
1}\frac{\#\pi_i}{n-k}f(\#\pi_i)\biggr] .
\]
This finally yields
\[
\mathbb P\bigl(d_{\mathrm{GH}}\bigl(T_{[n]}^{[k]},T_{[n]} \bigr)\geq\eta
a_n\bigr) \leq C_{2/\gamma} \eta^{-2/\gamma}\mathbb
E\biggl[\frac{a_{\#\pi_{(k+1)}}^{2/\gamma}(\#\pi
_{(k+1)})^{-1}}{a_n^{2/\gamma}n^{-1}}\biggr] .
\]
Since the sequence $(a_n
^{2/\gamma}n^{-1},n\geq1)$ is strictly positive and regularly
varying at $\infty$ with index 1, we get from Potter's bounds
(\cite{BGT}, Theorem 1.5.6),\vadjust{\goodbreak} the existence of a finite constant $C$ such
that $(a_k^{2/\gamma}k^{-1}) / (a_n
^{2/\gamma}n^{-1}) \leq C \sqrt{k/n}$ for
all $1 \leq k \leq n$. Hence,
\[
\mathbb P\bigl(d_{\mathrm{GH}}\bigl(T_{[n]}^{[k]},T_{[n]} \bigr)\geq\eta
a_n\bigr) \leq C C_{2/\gamma} \eta^{-2/\gamma}\mathbb
E\Biggl[\sqrt{\frac{\#\pi_{(k+1)} }{n}}\Biggr] .
\]
Note that the
quantity in the expectation is bounded by 1. By Lemma~\ref{sec:tightness-2}, it holds that $S^n_{k+1}/a_n \to
S_{k+1}$ in distribution as $n\to\infty$, where $S_{k+1}=\inf\{t\geq
0\dvtx[k]\cap\Pi_{(k+1)}(t)=\varnothing\}$. This convergence\vspace*{1pt} holds
jointly with that of $(n^{-1}\#\Pi^n_{(k+1)}(\lfloor a_nt\rfloor)$,
$t\geq0)$ to $(|\Pi_{(k+1)}(t)|,t\geq0)$ in the
Skorokhod space, whence we deduce that
\[
\limsup_{n\to\infty} \mathbb
P\bigl(d_{\mathrm{GH}}\bigl(T_{[n]}^{[k]},T_{[n]}\bigr)\geq\eta a_n\bigr)
\leq C C_{2/\gamma} \eta^{-2/\gamma} \mathbb E\bigl[
\sqrt{\bigl|\Pi_{(k+1)}(S_{k+1}-)}\bigr|\bigr] .
\]
Let
$S'_{k+1}=\inf\{t\geq0\dvtx\{2,3,\ldots,k+1\}\cap
\Pi_{(1)}(t)=\varnothing\}$, then by exchangeability,
\[
\limsup_{n\to\infty} \mathbb
P\bigl(d_{\mathrm{GH}}\bigl(T_{[n]}^{[k]},T_{[n]}\bigr)\geq\eta a_n\bigr)
\leq C C_{2/\gamma} \eta^{-2/\gamma} \mathbb
E\bigl[\sqrt{\bigl|\Pi_{(1)}(S'_{k+1}-)}\bigr|\bigr].
\]
Since the quantity in
the expectation goes to $0$ a.s. as $k\to\infty$ and is bounded
[indeed $S'_k\uparrow D_{\{1\}}$ a.s. and $|\Pi_{(1)}(D_{\{1\}}-)|=0$ by
(\ref{eq:1}) and Proposition~\ref{sec:self-simil-fragm-3}], we
conclude that for every $\eta>0$,
%
%e16 ###
\begin{equation}
\label{eq:25}
\lim_{k\to\infty}\limsup_{n\to\infty} \mathbb
P\bigl(d_{\mathrm{GH}}\bigl(T_{[n]}^{[k]},T_{[n]}\bigr)\geq\eta a_n
\bigr)=0 .
\end{equation}

It is now easy to get from this the convergence in distribution of
$a_n^{-1}T_{[n]}$ toward $\mathcal T_{\gamma,\nu}$ in
$(\T,d_{\mathrm{GH}})$, using the following lemma, together with
(\ref{eq:25}), Proposition~\ref{sec:conv-finite-dimens-1} and the
fact that $\mathcal{R}(\TT_{\gamma,\nu},[k])$ converges in distribution
in $(\T,d_{\mathrm{GH}})$ to $\TT_{\gamma,\nu}$ as $k\to\infty$
\cite{HM04}.
\begin{lemma}[(\cite{billingsley99}, Theorem 3.2)]
Let $X_n,X,X_n^k,X^k$ be random variables in a metric space
$(M,d)$. We assume that for every $k$, we have $X_n^k\to X^k$ in
distribution as $n\to\infty$, and $X^k\to X$ in distribution as $k
\to
\infty$. Finally, we assume that for every $\eta>0$,
\[
\lim_{k\to\infty}\limsup_{n\to\infty}\mathbb P\bigl(d(X_n^k,X_n)>\eta
\bigr)=0 .
\]
Then $X_n\to X$ in distribution as $n\to\infty$.
\end{lemma}
\begin{pf*}{Proof of Lemma~\ref{sec:tightn-grom-hausd}} Note
that if the statement holds for some $p>0$, it then holds for all $p'
\in(0,p)$. We can therefore assume in the following that $p>1/\gamma$,
and we let $\varepsilon>0$ be so that $p(\gamma-\varepsilon)>1$. The
main idea of the proof is to proceed by induction on $n$, using the
Markov branching property. We start with some technical preliminaries.

$\bullet$ First note that $\mathsf
{E}_n^q[\mathrm{ht}^r]<\infty$ for all $r>0$ and all $n \geq
1$. This can easily be proved by induction on $n$ ($r$ being
fixed) using the Markov branching property and the facts
that $q_n((n))<1$ and that $\mathrm{ht}=0$ almost-surely under~$\bP_1^q$.

$\bullet$ Second, we replace the sequence $(a_n,n\geq1)$ by
a ``nicer'' sequence $(\tilde a_n, n\geq1)$ such that $\tilde a_n \sim
a_n$, that is, $a_n/\tilde a_n \rightarrow1$ as $n \rightarrow
\infty$.
This step is trivial when $a_n = n^{\gamma}$; we then take $\tilde a_n
=a_n$. Since $a_n = n^{\gamma}\ell(n)$ with $\ell$ slowly varying at
$\infty$, it is
well known (see~\cite{BGT}, Theorem 1.3.1) that it can be written in
the form
\[
a_n = n^{\gamma}c(n)\exp\biggl(\int_1^n \eta(u) \,\mathrm du /u
\biggr),\qquad n\geq1,
\]
where $c(n) \rightarrow c>0$ as $n \rightarrow\infty$, and
$\eta$ is a measurable function that converges to
$0$ as $u \rightarrow\infty$. Define
\[
\tilde a_n = n^{\gamma}c\exp\biggl(\int_1^n \eta(u) \,\mathrm du /u
\biggr),\qquad n\geq1.
\]
We claim
that there exists an integer $n_{\varepsilon} \geq1$ such that for $n
\geq n_{\varepsilon}$,
%
%e17 ###
\begin{equation}
\label{bonnemajo}
\frac{\tilde{a}_k}{ \tilde{a}_n}\leq\biggl( \frac{k}{n}
\biggr)^{\gamma-\varepsilon} \qquad\forall1 \leq k \leq n.
\end{equation}
Indeed, let $u_{\varepsilon}$ be such that $| \eta(u) | \leq
\varepsilon$ for all $u \geq u_{\varepsilon}$. For $n \geq k \geq u_
{\varepsilon}$, we have
\[
\biggl|\int_k^n -\eta(u) \,\mathrm du /u \biggr| \leq
\varepsilon\int_k^n \mathrm du /u \leq\varepsilon\ln(n/k);
\]
hence
\[
\frac{\tilde{a}_k
}{\tilde{a}_n
} =\biggl(
\frac{k}{n} \biggr)^{\gamma} \exp\biggl(-\int_k^n \eta(u)
\,\mathrm du /u \biggr) \leq\biggl( \frac{k}{n}
\biggr)^{\gamma-\varepsilon}.
\]
Besides $\sup_{k \in\{1,\ldots, \lfloor u_{\varepsilon} \rfloor\}}
\tilde a_k k^{\varepsilon-\gamma}
/(\tilde a_n n^{\varepsilon-\gamma}) \leq1$ for all $n$
large enough (say $n \geq n'_{\varepsilon}$). Hence $
\tilde{a}_k/\tilde{a}_n \leq( k/n )^{\gamma
-\varepsilon} $
for all $n \geq n_{\varepsilon}=\max(n'_{\varepsilon},
u_{\varepsilon})$ and all $1 \leq k \leq n$.

$\bullet$ Since $a_n >0$ for all $n \geq1$
and $\tilde a_n \sim a_n$, there exists some $C>0$ such that $
a_n \geq C \tilde a_n$ for all $n \geq
1$. It is therefore sufficient to prove the existence of a finite
$C_p$ (a priori different from the one in the statement of the lemma)
such that
%
%e18 ###
\begin{equation}
\label{eqtilde}
\bP^q_n(\mathrm{ht} \geq x \tilde{a}_n
) \leq
\frac{C_p}{x^p}\qquad \forall x>0 \mbox{ and } n \geq1,
\end{equation}
to finish the proof of the lemma. In order to prove~(\ref{eqtilde}),
we will use the integer $n_{\varepsilon}$ introduced around
(\ref{bonnemajo}), and we will further assume, taking $n_{\varepsilon}$
larger if necessary, that $\tilde{a}_n \geq1$ for every $ n \geq
n_{\varepsilon}$. Introduce now $0<C^1_p<1$ such that
\[
(1-u)^{-p} \leq1+2pu\qquad \forall0 \leq u \leq C^1_p.
\]
Using (H$'$) and the fact that $q_n((n))<1$ for all $n \geq1$,
there exists also $C^2_p>0$ such that
%
%e19 ###
\begin{equation}
\label{C2}
\tilde{a}_n
\sum_{\lambda\in\mathcal P_{n}}
q_n(\lambda) \Biggl( 1-\sum_{i=1}^{p(\lambda)} \biggl( \frac
{\lambda_i}{n}\biggr)^{(\gamma-\varepsilon) p} \Biggr) \geq C^2_p\qquad
\forall n \geq1
\end{equation}
[recall that $(\gamma-\varepsilon)p>1$ and $\tilde a_n>0$ for all
$n \geq1$]. Last we let
\[
C_p(n_{\varepsilon}):=\max_{1 \leq n \leq
n_{\varepsilon}}
(\mathsf{E}_n^q[\mathrm{ht}^p]/\tilde{a}_n^p)<\infty,
\]
and
we set
\[
C_p:=\max
\bigl(C_p(n_{\varepsilon}),(1/C^1_p)^{p},
(2p/C^2_p)^p
\bigr)<\infty.
\]
Our goal is to prove by induction on $n \geq1$ that
{\renewcommand{\theequation}{$\mathrm{A}_n$}
\begin{equation}\label{An}\hypertarget{AN}
\bP_n^q(\mathrm{ht} < x \tilde{a}_n) \geq
1-\frac{C_p}{x^p}\qquad \mbox{for every } x>0.
\end{equation}}

\vspace*{-8pt}

\noindent Clearly,~(\ref{An}) holds for all $n \leq n_{\varepsilon}$
since $C_p \geq C_p(n_{\varepsilon})$ and $\bP_n^q (\mathrm{ht} \geq x
\tilde{a}_n ) \leq\mathsf{E}_n^q[\mathrm{ht}^p]/(x\tilde{a}_n )^p$.
Now assume that \hyperlink{AN}{($\mathrm{A}_k$)}
is satisfied for all $k \leq n-1$ for
some $n \geq n_{\varepsilon}$. For all $0 < x \leq C_p^{1/p}$, the
expected inequalities in~(\ref{An}) are obvious, so it remains
to prove them for $x> C_p^{1/p}$. To get~(\ref{An}), we will
prove by induction on $i \in\mathbb N$ that
{\renewcommand{\theequation}{$\mathrm{A}_{n,i}$}
\begin{equation}\label{Ani}\hypertarget{ANI}
\bP_n^q(\mathrm{ht} < x \tilde{a}_n )
\geq
1-\frac{C_p}{x^p} \qquad \mbox{for every } x \in
\biggl(0,\frac{i}{\tilde{a}_n } \biggr) ,
\end{equation}}

\vspace*{-8pt}

\noindent which will obviously lead to~(\ref{An}). Note first that
\hyperlink{ANI}{($\mathrm A_{n,1}$)}
holds since $1/\tilde{a}_n \leq1 \leq
C_p^{1/p}$. Assume next that~(\ref{Ani}) is true, and fix $x
\in(0,(i+1)/\tilde{a}_n )$. We can assume that $x >
C_p^{1/p}$ since~(\ref{An}) holds otherwise. Using the Markov
branching property and
the fact that \hyperlink{AN}{($\mathrm{A}_k$)}
holds for every $k \leq n-1$, as well
as~(\ref{Ani}), we then get
\begin{eqnarray*}
\bP_n^q(\mathrm{ht} < x \tilde{a}_n )&=& \sum_{\lambda\in\mathcal
P_n} q_n(\lambda) \prod_ {i=1}^{p(\lambda)}
\bP_{\lambda_i}^q(\mathrm{ht} < x \tilde{a}_n -1) \\ &\geq&
\sum_{\lambda\in\mathcal P_n} q_n(\lambda) \prod_{i=1}^{p(\lambda)}
\biggl( 1-\frac{C_p \tilde{a}_{\lambda_i}^p}{(x\tilde{a}_n
-1)^p}\biggr)^+ \\ &\geq& \sum_{\lambda\in\mathcal P_n}
q_n(\lambda) \Biggl( 1- \sum_{i=1}^{p(\lambda)}\frac{C_p
\tilde{a}_{\lambda_i}^p}{(x\tilde{a}_n-1)^p} \Biggr)
\end{eqnarray*}
using the notation $r^+=\max(r,0)$ and that for all sequences of
nonnegative terms $b_i, i \geq1$, $\prod_{i=1}^m (1-b_i)^+ \geq1-
\sum_{i=1}^m b_i$, for every $m \geq1$. Next, since $x\tilde a_n\geq
x >1/C_p^1>1$,
\[
\frac{1}{(x\tilde{a}_n-1)^p}=\frac{1}{(x\tilde{a}_n)^p
(1-1/(x\tilde{a}_n))^p} \leq\frac{1+2p/
(x\tilde{a}_n)}{(x\tilde{a}_n)^p},
\]
and then
\begin{eqnarray*}
&& \bP_n^q(\mathrm{ht} < x \tilde{a}_n) \\
&&\qquad\geq\sum_{\lambda
\in
\mathcal P_n} q_n(\lambda) - \frac{C_p}{x^{p}}\sum_{\lambda\in
\mathcal P_n} q_n(\lambda) \sum_{i=1} ^{p(\lambda)}
\biggl(\frac{\tilde{a}_{\lambda_i}}{\tilde{a}_n}\biggr)^p\\
&&\qquad\quad{}-\frac{2pC_p}{x^{p+1}\tilde{a}_n}\sum_ {\lambda\in
\mathcal P_n} q_n(\lambda)
\sum_{i=1}^{p(\lambda)}\biggl(\frac{\tilde{a}_{\lambda_i}}{\tilde
{a}_n}\biggr)^p
\\
&&\hspace*{-8pt}\qquad\mathop{\geq}_{\mathrm{by}\mbox{ }
\mbox{\fontsize{8.36}{8.36}\selectfont{(\ref{bonnemajo})}}} 1 -
\frac{C_p}{x^{p}}\sum_{\lambda\in\mathcal P_n} q_n(\lambda)
\sum_{i=1} ^{p(\lambda)} \biggl(\frac{\lambda_i
}{n}\biggr)^{(\gamma-\varepsilon) p}\\
&&\qquad\quad{}
-\frac{2pC_p}{x^{p+1}\tilde{a}_n}\sum_ {\lambda
\in\mathcal P_n} q_n(\lambda) \sum_{i=1}^{p(\lambda)}
\biggl(\frac{\lambda_i}{n}\biggr)^{(\gamma-\varepsilon) p}
\\
&&\qquad\geq1-\frac{C_p}{x^p} + \frac{C_p}{x^{p}}\sum_{\lambda
\in\mathcal P_n} q_n(\lambda) \Biggl( 1-
\sum_{i=1}^{p(\lambda)}
\biggl(\frac{\lambda_i
}{n}\biggr)^{(\gamma-\varepsilon) p}
\Biggr)\\
&&\qquad\quad{}
-\frac{2pC_p}{x^{p+1}\tilde{a}_n}
\sum_{\lambda\in\mathcal P_n} q_n(\lambda)
\sum_{i=1}^{p(\lambda)}\biggl(\frac{\lambda_i
}{n}\biggr)^{(\gamma-\varepsilon) p} .
\end{eqnarray*}
We then use~(\ref{C2}) and the fact that
$\sum_{i=1}^{p(\lambda)}(\lambda_i/n)^{(\gamma-\varepsilon) p}
\leq1$
[since $(\gamma-\varepsilon) p>1$] to get
\[
\bP_n^q(\mathrm{ht} < x \tilde a_n ) \geq
1-\frac{C_p}{x^p} + \frac{C_p}{x^{p}\tilde{a}_n} \biggl(
C^2_p-\frac{2p}x \biggr).
\]
By assumption, $x >C_p^{1/p} \geq2p/C^2_p$; hence
\[
\bP_n^q(\mathrm{ht} < x \tilde{a}_n) \geq1-\frac{C_p}{x^p}
\qquad\mbox{for every } x \in\biggl(0,\frac{i+1}{\tilde{a}_n} \biggr)
\]
as wanted.
\end{pf*}

%s4.4 ###
\subsection{Incorporating the measure}\label{sec:incorp-meas}

We now finish the proof of Theorem~\ref{sec:main-result-1}, by
improving the Gromov--Hausdorff convergence of Proposition
\ref{sec:tightn-grom-hausd-1} to a Gromov--Hausdorff--Prokhorov
convergence, when the uniform measure $\mu_n=\mu_{\partial T_n}$ on
leaves is added to $T_n$ in order to view it as an element of $\Tw$
rather than $\T$.\vadjust{\goodbreak}

We will use the fact (\cite{EW}, Lemma 2.3), that the convergence in
distribution of $a_n^{-1}T_n$ as $n\to\infty$ in $\T$ entails that the
laws of the random variables $a_n^{-1}T_n$ form a tight sequence of
probability measures on $\Tw$. Therefore, it suffices to identify the
limit as $\TT_{\gamma,\nu}$.

So let us assume that $a_n^{-1}T_n$ converges to
$(\TT',d',\rho',\mu')\in\Tw$ in distribution, when $n\to\infty$ along
some subsequence. Let
%. Let $I_1^n,I_2^n,\ldots,I_k^n$ be
%independent, uniform in $[n]$. These numbers are associated to leaves
%of $T_n$, say
$L_1^n,L_2^n,\ldots,L_k^n$ be $k$ i.i.d. uniform leaves of $T_n$.
Conditionally given the event that these leaves are pairwise distinct,
which occurs with probability going to $1$ as $n\to\infty$ with $k$
fixed, these leaves are just a uniform sample of $k$ distinct leaves
of $T_n$, so by Lemma~\ref{sec:mark-branch-trees-2} and
exchangeability, the subtree of $T_n$ spanned by the root and the
leaves $L_1^n,\ldots,L_k^n$ has same distribution as
$T_{[n]}^{[k]}$. By Proposition~\ref{sec:conv-finite-dimens-1}, we
know that $a_n^{-1}T_{[n]}^{[k]}$ converges in distribution to
$\mathcal{R}(\TT_{\gamma,\nu},[k])$ in $\T$.

A $k$-rooted compact metric space is an object of the form
$((X,d),x_1,\ldots,x_k)$ where $(X,d)$ is a compact metric space and
$x_1,\ldots,x_k\in X$. The set of $k$-rooted metric spaces can be
endowed with the $k$-rooted Gromov--Hausdorff distance
\begin{eqnarray*}
&&
d_{\mathrm{GH}}^{(k)}(((X,d),x_1,\ldots,x_k),
((X',d'),x'_1,\ldots,x'_k)) \\
&&\qquad=\inf_{\phi,\phi'}\max
_{1\leq
i\leq k}\dist(\phi(x_i),\phi'(x'_i))\vee
\dist_{\mathrm{H}}(\phi(X),\phi'(X')) ,
\end{eqnarray*}
where, as in the
definition of the Gromov--Hausdorff distance, the infimum is over
isometric embeddings $\phi,\phi'$ of $X,X'$ into some common space
$(M,\dist)$. Note in particular that
$d_{\mathrm{GH}}^{(1)}=d_{\mathrm{GH}}$. Now, the fact that
$a_n^{-1}T_n$ converges to $(\TT',d',\rho',\mu')$ in $\Tw$ implies
that the $k+1$-rooted space $(a_n^{-1}T_n,\rho,L_1^n,\ldots,L_k^n)$
converges in distribution to $(\TT',\rho',L_1,\ldots,L_k)$, where
$L_1,\ldots,L_k$ are i.i.d. with law $\mu'$ conditionally on the
latter. See~\cite{miermont09}, Proposition 10, for a proof and further
properties of the $k$-rooted Gromov--Hausdorff distance, which is
separable and complete.

If $(T,d,\rho)$ is a rooted $\R$-tree and $x_1,\ldots,x_k\in T$, the
union of geodesics from $\rho$ to the $x_i$'s,
\[
R(T,x_1,\ldots,x_k)=\bigcup_{i=1}^k\lbr \rho,x_i\rbr
\]
is in turn an
$\R$-tree rooted at $\rho$ with at most $k$ leaves, called the subtree
of $T$ spanned by $x_1,\ldots,x_k$ (the role of the root being
implicit).
\begin{lemma}\label{sec:incorp-meas-1}
Let $(\mathcal{A}_n,d_n,\rho_n),n\geq1$, be a sequence of rooted $\R
$-trees and
$x^n_1,\ldots,x^n_k$ be $k$ points in $\mathcal{A}_n$, such that
$((\mathcal{A}_n,d_n),\rho_n,x^n_1,\ldots,x_n^k)$ converges for~the
$k+1$-rooted Gromov--Hausdorff distance to a limit
$((\mathcal{A},d),\rho,x_1,\ldots,\break x_k)$. Then the subtree
$R(\mathcal{A}_n,x_1^n,\ldots,x_k^n)$ converges in $\mathscr{T}$ to
the subtree
$R(\mathcal{A},\break x_1,\ldots,x_k)$.
\end{lemma}

We will prove this lemma at the end of the section. By using the
Skorokhod representation theorem, we may assume that the convergence
of $(a_n^{-1}T_n,\rho,L_1^n,\ldots, L_k^n)$ to
$(\TT',\rho',L_1,\ldots,L_k)$ holds almost surely. This, together
with Lemma~\ref{sec:incorp-meas-1} and the discussion at the beginning
of this section, implies the joint convergence in distribution in $\T$
of $a_n^{-1}T_{[n]},a_n^{-1}T_{[n]}^{[k]}$ to $\TT',\TT'_k$, still
along the appropriate subsequence, where $\TT'_k$ is the subtree of
$\TT'$ spanned by $L_1,\ldots,L_k$. In particular, this identifies
the law of $\TT'_k$ as that of $\mathcal{R}(\TT_{\gamma,\nu},[k])$.
When $k\to\infty$, we already stressed that the latter trees converge
(in distribution in $\Tw$, with the uniform measure $\mu_k$ on the set
of its $k$ leaves) to $\TT_{\gamma,\nu}$. On the other hand, $\TT'_k$
converges a.s. to $(\TT'',d',\rho',\mu')$ in $\Tw$ as $k\to\infty$,
where $\TT''$ is the closure in $\TT'$ of
\[
\bigcup_{i=1}^\infty\lbr \rho',L_i\rbr  .
\]
But the joint convergence of
$T_{[n]},T_{[n]}^{[k]}$ in $\T$ along some subsequence and
%If this were not the case, we could find a leaf in $\TT'$ which is at
%positive distance from $\TT''$. This would mean that there is a leaf
%$x$ of $\TT'$ at positive distance from the support of $\mu'$. This
%would entail the existence of $\eps,\eta>0$ such that $$P(\min_{1\leq
%i\leq k}d'(x,L_i)>\eps\mbox{ for every }k\geq1)\geq\eps,$$
%which would contradict
(\ref{eq:25}) imply that for every $\eta>0$,
$\lim_{k\to\infty}\mathbb{P}(d_{\mathrm{GH}}(\TT'_k,\TT')>\eta)=0$.
So $\TT''=\TT'$ a.s., entailing that $(\TT',d',\rho',\mu')$ has
same law as $\TT_{\gamma,\nu}$. This identifies the limit of
$a_n^{-1}T_n$ in $\Tw$ as $\TT_{\gamma,\nu}$, ending the proof of
Theorem~\ref{sec:main-result-1}.

It remains to prove Lemma~\ref{sec:incorp-meas-1}. We only sketch the
argument, leaving the details to the reader. We use induction on
$k$. For $k=1$, the subtree $R(\mathcal{A}_n,x^n_1)$ is isometric to a real
segment $[0,d_n(\rho_n,x^n_1)]$ rooted at $0$. The $2$-rooted
convergence of $((\mathcal{A}_n,d_n),\rho_n,x^n_1)$ to $((\mathcal
{A},d),\rho,x_1)$
entails that $d_n(\rho_n,x^n_1)$ converges to $d(\rho,x_1)$, hence
that $R(\mathcal{A}_n,x^n_1)$ converges to $[0,d(\rho,x_1)]$ rooted
at $0$,
which is isometric to $R(\mathcal{A},x_1)$.

For the induction step, we use the general fact that if $\mathcal{A}$
is a
rooted $\R$-tree and
$x_1,\ldots,x_k,x_{k+1}\in\mathcal{A}$, then the distance between
$x_{k+1}$ and
the subtree of $\mathcal{A}$ spanned by $x_1,\ldots,x_k$ is equal to
\[
\delta_{k+1}=\min_{1\leq i\leq
k}\biggl(\frac{d(x_{k+1},x_i)+d(x_{k+1},\rho)-d(x_i,\rho)}{2}\biggr).
\]
Moreover, if $i\in\{1,2,\ldots,k\}$ is an index that realizes
this minimum, then the branchpoint $x_{k+1}\wedge x_i$ is at distance
$\delta_{k+1}$ from $x_{k+1}$ and is the ancestor of $x_i$ at height
(i.e., distance from $\rho$)
\[
h_{k+1}=d(\rho,x_{k+1})-\delta_{k+1} .
\]
In words, we get that $R(\mathcal{A},x_1,\ldots,x_{k+1})$ is obtained
from $R(\mathcal{A},x_1,\ldots,x_k)$ by grafting a segment with length
$\delta_{k+1}$ at the ancestor of $x_i$ with height $h_{k+1}$.

In our particular situation, and with some obvious notation, we get that
$R(\mathcal{A}_n,x^n_1,\ldots,x^n_{k+1})$ is obtained by grafting a
segment with
length $\delta_{k+1}^n$ to the ancestor of $x^n_{i_n}$ with height
$h_{k+1}^n$, where $i_n$ is some index in $\{1,\ldots,k\}$ that can
depend on $n$. Taking a subsequence if necessary, we may assume that
$i_n=i$ is
constant. The $k+2$-rooted convergence of
$((\mathcal{A}_n,d_n),\rho_n,x^n_1,\ldots,x^n_{k+1})$ to
$((\mathcal{A},d),\rho,x_1,\ldots,x_{k+1})$ entails that the $d_n$-distances
between elements of $\{\rho_n,x^n_1,\ldots,x^n_{k+1}\}$ converge to
the corresponding $d$-distances of elements in
$\{\rho,x_1,\ldots,x_{k+1}\}$. Consequently, it holds that
$\delta_{k+1}^n,h_{k+1}^n$ converge to $\delta_{k+1},h_{k+1}$, defined
as above. Together with the induction hypothesis stating that
$R(\mathcal{A}_n,\break x^n_1,\ldots,x^n_k)$ converges in $\mathscr{T}$ to
$R(\mathcal{A},x_1,\ldots,x_k)$, this entails easily that
$R(\mathcal{A}_n,\break x^n_1,\ldots, x^n_{k+1})$ converges in $\mathscr{T}$
to the
$\R$-tree obtained by grafting a segment with length $\delta_{k+1}$ to
the ancestor of $x_i$ with height $h_{k+1}$ in $R(\mathcal
{A},x_1,\ldots,x_k)$,
and this tree is $R(\mathcal{A},x_1,\ldots,x_{k+1})$. The result being
independent of the particular value of $i$ (selected by the choice of
a subsequence), the convergence holds without taking subsequences,
which concludes the proof.

%s4.5 ###
\subsection{\texorpdfstring{Proof of Theorem \protect\ref{sec:main-result-2}}{Proof of Theorem 6}}
\label{proof-theor2}

%In order to connect the construction above to the one in the previous
%paragraph, w

To pass from trees with $n$ vertices (with law $\bQ_n^q$) to trees
with laws of the form $\bP_n^{q'}$, with $n$ leaves, we introduce a
transformation on trees, in which every vertex which is not a leaf is
attached to an extra ``ghost'' neighbor, which is a leaf.

Precisely, if $\mathbbl{t}$ is a plane tree, then the
modification $\mathbbl{t}^\circ$ is defined as
\[
\mathbbl{t}^\circ=\mathbbl{t}\cup\bigcup_{u=(u_1,\ldots,u_k)\in
\mathbbl{t}\setminus
\partial\mathbbl{t}}\bigl\{\bigl(u_1,\ldots,u_k,c_u(\mathbbl{t})+1\bigr)\bigr\}.
\]
If we are given a
tree rather than a plane tree, then this construction performed on any
plane representative of the tree $\mathsf{t}$ will yield plane trees
in the
same equivalence class, which we call $\mathsf{t}^\circ$. Note that
\[
\#\partial\mathsf{t}^\circ=\#\mathsf{t} .
\]

We see $\mathsf{t}^\circ$ as an element of $\Mw$ (endowed with graph
distance and uniform distribution on $\partial\mathsf{t}^\circ$),
and view
$\mathsf{t}$ as an element of $\Mw$, by endowing it also with the graph
distance, but this time, with the uniform distribution $\mu_\mathsf
{t}$ on~$\mathsf{t}$. It
is easy to see, using the natural isometric embedding of $\mathsf{t}$ into
$\mathsf{t}^\circ$, that for every $a>0$,
%
%e20 ###
\setcounter{equation}{19}
\begin{equation}\label{eq:19}
d_{\mathrm{GHP}}(a\mathsf{t},a\mathsf{t}^\circ)\leq a .
\end{equation}
Let $(q_n,n\geq1)$ be, as in Section~\ref{sec:mark-branch-trees-1}, a
family of probability distributions, respectively, on $\P_n$, such that
$q_1((1))=1$. We introduce the family $q^\circ_n,n\geq1$, of
probability measures, respectively, on $\P_n$ by
$q^\circ_1(\varnothing)=1$, and
\[
q^\circ_{n+1}((\lambda,1))=q_n(\lambda) ,\qquad n\geq1,\lambda
\in
\P_{n} ,
\]
where
$(\lambda,1)=(\lambda_1,\ldots,\lambda_{p(\lambda)},1)\in\P_{n+1}$.

It is then immediate to show by induction that if $T_n$ has law
$\bQ^q_n$, then $T^\circ_n$ has law~$\bP^{q^\circ}_n$, with the
notation of Section~\ref{sec:mark-branch-trees}. We leave this
verification to the reader. In view of this and~(\ref{eq:19}), we see
that Theorem~\ref{sec:main-result-2} is a straightforward consequence
of the following statement.
\begin{lemma}\label{sec:proof-theorem-xxx}
If $(q_n,n\geq1)$ satisfies \textup{(H)} with either
$\gamma\in(0,1)$, or $\gamma=1$ and $\ell(n)\to0$ as $n\to\infty$,
then $(q^\circ_n,n\geq1)$ satisfies \textup{(H)}, with same
fragmentation pair $(-\gamma,\nu)$ and function $\ell$.
\end{lemma}
\begin{pf}
Let $f\dvtx\sd\to\R$ be a Lipschitz function with
uniform norm and Lipschitz constant bounded by $K$. Let also
$g(\bs)=(1-s_1)f(\bs)$. Then
\[
\biggl|f\biggl(\frac{(\lambda,1)}{n+1}\biggr)-f\biggl(\frac{\lambda
}{n}\biggr)\biggr|
\leq K\max\biggl(\sup_{1\leq i\leq p(\lambda)}\frac{\lambda
_i}{n(n+1)},\frac{1}{n+1}\biggr)
\leq\frac{K}{n+1}
,
\]
so that
\begin{eqnarray*}
|\ov{q}{}^{\circ}_{n+1}(g)-\ov{q}_n(g)|&\leq& \sum_{\lambda\in
\P_n}q_n(\lambda) \biggl|
\biggl(1-\frac{\lambda_1}{n+1}\biggr)f\biggl(\frac{(\lambda
,1)}{n+1}\biggr)-
\biggl(1-\frac{\lambda_1}{n}\biggr)f\biggl(\frac{\lambda}{n}
\biggr)\biggr|\\ &\leq&\sum_{\lambda\in
\P_n}q_n(\lambda)\biggl(\frac{K\lambda_1}{n(n+1)}+\frac
{K}{n+1}\biggr)\\ &\leq
& \frac{2K}{n+1}.
\end{eqnarray*}
Multiplying both sides by $n^\gamma\ell(n)$, we see that the
upper bound converges to $0$ as $n\to\infty$ under our
hypotheses. Since $n^\gamma\ell(n)\ov{q}_n(g)$ converges to $\nu(g)$
by~(H), we obtain the same convergence with $\ov{q}{}^\circ_n$
instead of $\ov{q}_n$. This yields the result.
\end{pf}

%s4.6 ###
\subsection{\texorpdfstring{Proof of Proposition \protect\ref{Propexemple}}{Proof of Proposition 7}}

Recall the notation $\Lambda^{(\mathbf s)}(n)$ for the decreasing
sequence of sizes of blocks restricted to $\{1,\ldots, n \}$ of a
random variable with paintbox distribution $\rho_{\mathbf s}(\mathrm d
\pi)$, with $\mathbf s \in\mathcal S^{\downarrow}$, $\sum_{i \geq
1}s_i=1$. Recall also that $\Lambda^{(\mathbf s)}(n)/n \rightarrow
\mathbf s$ in $\mathcal S^{\downarrow}$ almost-surely. Now set for
$\lambda\in\mathcal P_n$,
\begin{eqnarray*}
\tilde q_n(\lambda)&=&n^{-\gamma} \int_{\mathcal S^{\downarrow}}
\mathbb P \bigl(\Lambda^
{(\mathbf s)}(n)=\lambda\bigr) \mathbf1_{\{ n^{-\gamma/2} \leq1- s_1\}}
\nu(\mathrm d \mathbf s),\qquad
\lambda\neq(n), \\
\tilde q_n((n))&=& 1-\sum_{\lambda\in\mathcal P_n, \lambda\neq(n)}
\tilde q_n(\lambda).
\end{eqnarray*}
%
%(since $\nu$ is conservative, we do not need to define $\Lambda^{(
For $n$ large enough, say $n\geq n_0$,
\[
0<\sum_{\lambda\in\mathcal
P_n, \lambda\neq(n)} \tilde q_n(\lambda) \leq n^{-\gamma/2}
\int_{\mathcal S^{\downarrow}} (1-s_1)\nu(\mathrm d \mathbf s)\leq
1,
\]
hence $\tilde q_n$ defines a probability distribution on
$\mathcal P_n$ such that $\tilde q_n((n))<1$. Then set $q_n=\tilde
q_n$ for $n \geq n_0$, and for $n < n_0$, let $q_n$ be any distribution
on $\mathcal P_n$ such that $q_n((n))<1$.
%Clearly $q_n((n)) <1$ since $\nu((1,0,...))=0$.

Next, consider a continuous function $f \dvtx \mathcal
S^{\downarrow} \rightarrow\mathbb R_+$. For $n \geq n_0$, we have
%(by Fubini's theorem since $f$ is bounded on the compact space
%$\mathcal S^{\downarrow}$)
%
\begin{eqnarray*}
&&n^{\gamma}\sum_{\lambda\in\mathcal P_n} q_n(\lambda)
\biggl(1-\frac{\lambda_1}{n}\biggr)f\biggl( \frac{\lambda
}{n}\biggr) \\
&&\qquad=
\int_{\mathcal S^{\downarrow}} \mathbb E\biggl[
\biggl(1-\frac{\Lambda^{(\mathbf s)}_1(n)}{n} \biggr) f\biggl(\frac
{\Lambda^
{(\mathbf s)}(n)}{n} \biggr) \biggr] \mathbf1_{\{ n^{- \gamma/2}
\leq1- s_1\}} \nu(\mathrm d \mathbf s),
\end{eqnarray*}
which converges to $\int_{\mathcal S^{\downarrow}} f(\mathbf s)
(1-s_1) \nu(\mathrm d \mathbf s) $ as $n \rightarrow\infty$ by
dominated convergence. This completes the
proof.

%%%%%%%%%%%%%%%%%%%%%%%%%%%%%%%%%%%%%%%%%%%
%s5 ###
\section{Scaling limits of conditioned Galton--Watson trees}
\label{sec:cond-galt-wats}
%%%%%%%%%%%%%%%%%%%%%%%%%%%%%%%%%%%%%%%%%%%

Recall the notation of Section~\ref{sec:appl-galt-wats}. Since the
probability distribution $\mathrm{GW}_\xi$ enjoys the so-called branching
property, it holds that the conditioned versions $\gw^{(n)}$ are
Markov branching trees.
\begin{proposition}\label{sec:main-results-3}
\textup{(i)} One has $\gw^{(n)}=\bQ_n^q$ for every $n\geq1$, where the splitting
probabilities $q=(q_n,n\geq1)$ are defined by $q_1((1))=1$, and for
every $n\geq2$ and $\lambda=(\lambda_1,\ldots,\lambda_p)\in
\P_n$,
%
%e21 ###
\begin{equation}\label{eq:11}
q_n(\lambda) = \frac{p!}{\prod_{j\geq
1}m_j(\lambda)!}\xi(p)\frac{\prod_{i=1}^p\operatorname{
GW}_\xi(\#\mathbbl{t}=\lambda_i)}{\mathrm{GW}_\xi(\#\mathbbl{t}=
n+1)} .
\end{equation}

\textup{(ii)} On some probability space $(\Omega,\FF,\PP)$, let
$X_1,X_2,\ldots$ be i.i.d. with distribution $\PP(X_1=k)=\gw(\#
\mathbbl{t}=k)$,
and set $\tau_p=X_1+X_2+\cdots+X_p$. Then
%
%e22 ###
\begin{equation}
\label{eq:32}
q_n\bigl(p(\lambda)=p\bigr)=\xi(p)\frac{\PP(\tau_p = n)}{\PP(\tau_1 =
n+1)} ,
\end{equation}
and $q_n(\cdot| \{p(\lambda)=p\})$ is the law of the
nonincreasing rearrangement of $(X_1,\ldots,\break X_p)$ conditionally on
$X_1+\cdots+X_p = n$.
\end{proposition}
\begin{pf}
(i) Under $\gw$ (viewed as a law on plane trees),
conditionally on $c_\varnothing=p$, the $p$ (plane) subtrees born from
$\varnothing$ are independent with law $\gw$. For integers
$a_1,\ldots,a_p$ with sum $n$, the probability that these trees have
sizes equal to $a_1,\ldots,a_p$ is thus
$\prod_{i=1}^p\gw(\#\mathbbl{t}=a_i)$. Hence,
%
%e23 ###
\begin{equation}\label{eq:21}
\gw^{(n+1)}(c_\varnothing=p,\#\mathbbl{t}_i=a_i,1\leq i\leq p)=
\xi(p)\frac{\prod_{i=1}^p\gw(\#\mathbbl{t}=a_i)}{\gw(\#\mathbbl {t}= n+1)} ,
\end{equation}
and conditionally
on the event on the left-hand side, the subtrees born from the root
are independent with respective laws $\gw^{(a_i)},1\leq i\leq
p$. Letting $\lambda$ be the nonincreasing rearrangement of
$(a_1,\ldots,a_p)$ and re-ordering the subtrees by nonincreasing
order of size (with some convention for ties, e.g., taking them in
order of appearance according to the plane structure), we see that
these subtrees are independent with laws $\gw^{(\lambda_i)},1\leq
i\leq p$. Using the fact that there are $p!/\prod_{j\geq
1}m_j(\lambda)!$\vspace*{1pt}\vadjust{\goodbreak} compositions $(a_1,\ldots,a_p)$ of the integer $n$
corresponding to a partition $\lambda\in\P_n$, and viewing $\gw$ as a
law on $\bT$ instead of plane trees, the conclusion easily follows.

(ii) We have $q_n(p(\lambda)=p)=\gw^{(n+1)}(c_\varnothing=p)$, and the
wanted result is just an interpretation of~(\ref{eq:21}).
\end{pf}

On the same probability space $(\Omega,\FF,\PP)$ as in the previous
statement, we will also assume that $(S_r,r\geq0)$ is a random walk
with i.i.d. steps, each having distribution $\xi(i+1),i\geq-1$.
Then the well-known Otter--Dwass formula (or \textit{cyclic lemma})
(\cite{PitmanStFl}, Chapter 6), stating that $\PP(\tau_r=m)=(r/m)\PP(S_m=-r)$
for every $r,m\geq1$, allows us to rewrite
%
%e24 ###
\begin{equation}
\label{eq:14}q_n\bigl(p(\lambda)=p\bigr)=
\xi(p)\frac{{p}\PP(S_{n}=-p)/n}{\PP(S_{n+1}=-1)/({n+1})}=
\frac{n+1}{n}\hat{\xi}(p)\frac{\PP(S_{n}=-p)}{\PP(S_{n+1}=-1)}
,\hspace*{-28pt}
\end{equation}
where
$\hat{\xi}(p)=p\xi(p)$ is the size-biased distribution associated with
$\xi$ [this is a probability distribution by \eqref{eq:9}].

It is often convenient to work with size-biased orderings
of the sequence $(X_1,\ldots,X_p)$ rather than with its nonincreasing
rearrangement. Recall that if $(y_1,y_2,\ldots)$ is a nonnegative
sequence with $\sum_iy_i<\infty$, we define its size-biased ordering
in the following way. If all terms are zero, then we let $y^*_1=0$,
otherwise we let $i^*$ be a random variable with
\[
\PP(i^*=i)=\frac{y_i}{\sum_{j\geq1}y_j}
\]
and set
$y^*_1=y_{i^*}$. We then remove the $i^*$th term from the sequence
$(y_i,i\geq1)$ and resume the procedure, defining a random
re-ordering $(y^*_1,y^*_2,\ldots)$ of the sequence
$(y_1,y_2,\ldots)$. The size-biased ordering $(Y_1^*,Y_2^*,\ldots)$ of
a random sequence $(Y_1,Y_2,\ldots)$ with finite sum is defined
similarly, by first conditioning on $(Y_1,Y_2,\ldots)$. If $\mu$ is
the law of $(Y_1,Y_2,\ldots)$, we let $\mu^*$ be the law of
$(Y_1^*,Y_2^*,\ldots)$.

If $\mu$ is a probability distribution on $\sd$, then $\mu^*$ is a
probability distribution on the set
$\mathcal{S}_1=\{\mathbf{x}=(x_1,x_2,\ldots)\in[0,1]^\N\dvtx\sum
_{i\geq1}x_i\leq
1\}$ which is endowed with any metric inducing the product
topology---in particular, $\mathcal{S}_1$ is compact. Similarly, if
$\mu$ is
a nonnegative measure on $\sd$, we let
$\mu^*(f)=\int_\sd\mu(\d\bs)\mathbb{E}[f(\bs^*)]$, for every
nonnegative measurable $f\dvtx\mathcal{S}_1\to\R_+$, where $\bs^*$ is
the size-biased reordering of~$\bs$. The following statement is a
simple variation of
\cite{BertoinBook}, Proposition 2.3, replacing probability
distributions with finite measures.
\begin{lemma}\label{sec:cond-galt-wats-1}
Let $\mu_n,n\geq1$, and $\mu$ be finite measures on $\sd$, and assume
that $\mu$ is supported on $\{\bs\in\sd\dvtx\sum_is_i=1\}$. Then $\mu_n$
converges weakly to $\mu$ if and only if $\mu_n^*$ converges weakly to
$\mu^*$.
\end{lemma}

%s5.1 ###
\subsection{Finite variance case}\label{sec:finite-variance-case}

Here we assume that $\xi$ has finite variance $\sum_{p\geq
1}p(p-1)\xi(p)=\sigma^2<\infty$. In the\vspace*{1pt} proofs to come, $C$ will
denote a positive, finite constant with values that can differ from
line to line.\vadjust{\goodbreak}

Our goal is to check hypothesis (H) for the sequence $q$ of
(\ref{eq:11}), and for the measure $\nu=(\sigma/2)\nu_2$. Due to Lemma
\ref{sec:cond-galt-wats-1}, it suffices to show that
%
%e25 ###
\begin{equation}\label{eq:27}
n^{1/2}\bigl((1-s_1)\ov{q}_n(\d\bs)\bigr)^*
\mathop{\longrightarrow}^{(\mathrm{w})}_{n\to\infty}
(\sigma/2)\bigl((1-s_1)\nu_2(\d\bs)\bigr)^* .
\end{equation}
Now, for any nonnegative measure $\mu$ on $\sd$ and any nonnegative
continuous function $f$ on $\mathcal{S}_1$, one can check that
%
%e26 ###
\begin{equation}\label{eq:26}
\bigl((1-s_1)\mu(\d\bs)\bigr)^*(f)=
\int_{\mathcal{S}_1}\mu^*(\d\mathbf{x})(1-\max\mathbf
{x})f(\mathbf{x})
,
\end{equation}
where $\max\mathbf{x}=\max_{i\geq1}x_i$.
%Note that the function $\bx\mapsto
%continuous on $\triangle_1=\{\mathbf x \in\mathcal S_1 : \sum_i
%x_i=1\}$, and $\ov{q}_n^*(\mathcal S_1 \setminus\triangle_1)=0$, as
%well as $\nu_2^*(\mathcal S_1 \setminus\triangle_1)=0$.
%Further,
%$(1-\ov{\bx})/(1-x_1)\leq1$, so that $ | (1-\ov{\bx})f(\bx) | \leq
%C(1-x_1)$ for some finite $C$.
%$\nu_2^*$-almost-every continuous. This is a consequence of the fact
%that $\nu_2^*$-almost-every $\bx$ is proper so $\bx\mapsto\ov{\bx}$
%is $\nu_2^*$-almost-everywhere continuous, and the fact that
%$\nu_2^*(\{(1,0,0,\ldots)\})=0$ by definition. Further,
%$(1-\ov{\bx})/(1-x_1)\leq1$, so the function under consideration is
%also bounded.
Applying~(\ref{eq:26}) to $\mu=\ov{q}_n$ and $\nu_2$, we
conclude that~(\ref{eq:27}) is a consequence of the following statement.
%fact that
%$n^{1/2}(1-x_1)\ov{q}_n^*(\d\bx)$ converges weakly toward
%$(\sigma/2)(1-x_1)\nu_2^*(\d\bx)$, when these measure are restricted
%to $\triangle_1$.
%
\begin{proposition}\label{sec:proof-theor-refs-2}
Let $f\dvtx\mathcal{S}_1\to\R$ be a continuous function, and let
$g(\mathbf{x})=(1-\max\mathbf{x})f(\mathbf{x})$. Then
%
%e27 ###
\begin{equation}\label{eq:17}
\sqrt{n}\ov{q}_n^*(g)\mathop{\longrightarrow}_{n\to\infty}
\frac{\sigma}{\sqrt{2\pi}}\int_0^1\frac{\d
x}{x^{1/2}(1-x)^{3/2}}g(x,1-x,0,\ldots) .
\end{equation}
\end{proposition}

In summary, Theorem
\ref{sec:main-results-1} in Case 1 is a consequence of this
statement, and Theorem~\ref{sec:main-result-2}.

Proposition~\ref{sec:proof-theor-refs-2} will be proved in a couple of
steps. A difficulty that we will have to be careful about is that
$\mathbf{x}\mapsto\max\mathbf{x}$ is not continuous on $\mathcal
{S}_1$. Fix~$f$,
as in the statement. Note that $0\leq1-\max\mathbf{x}\leq1-x_1$ for every
$\mathbf{x}\in\mathcal{S}_1$, so that $|g(\mathbf{x})|\leq
c(1-x_1)$ for every
$\mathbf{x}\in\mathcal{S}_1$ for some finite $c>0$, a fact that will be
useful. In the sequel, to simplify things, we will assume $c=1$
without loss of generality.

First, note that combining (ii) in Proposition
\ref{sec:main-results-3} with a size-biased ordering, it holds that
%
%e28 ###
\begin{equation}\label{eq:22}
\ov{q}^{}*_n(g)=\sum_{p\geq1}q_n\bigl(p(\lambda)=p\bigr)
\E\biggl[g\biggl(\frac{(X^*_1,\ldots,X^*_p,0,\ldots)}{n}\biggr)
\bigg| \tau_p = n\biggr] .
\end{equation}

\begin{lemma}\label{sec:proof-theor-refs-1}
For every $\eps>0$,
\[
\sqrt{n}q_n\bigl(p(\lambda)>\eps
\sqrt{n}\bigr)\mathop{\longrightarrow}_{n\to\infty}0 .
\]
\end{lemma}
\begin{pf}
From~(\ref{eq:14}), the local limit theorem in
the finite-variance case
%
%e29 ###
\begin{equation}\label{eq:16}
\sup_{p\in\Z}\biggl|\sqrt{n}\PP(S_n=-p)-\frac{1}{\sqrt{2\pi
\sigma^2}}
\exp\biggl(-\frac{p^2}{2n\sigma^2}\biggr)\biggr|\mathop
{\longrightarrow}_{n\to\infty}0
\end{equation}
shows that $q_n(p(\lambda)=p)\leq C \hat{\xi}(p)$ for every $n,p$. Now
\[
\sum_{k\geq0}\hat{\xi}((k,\infty))<\infty,
\]
because
$\hat{\xi}$ has finite mean. Since $\hat{\xi}((k,\infty))$ is
nonincreasing, this entails that
$\hat{\xi}((k,\infty))=o(k^{-1})$. Hence the result.
\end{pf}
\begin{lemma}\label{sec:proof-theor-refs-3}
One has
\[
\lim_{\eta\downarrow
0}\limsup_{n\to\infty}\sqrt{n}\ov{q}_n^*\bigl(|g|\ind_{\{x_1>1-\eta\}
}\bigr)=0
\quad \mbox{and}\quad
\lim_{n\to\infty}\sqrt{n}\ov{q}_n^*\bigl(\ind_{\{x_1<n^{-7/8}\}}\bigr)=0 .
\]
\end{lemma}
\begin{pf}
Let $\eta>0$. Since $|g(\mathbf{x})|\leq
(1-x_1)$, we
get using~(\ref{eq:22}), that\break
$\sqrt{n}\ov{q}_n^*(|g|\times\ind_{\{x_1>1-\eta\}})$ is bounded from above
by
\begin{eqnarray*}
&&n^{1/2}\sum_{p\geq1
%1\leq p\leq\eps n^{1/2}
}q_n\bigl(p(\lambda)=p\bigr)\\
&&\hspace*{34.8pt}{}\times\sum_{
(1-\eta)n\leq m_1\leq n}\biggl(1-\frac{m_1}{n}\biggr)
\frac{pm_1}{n}\frac{\PP(X_1=m_1)\PP(\tau_{p-1} = n-m_1)}{\PP(\tau
_p =
n)}
% +o(1)
,
\end{eqnarray*}
where we used the fact (left as an exercise to the reader) that
\[
\PP(X_1^*=m | X_1+\cdots+X_p = n)=\frac{pm}{n}
\frac{\PP(X_1=m)\PP(X_2+\cdots+X_p = n-m)}{\PP(X_1+\cdots+X_p =
n)} .
\]
%
%The $o(1)$ term accounts for the fact that we restricted the sum to
%$1\leq p\leq\eps n^{1/2}$, which costs at most $o(n^{-1/2})$ by Lemma
%Using this for $p=1$,
%and with the help of~(\ref{eq:16}), we obtain that $\PP(\tau_1 = n)\sim
%(\sigma\sqrt{2\pi})^{-1}n^{-3/2}$ as $n\to\infty$. This
By replacing $q_n(p(\lambda)=p)$ by its value \eqref{eq:32}, and using
the cyclic lemma, we
obtain the upper bound
\begin{eqnarray*}
&&\sqrt{n}\ov{q}_n^*\bigl(|g|\ind_{\{x_1>1-\eta\}}\bigr)\\
&&\qquad\leq\frac{n+1}{n}\sum_{p\geq
1}p(p-1)\xi(p)\sum_{(1-\eta)n\leq m_1\leq
n}\frac{\PP(S_{m_1}=-1)\PP(S_{n-m_1}=-p+1)}{\sqrt{n}\PP(S_{n+1}=-1)}
%+o(1)
.
\end{eqnarray*}
%
%Since $1\leq p\leq\eps n^{1/2}$,
Now,~(\ref{eq:16}) implies that $\sqrt{m_1}\PP(S_{m_1}=-1)$ and
$\sqrt{n-m_1}\PP(S_{n-m_1}=-p+1)$ are bounded from above by positive
constants that are independent of $n,m_1,p$,
%(but may depend on $\eps$
%in the last case),
while $\sqrt{n}\PP(S_{n+1}=-1)$ converges to a
positive limit. Consequently, the bound is
\[
C
%_\eps
\sum_{p\geq1}p(p-1)\xi(p)\frac{1}{n}\sum_{(1-\eta)n< m_1\leq
n}\frac{1}{\sqrt{({m_1}/{n})(1-{m_1}/{n})}}
,
\]
and this converges to
$C
% _\eps
\sigma^2\int_{1-\eta}^1(x(1-x))^{-1/2}\,\d x$ as $n\to\infty$. In
turn, this goes to $0$ as $\eta\to0$.
%, for fixed $\eps$.
The second
limit is obtained in a similar way, writing
$\sqrt{n}\ov{q}_n^*(\ind_{\{x_1<n^{-7/8}\}})$ as
\begin{eqnarray*}
&&
n^{1/2}\sum_{p\geq1} q_n\bigl(p(\lambda)=p\bigr)\sum_{1\leq
m_1\leq
n^{1/8}}\frac{pm_1}{n}\frac{\PP(X_1=m_1)\PP(\tau_{p-1} = n-m_1)}{
\PP(\tau_p = n)}\\
&&\qquad\leq Cn^{-1/2}\sum_{p\geq1}p(p-1)\xi(p)\sum
_{1\leq m_1\leq
n^{1/8}}\PP(S_{m_1}=-1)\frac{\PP(S_{n-m_1}=-p+1)}{\PP(S_{n+1}=-1)}
\\
&&\qquad\leq C n^{-3/8}\sum_{p\geq1}p(p-1)\xi(p)
\end{eqnarray*}
for some finite constant $C$, where we used the local limit theorem in
the last step, and bounded $\PP(S_{m_1}=-1)$ by $1$ (we could also
bound it by $Cm_1^{-1/2}$ and obtain a better bound, but we do not
need it).
\end{pf}
\begin{lemma}\label{sec:proof-theor-refs-4}
For every $\eta>0$, it holds that
\[
\lim_{n\to\infty}\sqrt{n}\ov{q}{}^*_n\bigl(\ind_{\{x_1+x_2<1-\eta\}
}\bigr)=0.
\]
\end{lemma}
\begin{pf}
Fix $\eps>0$. Then by \eqref{eq:22},
\begin{eqnarray*}
&&\sqrt{n}\ov{q}{}^*_n\bigl(\ind_{\{x_1+x_2<1-\eta\}}\bigr)\\
&&\qquad\leq\sqrt{n}\sum_{2\leq p\leq\eps n^{1/2}}q_n\bigl(p(\lambda)=p\bigr)
\PP\bigl(X_1^*+X_2^*<(1-\eta)n|\tau_p = n\bigr)\\
&&\qquad\quad{}+\sqrt{n}q_n\bigl(p(\lambda)> \eps
\sqrt{n}\bigr)
\end{eqnarray*}
as one can note that the $p=1$ term in the sum is zero. The last
quantity being $o(1)$ by Lemma~\ref{sec:proof-theor-refs-1}, we only
have to concentrate on the first term of the right-hand side. Using
\eqref{eq:32}, we rewrite it as
%
%e30 ###
\begin{equation}\label{eq:33}
\sqrt{n}\sum_{2\leq p\leq\eps n^{1/2}}
\xi(p) \sum_{m_1+m_2<
(1-\eta)n}\frac{\PP(X_1^*=m_1,X_2^*=m_2,\tau_p = n)}{\PP(\tau_1 =
n+1)}
.
\end{equation}
Next, using the fact (also left to the reader) that
\begin{eqnarray*}
&&\PP(X_1^*=m_1,X_2^*=m_2,\tau_p = n)\\
&&\qquad=\frac{pm_1}{n}\PP(X_1=m_1)\frac{(p-1)m_2}{n-m_1}\PP(X_2=m_2)\\
&&\qquad\quad{}\times\PP
(X_3+\cdots
+X_p = n-m_1-m_2) ,
\end{eqnarray*}
and using the cyclic lemma, we can bound \eqref{eq:33} by
\begin{eqnarray*}
&&\frac{(n+1)\sqrt{n}}{n}
\sum_{2\leq p\leq\eps n^{1/2}}
p^3\xi(p) \\
&&\qquad{}\times\sum_{m_1+m_2< (1-\eta)n}
\frac{\PP(S_{m_1}=-1)\PP(S_{m_2}=-1)\PP
(S_{n-m_1-m_2}=-p+2)}{(n-m_1)(n-m_1-m_2)\PP(S_{n+1}=-1)}.
\end{eqnarray*}
If $m_1+m_2< n(1-\eta)$, then $n-m_1\geq n-m_1-m_2\geq\eta n$. In this
case, the local limit theorem~(\ref{eq:16}) implies
\[
\frac{\PP(S_{m_1}=-1)\PP(S_{m_2}=-1)\PP
(S_{n-m_1-m_2}=-p+2)}{(n-m_1)(n-m_1-m_2)\PP(S_{n+1}=-1)}\leq
\frac{C}{\eta^{5/2}n^2\sqrt{m_1m_2}} .
\]
Note that the constant
$C$ here does not
depend on $p,\eps$. Consequently, we obtain the bound
\begin{eqnarray*}
&&\sqrt{n}\ov{q}{}^*_n\bigl(\ind_{\{x_1+x_2<(1-\eta)n\}}\bigr)\\[-2pt]
&&\qquad\leq
\frac{C}{\eta^{5/2}\sqrt{n}}\sum_{2\leq p\leq\eps
n^{1/2}}p^3\xi(p)\frac{1}{n^2}\sum_{m_1+m_2<(1-\eta)n}
\sqrt{\frac{n}{m_1}\frac{n}{m_2}} +o(1)\\[-2pt]
&&\qquad\leq
\frac{C\eps}{\eta^{5/2}}\sum_{p\geq1}p^2\xi(p)\int_{x_1+x_2\leq
1}\frac{\d x_1\,\d x_2}{\sqrt{x_1x_2}}+o(1) ,
\end{eqnarray*}
where $C$ is still independent of $p,\eps$. The first term on the
right-hand side is finite and does not depend on $n$ anymore, and it goes
to $0$ as $\eps\to0$, entailing the result.
\end{pf}
\begin{lemma}\label{sec:proof-theor-refs-5}
There exists a function $\beta_\eta=o(\eta)$ as
$\eta\downarrow0$, so that
%For every $\eps>0$, there exists some $\eta_0>0$ such that for every
%$\eta\in(0,\eta_0)$, there exists $\eta'>0$ such that for every $n$
%large enough,
%
\begin{eqnarray*}
&&
\lim_{\eta\downarrow0}\liminf_{n\to\infty}\sqrt{n}
\ov{q}{}^*_n\bigl(g\ind_{\{x_1<1-\eta,x_1+x_2>1-\beta_\eta\}}\bigr) \\[-2pt]
&&\qquad=\lim_{
\eta\downarrow0}\limsup_{n\to\infty}\sqrt{n}
\ov{q}{}^*_n\bigl(g\ind_{\{x_1<1-\eta,x_1+x_2>1-\beta_\eta\}}\bigr)\\[-2pt]
&&\qquad=
\frac{\sigma}{\sqrt{2\pi}}
\int_0^1\frac{g((x,1-x,0,\ldots))}{x^{1/2}(1-x)^{3/2}}\,\d x .
\end{eqnarray*}
\end{lemma}
\begin{pf}
The proof is similar to the previous ones, but
technically more tedious, so we will only sketch the details.
% First, we note that $\sqrt{n}\ov{q}_n^*(\{\bx:x_1<n^{-3/4}\})\leq
%n^{1/2}n^{-3/4}=o(1)$.
Fix $\eta>0$, and consider $\eta'\in(0,\eta)$ and $\eps>0$. Then, by
decomposing with respect to the events $\{p(\lambda)>\eps\sqrt{n}\}$
and $\{\mathbf{x}\dvtx x_1\leq n^{-7/8}\}$, we obtain, using Lemma
\ref{sec:proof-theor-refs-1} and the second limit of Lemma
\ref{sec:proof-theor-refs-3},
\begin{eqnarray*}
\hspace*{-4pt}&&\sqrt{n}\ov{q}{}^*_n\bigl(g\ind_{\{x_1<1-\eta,x_1+x_2>1-\eta'\}}\bigr)\\[-2pt]
\hspace*{-4pt}&&\qquad=o(1)\\[-2pt]
\hspace*{-4pt}&&\qquad\quad{}+
\sqrt{n}\sum_{2\leq p\leq\eps n^{1/2}}q_n\bigl(p(\lambda)=p\bigr)\\[-2pt]
\hspace*{-4pt}&&\qquad\quad\hspace*{67pt}{}\times
\mathop{\sum_{n^{1/8}\leq m_1\leq(1-\eta)n}}_{(1-\eta')n\leq
m_1+m_2\leq n}\hspace*{-2pt}\E\bigl[g\bigl((m_1,m_2,X_3^*,\ldots,X^*_p,0,\ldots)/n\bigr)
|\\[-2pt]
\hspace*{-4pt}&&\hspace*{-2pt}\qquad\quad\hspace*{200pt}\tau_p =
n,X_1^*=m_1,X^*_2=m_2\bigr]\\[-2pt]
\hspace*{-4pt}&&\hspace*{-4.6pt}\qquad\quad\hspace*{151.4pt}{}\times
\frac{pm_1}{n}\frac{(p-1)m_2/n}{1-m_1/n}\PP(X_1=m_1)\\[-2pt]
\hspace*{-4pt}&&\hspace*{-4.6pt}\qquad\quad\hspace*{151.4pt}{}\times\PP(X_2=m_2)\frac{\PP(\tau_{p-2} = n-m_1-m_2)}{\PP(\tau_p = n)} .
\end{eqnarray*}
We now give a lower bound of the $\liminf$ of this as
$n\to\infty$. Obtaining an appropriate upper bound for the $\limsup$
is similar
and easier.\vadjust{\goodbreak}

Note that if $x_1+x_2\geq1-\eta'$ and $x_1\leq1-\eta$, we have that
$(1-x_1-x_2)/(1-x_1) \leq\eta'/\eta$, and then $x_2/(1-x_1) \geq
1-\eta'/\eta$. Next, by \eqref{eq:16}, we can always choose
$\eta'$ small enough so that $\PP(X_1=m_2)/\PP(X_1 = n-m_1)\geq
1-\eta$ for every $n$ large enough, where $m_1,m_2$ are as in the
above sum.

Also, still by \eqref{eq:16}, and using \eqref{eq:32}, we can choose
$\eps$ small
enough so that for every $1\leq p\leq\eps n^{1/2}$ and every $n$
large, we have
\[
q_n\bigl(p(\lambda)=p\bigr)/\hat{\xi}(p)\geq(1-\eta)
\quad\mbox{and}\quad
\bigl(p^{-1}n^{3/2}\PP(\tau_p = n)\bigr)^{-1}\geq(1-\eta)\sigma
\sqrt{2\pi} .
\]
A
third use of \eqref{eq:16} entails that
\[
m_1^{3/2}\PP(X_1=m_1)\wedge m_2^{3/2}\PP(X_2=m_2)\geq
(1-\eta)/\sigma\sqrt{2\pi}
\]
for every $n$ large and $m_1\geq
n^{1/8}$, $m_2\geq(\eta-\eta')n$.

Finally, we use the fact that $f$ is uniformly continuous on
$\mathcal{S}_1$, while $\max\mathbf{x}= x_1\vee x_2$ on the set $\{
\mathbf{x}\in
\mathcal{S}_1\dvtx x_1+x_2>3/4\}$. Consequently, the function
$g(\mathbf{x})=(1-\max\mathbf{x})f(\mathbf{x})$ is uniformly
continuous on the latter set.
Therefore, we can choose $\eta'<1/4$ small enough so that
\[
\bigl|g\bigl((m_1,m_2,m_3,\ldots)/n\bigr)- g\bigl((m_1,n-m_1,0,\ldots)/n\bigr)\bigr|\leq\eta
\]
for every $(m_1,m_2,\ldots)$ with sum $n$, such that $m_1+m_2\geq
(1-\eta')n$. Putting things together, for every $\eta>0$, we can
choose $\eta'=:\beta_\eta,\eps$ small so that for every $n$ large
enough, $\sqrt{n}\ov{q}{}^*_n(g\ind_{\{x_1<1-\eta,x_1+x_2>1-\eta'\}})$
is greater than or equal to\looseness=-1
\begin{eqnarray*}
\hspace*{-4pt}&&(1-\eta)^5(1-\eta'/\eta)\\
\hspace*{-4pt}&&\hspace*{-2pt}\qquad{}\times\sum_{2\leq p\leq\eps
n^{1/2}}(p-1)\hat{\xi}(p)\frac{1}{n}\\
\hspace*{-4pt}&&\qquad\hspace*{48.7pt}{}\times\sum_{n^{1/8}\leq m_1\leq
(1-\eta)n} \bigl(g\bigl((m_1,n-m_1,0,\ldots)/n\bigr)-\eta\bigr)\frac{m_1}{n}\\
\hspace*{-4pt}&&\qquad\hspace*{48.7pt}\hspace*{74.8pt}{}\times
\frac{1}{\sigma\sqrt{2\pi}((m_1/n)(1-m_1/n))^{3/2}}
\\
\hspace*{-4pt}&&\qquad\hspace*{48.7pt}\hspace*{74.8pt}{}\times\sum_{(1-\eta')n-m_1\leq m_2\leq n-m_1}\PP(\tau_{p-2} = n-m_1-m_2).
\end{eqnarray*}\looseness=0
Finally, the last sum is $\PP(\tau_{p-2}\in[0,\eta'n])\geq
\PP(\tau_{\lfloor\eps\sqrt{n}\rfloor}\in[0,\eta' n])$, and this
can be made larger than $1-\eta$ when $n$ is large enough, by
choosing $\eps$ even smaller than before if necessary. Indeed, as is
well known, and again a consequence of~\eqref{eq:16},
$\tau_{\lfloor\sqrt{a}\rfloor}/a$ converges\vspace*{1pt} in distribution as
$a\to\infty$ to a stable random variable with index $1/2$. Taking the
$\liminf$ in $n$ and using a convergence of Riemann sums, yields
\begin{eqnarray*}
&&\liminf_{n\to\infty}\sqrt{n}\ov
{q}^*_n\bigl(g\ind_{\{x_1<1-\eta,x_1+x_2>1-\eta'\}}\bigr)\\
&&\qquad\geq(1-\eta)^6(1-\eta'/\eta)\\
&&\qquad\quad{}\times\sum_{p\geq
2}(p-1)\hat{\xi}(p)\int_0^{1-\eta}
\frac{\d
x}{\sigma\sqrt{2\pi}x^{1/2}(1-x)^{3/2}}\bigl(g(x,1-x,0,\ldots)-\eta\bigr).
\end{eqnarray*}
One concludes using the fact that $\sum_{ p\geq
2}(p-1)\hat{\xi}(p)=\sigma^2$.
\end{pf}

We can now finish the proof of Proposition
\ref{sec:proof-theor-refs-2}. Simply write
\[
\bigl|\ov{q}_n^*(g)-\ov{q}_n^*\bigl(g\ind_{\{x_1<1-\eta,x_1+x_2>1-\eta'\}
}\bigr)\bigr|\leq
\ov{q}_n^*\bigl(|g|\ind_{\{x_1\geq
1-\eta\}}\bigr)+\ov{q}_n^*\bigl(|g|\ind_{\{x_1+x_2\leq1-\eta'\}}\bigr) .
\]
Now
fix $\eps>0$, and using Lemmas~\ref{sec:proof-theor-refs-3} and
\ref{sec:proof-theor-refs-5}, choose $\eta,\eta'$ in such a way that
$\sqrt{n}\ov{q}_n^*(|g|\ind_{\{x_1\geq1-\eta\}})\leq\eps/2$ and
\[
\biggl|\sqrt{n}\ov{q}_n^*\bigl(g\ind_{\{x_1<1-\eta,x_1+x_2>1-\eta'\}}\bigr)-
\frac{\sigma}{\sqrt{2\pi}}
\int_0^1\frac{g((x,1-x,0,\ldots))}{x^{1/2}(1-x)^{3/2}}\,\d x
\biggr|\leq
\eps/2
\]
for every $n$ large. For this choice of $\eta,\eta'$, we then
have for every $n$ large enough,
\[
\biggl|\sqrt{n}\ov{q}_n^*(g)-\frac{\sigma}{\sqrt{2\pi}}
\int_0^1\frac{g((x,1-x,0,\ldots))}{x^{1/2}(1-x)^{3/2}}\,\d x
\biggr|\leq
\eps+ \sqrt{n}\ov{q}_n^*\bigl(|g|\ind_{\{x_1+x_2\leq1-\eta'\}}\bigr) ,
\]
and the
upper-bound converges to $\eps$ as $n\to\infty$ by Lemma
\ref{sec:proof-theor-refs-4}. Since $\eps$ was arbitrary, this proves
Proposition~\ref{sec:proof-theor-refs-2}, hence implying Theorem
\ref{sec:main-results-1} in Case 1.

%s5.2 ###
\subsection{Stable case}\label{sec:stable-case}

Assume that $\xi(p)\sim cp^{-\alpha-1}$ for some $\alpha\in(1,2)$ and
$c>0$. Theorem~\ref{sec:main-results-1} in Case 2 will follow
if we can show that hypothesis (H) holds for
$\gamma=1-1/\alpha$, $\ell\equiv(\alpha(\alpha-1)/(c
\Gamma(2-\alpha)) )^{1/\alpha}$ and the dislocation measure
$\nu_\alpha$.
%Since the latter is a.s. carried by proper sequences,
A~similar reasoning as in the beginning of the previous section shows
that it suffices to prove the following statement.
\begin{proposition}\label{sec:proof-theor-refs-7}
If $f\dvtx\mathcal{S}_1\to\R$ is a continuous function bounded by $1$,
and $g(\mathbf{x})= (1-\max\mathbf{x})f(\mathbf{x})$, then
\[
n^{1-1/\alpha}\ov{q}_n^*(g)\mathop{\longrightarrow}_{n\to\infty}
\biggl(c\frac{\Gamma(2-\alpha)}{\alpha(\alpha-1)}\biggr)^{1/\alpha}
\nu_\alpha^*(g) .
\]
\end{proposition}

One will note that the function $g$ of the statement is continuous
$\nu_\alpha^*$-a.e., since $\mathbf{x}\mapsto\max\mathbf{x}$ is
continuous at every
point $\mathbf{x}$ with sum $1$. Now,
%
%e31 ###
\begin{eqnarray} \label{eq:15}
\ov{q}_n^*(g)&=&\sum_{p\geq
1}q_n\bigl(p(\lambda)=p\bigr)\E\biggl[g\biggl(\frac{(X_1^*,\ldots
,X_p^*,0,\ldots)}{n}\biggr)
\bigg| \tau_p = n\biggr]\nonumber\\ &=&n^{1/\alpha}\int_0^\infty
\d
x q_n\bigl(p(\lambda)=\lceil n^{1/\alpha}x \rceil\bigr)\\
&&\qquad\quad\hspace*{5.2pt}{}\times
\E\biggl[g\biggl(\frac{(X_1^*,\ldots,X_{\lceil n^{1/\alpha}
x\rceil}^*,0,\ldots)}{n}\biggr) \bigg| \tau_{\lceil
n^{1/\alpha} x\rceil} = n\biggr].\nonumber
\end{eqnarray}
Recall the notation around~(\ref{eq:14}). The random walk $S_n$ is now
such that $(S_{\lfloor nt\rfloor}/n^{1/\alpha}$, $t\geq0)$ converges in
distribution in the Skorokhod space to a spectrally positive stable
L\'{e}vy process $(Y_t,t\geq0)$ with index $\alpha$ and L\'{e}vy
measure $c\,\d
x/x^{1+\alpha}\ind_{\{x>0\}}$. Its Laplace transform is given by
$\E[\exp(-\lambda Y_t)]=\exp(tc'\lambda^\alpha)$, where
$c'=c\frac{\Gamma(2-\alpha)}{\alpha(\alpha-1)}$. The
Gnedenko--Kolmogorov local limit theorem also yields
\[
n^{1/\alpha}\PP(S_n=k)=p_1(k/n^{1/\alpha})+\eps(n,k),
\]
where
$\sup_k|\eps(n,k)|\to0$ as $n\to\infty$, and $p_t$ is the density
of $Y_t$.
This, together with~(\ref{eq:14}) and our
hypothesis on the asymptotic behavior of $\xi$, entails that
\[
q_n\bigl(p(\lambda)=\lceil n^{1/\alpha}x\rceil\bigr)\sim
cn^{-1}x^{-\alpha}\frac{p_1(-x)}{p_1(0)} .
\]
%
%=\frac{Cn^{-1}}{p_1(0)}x^{-\alpha-1}Q_x(1)
%.$$
%
%On the other hand, classical convergence results for random walks
%toward stable subordinators entail the convergence in distribution in
%the Skorokhod space
%$$(\frac{\tau_{\lfloor n^{1/\alpha} x\rfloor}}{n}, x\geq
%0)\stackrel{(d)}{\underset{n\to\infty}{\longrightarrow}}(T_x,x
%0) ,$$ where
Let us now focus on the random variables $X_1,X_2,\ldots$ and
$\tau_p=X_1+\cdots+X_p$. We have $\PP(X_1 = n) = n^{-1}\PP
(S_n=-1)\sim
n^{-1-1/\alpha}p_1(0)$, which gives that $X_1$ is in the domain of
attraction of a stable random variable with index $1/\alpha$. More
specifically, one has that $(\tau_{\lfloor n
x\rfloor}/n^{\alpha},x\geq0)$ converges in the Skorokhod space to a
stable subordinator $(T_y,y\geq0)$ with index $1/\alpha$ and L\'{e}vy
measure
%
%e32 ###
\begin{equation}\label{eq:28}
p_1(0)\,\frac{\d
x}{x^{1+1/\alpha}}\ind_{\{x>0\}} .
\end{equation}
Its Laplace transform is given by
\[
\E[\exp(-\lambda T_x)]=\exp\bigl(-xp_1(0)\alpha
\Gamma(1-1/\alpha)\lambda^{1/\alpha}\bigr) .
\]
On the other hand, $T_x$
has same distribution as the first hitting time of $-x$ by $(Y_t,t\geq
0)$ (because a similar statement is true of $\tau_p$ and $S_n$), which
identifies the Laplace exponent of $T_1$ as $(\lambda/c')^{1/\alpha}$,
see~\cite{BertoinLevy}, Chapter VII. This yields
%
%e33 ###
\begin{equation}\label{eq:29}
p_1(0)=\frac{1}{ \alpha\Gamma(1-1/\alpha)(c')^{1/\alpha}}=
\frac{1}{\alpha\Gamma(1-1/\alpha)}
\biggl(\frac{\alpha(\alpha-1)}{c\Gamma(2-\alpha)}\biggr)^{1/\alpha}.
\end{equation}
Let $Q_y$ be the probability density function of $T_y$. By
\cite{BertoinLevy}, Corollary VII.3 (which is also called the cyclic
lemma) we have $tQ_x(t)=xp_t(-x)$, while the Gnedenko--Kolmogorov local
limit theorem states that
\[
p^{\alpha}\PP(\tau_p = n)=Q_1(n/p^\alpha)+\eps'(p,n) ,
\]
where
$\sup_n|\eps'(p,n)|\to0$ as $p\to\infty$. Finally, the scaling
relation $Q_x(t)=x^{-\alpha}Q_1(tx^{-\alpha})$ holds and will be
useful in the sequel.
\begin{lemma}\label{sec:proof-theor-refs-6}
The sequence $(X_1^*,\ldots,X^*_{\lceil n^{1/\alpha}x\rceil})/n$
conditioned on $\tau_{\lceil n^{1/\alpha}x\rceil} = n$ converges in
distribution to\vadjust{\goodbreak} a random sequence $(\Delta_1^*,\Delta_2^*,\ldots)$,
defined inductively by
\begin{eqnarray*}
&&\PP\Biggl(\Delta_{i+1}^*\in\d z \Big|
\Delta^*_1,\ldots,\Delta^*_i,\sum_{j=1}^i\Delta_j^*=y\Biggr)\\[-2pt]
&&\qquad=\frac
{p_1(0)x}{
z^{1/\alpha}}\frac{Q_x(1-y-z)}{Q_x(1-y)}\,\d z ,\qquad 0\leq z\leq
1-y .
\end{eqnarray*}
\end{lemma}
\begin{pf}
The case $i=0$ is obtained by using the local
limit theorem in
\begin{eqnarray*}
&&
n\PP\bigl(X_1^*=\lfloor nz\rfloor| \tau_{\lceil
n^{1/\alpha}x\rceil} = n\bigr)\\[-2pt]
&&\qquad=\lceil n^{1/\alpha}x\rceil\lfloor
nz\rfloor\PP(X_1=\lfloor nz\rfloor)\frac{\PP(\tau_{\lceil
n^{1/\alpha}x\rceil-1} = n-\lfloor nz\rfloor)}{\PP(\tau_{\lceil
n^{1/\alpha}x\rceil} = n)} .
\end{eqnarray*}
One then argues by induction, in
an elementary way. Details are left to the reader.\vspace*{-2pt}
\end{pf}

%Note that $p_1(0)$ is the limit of $Q_x(1)/x$ as $x\to0$, and this
%is equal to the density of the L\'{e}vy measure
The limiting sequence $(\Delta^*_i,i\geq1)$ has same distribution as
the sequence of jumps of the subordinator $(T_y,0\leq y\leq x)$,
conditionally given $T_x=1$, and arranged in size-biased order; see
\cite{PitmanStFl}, Chapter 4, or~\cite{BertoinBook}. We will denote by
$\Delta T_{[0,x]}^*$ this randomly ordered sequence of jumps.
%. From this, it follows that the jumps of $\tau$ converge to those
%of $T$ (both in the unconditioned and conditioned settings), so
%that $$E[f(\frac{(X_1,\ldots,X_{\lceil n^{1/\alpha}
%x\rceil})^\da}{n}) | \tau_{\lceil n^{1/\alpha} x\rceil}
%= n]
%T_{[0,x]}) | T_x=1] .$$
Hence, provided we have the right to apply dominated convergence in
(\ref{eq:15}), we obtain, using $xp_1(-x)=Q_x(1)$,
%
%e34 ###
\begin{equation}\label{eq:30}
n^{1-1/\alpha}\ov{q}{}^*_n(g)\mathop{\longrightarrow}_{n\to\infty}
\frac{c}{p_1(0)}\int_0^\infty\frac{\d x}{x^{\alpha+1}}Q_x(1)\E
\bigl[g\bigl(\Delta
T_{[0,x]}^*\bigr) | T_x=1\bigr] .
\end{equation}
Using scaling for the subordinator $(T_y,y\geq0)$,
%we know that $(T_{xy},0\leq y\leq1)$ has same law as
%$(x^{\alpha}T_y,0\leq y\leq1)$, so
the previous integral can be rewritten as
\[
\frac{c}{p_1(0)}\int_0^\infty\frac{\d
x}{x^{2\alpha+1}}Q_1(x^{-\alpha})\E\bigl[g\bigl(x^\alpha\Delta T_{[0,1]}^*\bigr)
|
T_1=x^{-\alpha}\bigr] ,
\]
and changing variables $u=x^{-\alpha}$ shows
that this is equal to
\[
\frac{c}{\alpha p_1(0)}\int_0^\infty Q_1(u)\,\d u \E\bigl[ug\bigl(\Delta
T_{[0,1]}^*/u\bigr) | T_1=u\bigr]=\frac{c}{\alpha p_1(0)}\E\bigl[T_1g\bigl(\Delta
T_{[0,1]}^*/T_1\bigr)\bigr] .
\]
Finally, the sequence $\Delta T_{[0,1]}$ of
jumps of $T$ before time $1$ is the sequence of atoms of a Poisson
measure with intensity given by~(\ref{eq:28}). Using~(\ref{eq:29}), it
thus has same distribution as
$\alpha(\alpha-1)c^{-1}\Gamma(2-\alpha)^{-1}(\Delta_1,\Delta
_2,\ldots)$, as
defined in Section~\ref{sec:appl-galt-wats}. Using the notation
therein and~(\ref{eq:29}), we get after rearrangements
\begin{eqnarray*}
&&
\frac{c}{\alpha p_1(0)}\E\bigl[T_1g\bigl(\Delta T^*_{[0,1]}/T_1\bigr)\bigr]\\[-2pt]
&&\qquad=
\biggl(c\frac{\Gamma(2-\alpha)}{ \alpha(\alpha-1)}\biggr)^{1/\alpha}
\frac{\alpha^2\Gamma(2-1/\alpha)}{ \Gamma(2-\alpha)}\E\biggl[
Tg\biggl(\frac{\Delta_i^*}{T},i\geq
1\biggr)\biggr]\\[-2pt]
&&\qquad=\biggl(c\frac{\Gamma(2-\alpha)}{
\alpha(\alpha-1)}\biggr)^{1/\alpha}\nu_\alpha^*(g)
\end{eqnarray*}
as wanted. It remains to justify that the convergence~(\ref{eq:30})
is indeed dominated. To this end, using~(\ref{eq:15}) and the fact
that $q_n(p(\lambda)=\lceil n^{1/\alpha}x\rceil)\leq C\lceil
n^{1/\alpha}x\rceil^{-\alpha}$, it suffices to show that the
expectation term in this equation is bounded by $C\lceil
n^{1/\alpha}x\rceil/n^{1/\alpha}$ for some $C$ independent of $n$, and
for $x\in[0,1]$. In turn, since $|g(\mathbf{x})|\leq(1-x_1)$, it suffices
to substitute this upper-bound to~$g$. Now, we have $\PP(X_1=m)\leq
Cm^{-1-1/\alpha}$ for every $m$, so that
\begin{eqnarray*}
&&\E\biggl[\biggl(1-\frac{X_1^*}{n}\biggr)\bigg| \tau
_{\lceil
n^{1/\alpha}x\rceil} = n\biggr]\\
&&\qquad=\sum_{m=1}^n\biggl(1-\frac
{m}{n}\biggr)\lceil
n^{1/\alpha}x\rceil\frac{m}{n}\PP(X_1=m)\frac{\PP(\tau_{\lceil
n^{1/\alpha}x\rceil-1} = n-m)}{\PP(\tau_{\lceil
n^{1/\alpha}x\rceil} = n)}\\
&&\qquad\leq
\sum_{m=1}^n\biggl(1-\frac{m}{n}\biggr)\lceil n^{1/\alpha}x\rceil
\frac{m}{n}\PP(X_1=m)\\
&&\qquad\hphantom{\leq
\sum_{m=1}^n}{}\times\frac{(({\lceil
n^{1/\alpha}x\rceil-1})/({n-m}))\PP(S_{n-m}=-\lceil
n^{1/\alpha}x\rceil+1) }{({\lceil
n^{1/\alpha}x\rceil}/{n})\PP(S_{n}=-\lceil
n^{1/\alpha}x\rceil)}\\
&&\qquad\leq C\frac{\lceil
n^{1/\alpha}x\rceil}{n^{1/\alpha}}
\frac{1}{n}\sum_{m=1}^n\frac{1}{({m/n}
)^{1/\alpha}
(1-{m/n})^{1/\alpha}} \\
&&\qquad\leq C\frac{\lceil
n^{1/\alpha}x\rceil}{n^{1/\alpha}} ,
\end{eqnarray*}
where we have used that $(n-m)^{1/\alpha}\PP(S_{n-m}=-\lceil
n^{1/\alpha}x\rceil+1)$ is uniformly bounded (in $n,m,x$) and that
$n^{1/\alpha}\PP(S_{n}=-\lceil n^{1/\alpha}x\rceil)$ is uniformly
bounded away from $0$ for $x\in[0,1]$. This is the wanted bound,
concluding the proof of Proposition~\ref{sec:proof-theor-refs-7},
hence of Theorem~\ref{sec:main-results-1}.

%%%%%%%%%%%%%%%%%%%%%%%%%%%
%s6 ###
\section{Scaling limits of uniform unordered trees}
\label{SectionPolya}
%%%%%%%%%%%%%%%%%%%%%%%%%%%

In this section, we fix once and for all an integer $m \in
\{2,\ldots,\infty\}$ and consider trees in which every vertex has at
most $m$ children. We use the notation of Section
\ref{sec:appl-polya-trees} and let $T_n$ be uniformly distributed in
$\bT^{(m)}_{n}$, for $n \geq1$.

The first difficulty we have to overcome is that the sequence $(T_n,n
\geq1)$ \textit{is not} Markov branching as defined in Section
\ref{sec:mark-branch-trees-1}. We will therefore start in
Section~\ref{Polyacoupling} by coupling this sequence with a family of Markov
branching trees that are asymptotically close to $T_n,n\geq1$, and then
check in Section~\ref{PolyaH2} that the coupled trees satisfy (H).

Let us fix some notation. For $\mathsf{t}\in\bT_n^{(m)}$, we can write
$\mathsf{t}=\langle\mathsf{t}^{(1)},\ldots,\mathsf{t}^{(k)}\rangle
$ with
$\sum_{i=1}^k\#\mathsf{t}^{(i)} = n-1$, and we let $\lambda(\mathsf
{t})\in\P_{n-1}$
be the partition obtained by arranging in decreasing order the
sequence $(\#\mathsf{t}^{(1)},\ldots,\#\mathsf{t}^{(k)})$ (of
course, this does not
depend on the labeling of the trees $\mathsf{t}^{(1)},\ldots,\mathsf
{t}^{(k)}$).
Let $\bF_j(k)$ be the set of\vadjust{\goodbreak} multisets\footnote{Recall that a
multiset with $k$ elements in some
set $A$ is an element of the quotient set $A^k/\mathfrak{S}_k$,
where $\mathfrak{S}_k$ acts in the natural way by permutation of
components.} with $k$ elements in $\bT^{(m)}_j$. By convention, we
set $\bF_j(0)=\{\varnothing\}$.
Then, for $\lambda\in\P_{n-1}$ with $p(\lambda)\leq m$, we have a
bijection
%
%e35 ###
\begin{equation}\label{eq:31}
\bigl\{\mathsf{t}\in
\bT_n^{(m)}\dvtx\lambda(\mathsf{t})=\lambda\bigr\}\equiv\prod
_{j=1}^{n-1}\bF_j(m_j(\lambda))
,
\end{equation}
obtained by grouping the subtrees of $\mathsf{t}$ born from
the root with size $j$ into a multiset, denoted by $\mathsf
{f}_j(\mathsf{t})$, of
$m_j(\lambda)$ trees. From this, we deduce that $\mathsf
{f}_j(T_n),1\leq
j\leq n-1$, are independent uniform random elements in
$\bF_j(m_j(\lambda))$ conditionally given~$\lambda(T_n)$.
However,\vspace*{-1pt}
the uniform random element in $\bF_j(k)$ has a different distribution
from the multiset induced by $k$ i.i.d. uniform elements in $\bT^{(m)}_j$,
as soon as $k\geq2$. This is what prevents $T_n$ from enjoying the
Markov branching property, that is, from having law $\bQ^q_n$, where
for $n\geq1$,
$q_n$ is the law of $\lambda(T_{n+1})$.

Letting $\bFF_j(k)=\#\bF_j(k)$, the previous bijection yields
\[
\mathbf S_n^{(\lambda)}:=\#\bigl\{\mathsf{t}\in
\bT^{(m)}_n\dvtx\lambda(\mathsf{t})=\lambda\bigr\} = \prod_{j=1}^{n-1}
\bFF_j(m_j(\lambda)) .
\]
When $p(\lambda) > m$, we set $\mathbf
S_n^{(\lambda)}=0$. Of course, letting $\bfT^{(m)}_n=\#\bT^{(m)}_n$,
we also have
\[
\bfT^{(m)}_n = \sum_{\lambda\in\mathcal P_{n-1}}
\mathbf S_n^{(\lambda)} .
\]
Using the obvious fact that $\bFF_j(k)\leq\bfT^{(m)}_j \bFF_j(k-1)$,
we obtain the rough but useful bound
%
%e36 ###
\begin{equation}
\label{bound1}
\mathbf S_{n}^{(\lambda)} \leq\bfT^{(m)}_{\lambda_1}
\mathbf S_{n-\lambda_1}^{(\lambda_2,\lambda_3,\ldots,\lambda
_{p(\lambda)})}.
\end{equation}

We recall the key result~(\ref{eq:18}) of Otter
\cite{otter48}, which is used throughout the proofs below:
\[
\mathbf T_n^{(m)} \mathop{\sim}_{n\to\infty} \kappa_m
\frac{\rho_m^n}{n^{3/2}}.
\]
Setting $\bfT^{(m)}_0=1$ by convention, we obtain that for
$\rho=\rho_m>1$, and two constants $K\geq1\geq\mathrm{k}>0$,
%
%e37 ###
\begin{equation}\label{eq:20}
\bfT^{(m)}_n \leq K \frac{\rho^n}{n^{3/2}} ,\qquad n\geq0,\qquad
\bfT^{(m)}_n \geq\mathrm{k} \frac{\rho^n}{n^{3/2}} ,\qquad n \geq
1 .
\end{equation}
Note that we also have $\bfT^{(m)}_n \leq K \rho^n$ for all $n \geq
0$. Last, we let $\kappa=\kappa_m$.

%%%%%%%%%%%%%
%s6.1 ###
\subsection{Coupling}
\label{Polyacoupling}
%%%%%%%%%%%%%

Let $\varrho_n$ be the uniform probability
distribution over $\bT^{(m)}_n$, and let $q_n=\lambda_*\varrho
_{n+1}$ be
the law of the partition of $n$ induced by the subtrees born from the root\vadjust{\goodbreak}
of a $\varrho_{n+1}$-distributed tree. For every $n\geq1$, we want to
construct a pair of random variables $(T_n,T'_n)$ on some probability
space $(\Omega,\FF,\PP)$, such that:
\begin{itemize}
\item
$T_n$ has law $\varrho_n$;
\item
$T'_n$ has law $\bQ^q_n$;
\item
for every $\eps>0$,
$\lim_{n\to\infty}\mathbb E
[d_{\mathrm{GHP}}(n^{-\varepsilon}T_n,n^{-\varepsilon}T'_n ) ]=0$.
\end{itemize}

Recall that if $T_n$ has distribution $\varrho_n$, and conditionally on
$\lambda(T_n)=\lambda$, then $\mathsf{f}_j(T_n),1\leq j\leq n-1$, are
independent, respectively, uniform in $\bF_j(m_j(\lambda))$.
We are going to need the following fact.
\begin{lemma}\label{sec:coupling-1}
For every $j,k\geq1$, let $F_j$ be uniform in $\bF_j(k)$ and
$\overline F_j$
be the multiset induced by an i.i.d. sequence of $k$ random variables
with law $\varrho_j$. Let $A_j$ be the set of elements in $\bF_j(k)$ with
components that are pairwise distinct. Then:

\begin{longlist}
\item
one has $\PP(F_j\in A_j)\leq\PP(\overline F_j\in A_j)$;

\item
the conditional distributions of $F_j$ and $
\overline F_j$
given $A_j$ are equal.
\end{longlist}
\end{lemma}
\begin{pf}
For a finite set $A$, the number of multisets
with $k$ elements is $\#(A^k/\mathfrak{S}_k)\geq\#A^k/k!$. Then
\begin{eqnarray*}
\PP(F_j\in
A_j)&=&\frac{\#\bT^{(m)}_j(\#\bT^{(m)}_j-1)\cdots(\#\bT
^{(m)}_j-k+1)}{k!\#\bF_j(k)}\\ &\leq&
\frac{\#\bT^{(m)}_j(\#\bT_j^{(m)}-1)\cdots(\#\bT^{(m)}_j-k+1)}{(\#
\bT^{(m)}_j)^k}\\ &=&\PP( \overline F_j\in
A_j) .
\end{eqnarray*}
This gives (i). Property (ii) is also obtained by counting: on the
event $A_j$, the probability that $F_j$ equals some given (multi)set
$S\in\bF_j(k)$ with all distinct elements is $\#\bF_j(k)^{-1}$, while
the probability that $\overline F_j$ equals the same set $S$ is
$k!(\#\bT^{(m)}_j)^{-k}$. Dividing by $\PP(F_j\in A_j)$ and $\PP
(\overline F_j\in A_j)$,
respectively, gives the same result.
\end{pf}

The previous statement allows us to construct a coupling between $F_j$
and $ \overline F_j$, in the following way.
Let $\mathsf{f}\in\bF_j(k)$.
Consider three independent random
variables $\mathsf{f}'',\mathsf{f}''' ,B$,
such that the law of
$\mathsf{f}''$ is the law of $ \overline F_j$ conditionally given
$A_j$, the law of $\mathsf{f}'''$ is the law of $ \overline F_j$
conditionally given $A^c_j$ and $B$ is an
independent Bernoulli random variable with $\PP(B=1)=\PP( \overline
F_j\in
A_j^c)/\PP(F_j\in A_j^c)$, which is indeed in $[0,1]$ by (i) in Lemma
\ref{sec:coupling-1}.
Set
\[
\mathsf{f}'=\cases{
\mathsf{f}, &\quad if
$\mathsf{f}\in A_j$,\cr
\mathsf{f}'', &\quad if $\mathsf{f}\notin A_j$ and $B=0$,\cr
\mathsf{f}''', &\quad if $\mathsf{f}\notin A_j$ and $B=1$.}
\]
We let $K_j(\mathsf{f},\cdot)$ be the law of the multiset $\mathsf{f}'$
thus obtained, hence defining a Markov kernel on $\bF_j(k)$. We say
that the random variables $F,F'$ are \textit{naturally coupled} if
$(F,F')$ has law $\varrho(\d\mathsf{f})K_j(\mathsf{f},\d\mathsf
{f}')$, where $\varrho$ is the
law of $F$ on $\bF_j(k)$. Using (ii) in Lemma~\ref{sec:coupling-1}, it
is then easy to obtain the next result.
\begin{lemma}\label{sec:coupling-2}
If $F_j$ is uniform in $\bF_j(k)$, and $(F_j,F'_j)$ are naturally
coupled, then the law of $F'_j$ is that of the multiset induced by $k$
i.i.d. uniform elements in $\bT_j$.
\end{lemma}

Next, we define a Markov kernel $K(\mathsf{t},\cdot)$ on $\bT
^{(m)}$, in an
inductive way. Let $K(\bullet,\{\bullet\})=1$. Assume that the measure
$K(\mathsf{t},\cdot)$ on $\bT_{\#\mathsf{t}}^{(m)}$ has been
defined for every
$\mathsf{t}\in\bT^{(m)}$ with $\#\mathsf{t}\leq n-1$. Take $\mathsf
{t}\in\bT^{(m)}_n$, and
let $\lambda=\lambda(\mathsf{t}),p=p(\lambda)$. Let $\mathsf
{f}_j(\mathsf{t})\in
\bF_j(m_j(\lambda)),1\leq j\leq n-1$, be the multisets of trees born
from the root of $\mathsf{t}$, respectively, with size $j$. Let
$\mathsf{f}_j'(\mathsf{t})$
be independent random multisets, respectively, with law
$K_j(\mathsf{f}_j(\mathsf{t}),\cdot)$. We relabel\vspace*{1pt} the $p$ elements
of the multisets
$\mathsf{f}'_j(\mathsf{t}),1\leq j\leq n-1$, as $\mathsf
{t}_{(1)},\ldots,\mathsf{t}_{(p)}$, in
nonincreasing order of size, so that $\#\mathsf{t}_{(i)}=\lambda_i$---if
there is some $j$ with $m_j(\lambda)\geq2$, we arrange the trees with
same size in exchangeable random order. All these trees are in
$\bT^{(m)}$ and have at most $n-1$ vertices. By the induction
hypothesis, conditionally on this family, we can\vspace*{1pt} consider another
family $\mathsf{t}_{(1)}',\ldots,\mathsf{t}_{(p)}'$ of independent
trees with
respective laws $K(\mathsf{t}_{(i)},\cdot)$. Let $K(\mathsf{t},\cdot
)$ be the law of
the tree $\langle\mathsf{t}'_{(i)},1\leq i\leq p \rangle$. This procedure
allows us to define the Markov kernel $K(\mathsf{t},\cdot)$ for every
tree in
$\bT^{(m)}$.

We say that the random trees $(T,T')$, defined on a common probability
space, are \textit{naturally coupled} if the law of $(T,T')$ is
$\varrho(\d\mathsf{t})K(\mathsf{t},\d\mathsf{t}')$, where
$\varrho$ is the law of $T$. Is is easy
to see that for every random variable $T$ on $\bT^{(m)}$ with law
$\varrho$,
then, possibly at the cost of enlarging the probability space
supporting $T$, one can construct a random variable $T'$ so that
$(T,T')$ is naturally coupled.
\begin{proposition}\label{Lemmapetito}
Let $T_n$ have law $\varrho_n$ and $(T_n,T'_n)$ be naturally coupled.
Endow these trees, respectively,
with the measures $\mu_{T_n}$ and $\mu_{T'_n}$.
Then:

\begin{longlist}
\item
the tree $T'_n$ has distribution $\bQ^q_n$, for
every $n\geq1$;

\item for all $a>0$, the Gromov--Hausdorff--Prokhorov
distance between $aT_n$ and $aT_n'$ is at most $2aj^*$ where $j^*$ is
the supremum integer $j\geq1$ so that there exist two subtrees of
$T_n$ with size $j$, born from the same vertex and which are equal
(with the convention $\sup\varnothing=0$);

\item for all $\varepsilon>0$, $\mathbb E
[d_{\mathrm{GHP}}(n^{-\varepsilon}T_n,n^{-\varepsilon}T'_n ) ]
\rightarrow0$ as $n \rightarrow\infty$.
\end{longlist}
\end{proposition}
\begin{pf}
We prove (i) by induction. For $n=1$ the property is
obvious. Assume the property holds for every index up to $n-1$, and
condition on $\lambda(T_n)=\lambda$, which by definition has
probability $q_{n-1}(\lambda)$. As noticed before Lemma
\ref{sec:coupling-1}, the multisets $F_j=\mathsf{f}_j(T_n),1\leq
j\leq n-1$,
are independent, respectively, uniform in
$\bF_j(m_j(\lambda))$. Conditionally on $F_j,1\leq j\leq n-1$, let
$F'_j,1\leq j\leq n-1$, be independent with respective laws
$K_j(F_j,\cdot)$. By Lemma~\ref{sec:coupling-2}, we obtain that $F'_j$
is the multiset induced by a sequence of $m_j(\lambda)$ i.i.d. random
variables, with law $\varrho_j$. Consequently,\vspace*{1pt} if we relabel the elements
of $F'_j,1\leq j\leq n-1$, as $T_{(1)},\ldots,T_{(p)}$ in
nonincreasing order of size (and exchangeable random order for trees
with same size), then we obtain that these trees are independent,
respectively, with distribution $\varrho_{\lambda_j}$. Since,\vspace*{1pt} by
definition of $K$, the natural coupling $(T_n,T'_n)$ is obtained by
letting $T'_n=\langle T'_{(i)},1\leq i\leq p\rangle$ where
$(T_{(i)},T'_{(i)})$ are naturally\vspace*{1pt} coupled, we readily obtain the
Markov branching property, with branching laws $(q_n,n\geq1)$.

For (ii), we again apply an induction argument. The statement is
trivial for \mbox{$n=1$}. Now, in the first step of the natural coupling, the
action of the Markov kernel $K_j$ on $\mathsf{f}_j(T_n)$ leaves it unchanged
if $\mathsf{f}_j(T_n)\in A_j$, that is, if there are no ties in the multiset
$\mathsf{f}_j(T_n)$. Consequently, with the same notation as in the previous
paragraph, a subtree of $T_n$ born from the root that appears with
multiplicity $1$ will also appear among $T_{(1)},\ldots,T_{(p)}$.

Moreover, subtrees that are replaced are always replaced by trees with
the same number of vertices and a tree with $j$ vertices and
edge-lengths $a$ has height at most~$aj$. So the
Gromov--Hausdorff--Prokhorov distance between two trees with
edge-lengths $a$ that both decompose above the root in subtrees of
same size $j$ is at most $2aj$ (it is implicit in this proof that all
trees are endowed with the uniform measure on their vertices). We now
appeal to the following elementary:

\textit{Fact}. Let $\mathsf{t},\mathsf{t}'$ be such that
$k=p(\lambda(\mathsf{t}))=p(\lambda(\mathsf{t}'))$, and let
$\mathsf{t}=\langle
\mathsf{t}_{(1)},\ldots,\mathsf{t}_{(k)}\rangle$ and $\mathsf
{t}'=\langle
\mathsf{t}'_{(1)},\ldots,\mathsf{t}'_{(k)}\rangle$ with $\# \mathsf
{t}_{(i)}=\# \mathsf{t}'_{(i)}$ for $1 \leq i \leq k$. Then for every $a>0$,
\[
d_{\mathrm{GHP}}(a\mathsf{t},a\mathsf{t}')\leq\max_{1\leq i\leq
k}d_{\mathrm{GHP}}\bigl(a\mathsf{t}_{(i)},a\mathsf{t}_{(i)}'\bigr) .
\]

From this we deduce that the Gromov--Hausdorff--Prokhorov distance
between $aT_n$ and $aT'_n$ is at most
\[
(2a\sup\{1\leq j \leq
n-1\dvtx F_j\in A_j^c\} )\vee\mathop{\sup_{1\leq j\leq n-1}}_{F_j\in
A_j}\sup_{i\dvtx \#T_{(i)}=j}d_{\mathrm{GHP}}\bigl(aT_{(i)},aT'_{(i)}\bigr)
,
\]
where $(T_{(i)},T'_{(i)})$ is the natural coupling. The induction
hypothesis allows us to conclude.

Last, for (iii), fix $\varepsilon\in(0,1)$. The
Gromov--Hausdorff--Prokhorov distance between $n^{-\varepsilon}T_n$ and
$n^{-\varepsilon}T'_n$ is bounded from above by $2n^{1-\varepsilon}$
for all $n \geq1$. Next, for $\gamma\in(0,1)$, let $A_n^{\gamma}$
be the subset of trees of $\bT^{(m)}_{n}$ that have at least two
subtrees born from the same vertex that are identical, and with size
larger than $n^{\gamma}$. By (ii), when $T_n \notin A_n^{\gamma}$,
$d_{\mathrm{GHP}}(n^{-\varepsilon}T_n ,n^{-\varepsilon}T'_n) \leq
2n^{\gamma-\varepsilon}$. Hence,
\begin{eqnarray*}
\mathbb E
[d_{\mathrm{GHP}}(n^{-\varepsilon}T_n,n^{-\varepsilon}T'_n )
] &=&\mathbb E
\bigl[d_{\mathrm{GHP}}(n^{-\varepsilon}T_n,n^{-\varepsilon}T'_n )
\mathbf1_{\{T_n \in A^{\gamma}_n\}}\bigr]\\
&&{}+ \mathbb
E
\bigl[d_{\mathrm{GHP}}(n^{-\varepsilon}T_n,n^{-\varepsilon}T'_n )
\mathbf1_{\{ T_n \notin A^{\gamma}_n \}}\bigr] \\ &\leq
& 2n^{1-\varepsilon}\mathbb P(T_n \in A_n^{\gamma})+ 2n^{\gamma
-\varepsilon}.
\end{eqnarray*}
Taking $\gamma<\varepsilon$ and using Lemma~\ref{negligible}
following right below, we get the
result.\vspace*{-2pt}
\end{pf}
\begin{lemma}
\label{negligible}
For $\gamma\in(0,1)$, let $A_n^{\gamma}$ be the subset of $\bT
^{(m)}_{n}$ of trees $\mathsf t$ that have at least one vertex $v$ such
that at least two subtrees born from $v$ are equal and have at least
$n^{\gamma}$ vertices. Then,
\[
\mathbb P(T_n \in A_n^{\gamma})=O(\rho^{-n^{\gamma}}n^{5/2})
\qquad\mbox{as } n
\rightarrow\infty.\vspace*{-2pt}
\]
\end{lemma}

This lemma will be an easy consequence of the following result.
For every tree $\mathsf{t}$ and any vertex $v$ of $\mathsf{t}$, we
let $\mathsf t^{(v)}$
denote the subtree of $\mathsf t$ rooted at~$v$.
When $v^*$ is taken uniformly at random among the vertices of
$\mathsf t$, we set $\mathsf t^{(*)}:=\mathsf t^{(v^*)}$.\vspace*{-2pt}
\begin{lemma}
\label{lemmuniform}
The distribution of $T_n^{(*)}$ conditionally on $\# T_n^{(*)}=k$ is
uniform on $\bT^{(m)}_k$, for every $1 \leq k \leq n$.\vspace*{-2pt}
\end{lemma}

Note that the event $\{\# T_n^{(*)}=k\}$ has a strictly positive
probability for all \mbox{$1 \leq k
\leq n$}.\vspace*{-2pt}
\begin{pf*}{Proof of Lemma~\ref{lemmuniform}}
Let $k \in\{1,\ldots,n\}$. For all
$\mathsf{t}_0
\in\bT^{(m)}_k$, using that $\mathbb P(T_n =\mathsf t)=1/\mathbf
T_n^{(m)}$ for all $\mathsf{t}\in\bT_n^{(m)}$,
\begin{eqnarray*}
\mathbb P\bigl(T_n^{(*)} = \mathsf{t}_0 \bigr)&=&
\sum_{\mathsf t \in\bT^{(m)}_n} \mathbb
P\bigl(T_n^{(*)} = \mathsf t_0 | T_n =\mathsf t\bigr)
\mathbb P(T_n =\mathsf t) \\[-1pt]
&=& \frac{1}{\mathbf T^{(m)}_n}\sum
_{\mathsf t \in
\bT^{(m)}_n} \mathbb P\bigl(\mathsf{t}^{(*)}=\mathsf{t}_0\bigr)\\[-1pt]
&=&\frac{1}{n\bfT^{(m)}_n}\sum_{\mathsf{t}\in\bT^{(m)}_n,v\in
\mathsf{t}}\ind_{\{\mathsf{t}^{(v)}=\mathsf{t}_0\}}.
\end{eqnarray*}
This quantity is independent of $\mathsf{t}_0 \in\bT^{(m)}_k$
because there is an obvious bijection
between the sets $\{(\mathsf{t},v)\dvtx\mathsf{t}\in\bT^{(m)}_n,v\in
\mathsf{t},\mathsf{t}^{(v)}=\mathsf{t}_0\}$ and
$\{(\mathsf{t},v)\dvtx\mathsf{t}\in\bT^{(m)}_n,v\in\mathsf{t},\mathsf
{t}^{(v)}=\mathsf{t}_1\}$ for $\mathsf{t}_1\in\bT^{(m)}_k$.
Hence the result.\vspace*{-2pt}
\end{pf*}
\begin{pf*}{Proof of Lemma~\ref{negligible}} Let
$A_{n}^{\gamma}(k)$ be the subset of trees of $ \bT^{(m)}_k$ whose
decomposition above the root gives birth to at least two identical
subtrees with size\vadjust{\goodbreak} greater than $ n^{\gamma}$, $k \leq n$. We first
give an upper bound for the probability $\mathbb P(T_k \in
A_n^{\gamma}(k))$. To do so, we bound from above the number of trees
of $ \bT^{(m)}_k$ that decompose in at least two identical subtrees of
size $i $ [$i \leq(k-1)/2$]: there are $ \mathbf T^{(m)}_i $ choices
for the tree with size $i$ appearing twice. Then, there are $\mathbf
T^{(m)}_{k-2i} $ forests with $k-1-2i$ vertices. Gluing the twin trees
and a forest with $k-1-2i$ vertices to a common root gives a tree with
$k$ vertices (and a root branching in possibly more than $m$ subtrees)
and all trees in $\bT^{(m)}_{k}$ with at least two subtrees with size
$i$ can be obtained in this way. From this we deduce that the
cardinality of $A_{n}^{\gamma}(k)$ is at most $\sum_{i =
n^{\gamma}}^{(k-1)/2} \mathbf T^{(m)}_i \mathbf T^{(m)}_{k-2i} $. In
particular, using~(\ref{eq:20}) and the fact that $\rho>1$,
\begin{eqnarray*}
\mathbb P\bigl(T_k \in A_{n}^{\gamma}(k) \bigr) &\leq&\frac{1}{
\mathbf T^{(m)}_k }\sum_{i = n^{\gamma}}^{(k-1)/2} \mathbf T^{(m)}_i
\mathbf T^{(m)}_{k-2i}\\
&\leq& C \frac{k^{3/2}}{\rho^k} \sum_{i =
n^{\gamma}}^{(k-1)/2}\frac{\rho^{i}\rho^{k-2i}}{i^{3/2}} \\
&\leq& C
n^{3/2}\rho^{-n^{\gamma}},
\end{eqnarray*}
where $C$ is a generic constant independent of $n$ and $k \leq
n$. Now, in the following lines, given $T_n$, we let $v_1,v_2,\ldots,v_n$
denote its vertices labeled uniformly at random,
\begin{eqnarray*}
\mathbb P(T_n \in A_n^{\gamma} ) &\leq& \mathbb E\biggl[
\sum_{v \in T_n}\mathbf1_{\{T_n^{(v)} \in A^{\gamma}_{n}(\#
T_n^{(v)} )\}}\biggr] \\
&=&\mathbb E\Biggl[ \sum_{i=1}^n\mathbf
1_{\{T_n^{(v_i)} \in A_{n}^{\gamma}(\# T_n^{(v_i)} )\}}\Biggr] \\ &=& n
\mathbb P\bigl( T_n^{(*)} \in A_{n}^{\gamma}\bigl(\#
T_n^{(*)} \bigr)\bigr) \\
&&\hspace*{-30.58pt}\displaystyle \mathop{=}_{\mathrm{by}\ \mathrm{Lemma}\mbox{ }
\mbox{\fontsize{8.36}{8.36}\selectfont{\ref{lemmuniform}}}} n \sum_{k=1}^n
\mathbb
P\bigl( T_k \in A_{n}^{\gamma}(k) \bigr) \mathbb P\bigl(\#
T_n^{(*)}=k\bigr) \\ &\leq& C
n^{5/2}\rho^{-n^{\gamma}} \sum_{k=1}^n \mathbb P\bigl(\#
T_n^{(*)}=k \bigr) \\
&=& C n^{5/2}\rho^{-n^{\gamma}}.
\end{eqnarray*}
\upqed
\end{pf*}

%%%%%%%%%%%%%%%%%%%%%%%%%%
%s6.2 ###
\subsection{Hypothesis \textup{(H)} and conclusion}
\label{PolyaH2}
%%%%%%%%%%%%%%%%%%%%%%%%%%

It remains to check that the family of probability distributions
on $\mathcal P_{n}, n \geq1$, defined by
\[
q_n(\lambda)=\mathbb P\bigl(\lambda(T_{n+1})=\lambda
\bigr)=\mathbb
P\bigl(\lambda(T'_{n+1})=\lambda\bigr)=\frac{\mathbf
S_{n+1}^{(\lambda)}}{\mathbf T^{(m)}_{n+1}}\qquad
\forall\lambda\in\mathcal P_{n}
\]
satisfies the assumption
(H) with $\gamma=1/2$, $\ell\equiv1$ and $\nu$
proportional to the Brownian dislocation measure $\nu_2$. For this we
recall and fix some more notation:
\begin{enumerate}[$\ast$]
\item[$\ast$] $\wt\bT^{(m)} _{n} $ is the subset of $\bT^{(m)}_{n}$
of trees with root degree less or equal to $m-2$;
\item[$\ast$] $\wt{\mathbf{T}}^{(m)} _{n} $ is the cardinality of
$\wt{\bT}^{(m)} _{n} $, and $\psi^{(m)}(x)=\sum_{n\geq
1}\bfT^{(m)}_nx^n$, $\wt{\psi}^{(m)}(x)=\sum_{n\geq
1}\wt{\bfT}^{(m)}_nx^n$;
\item[$\ast$] for $\lambda=(\lambda_1,\lambda_2,\ldots) \in
\mathcal
P_n$, $\lambda_{\mathrm r}:=\sum_{i=3}^{\infty}\lambda_i = n
-\lambda_1-\lambda_2$.
\end{enumerate}
The main result of this section follows.
\begin{proposition}
\label{PropH2}
For all $2 \leq m \leq\infty$, and all continuous functions $f \dvtx
\mathcal S^{\downarrow} \rightarrow\mathbb R$ such that $|f(\mathbf
s)| \leq1-s_1$ for $\mathbf s \in\mathcal S^{\downarrow}$,
\[
\sqrt n \sum_{\lambda\in\mathcal P_n}
f\biggl(\frac{\lambda}{n}\biggr) \frac{\mathbf
S^{(\lambda)}_{n+1}}{\mathbf T_{n+1}^{(m)}}
\mathop{\longrightarrow}_{n\to\infty}\kappa\wt{\psi
}^{(m)}(1/\rho)
\int_{1/2}^1 \frac{f(x,1-x,0,\ldots)
}{x^{3/2}(1-x)^{3/2}} \,\mathrm dx.
\]
\end{proposition}

Note that $\wt{\psi}^{(m)}(1/\rho)
$ is finite, since $\wt{\mathbf{T}}^{(m)} _{n} \leq\mathbf
T^{(m)}_{n} \leq K \rho^{n}/n^{3/2}$. This sum is explicit in terms of
$\kappa$ and $\rho$ when $m=2$ or $m=\infty$. See Section \ref
{sec:appl-polya-trees} for details.

With this proposition, it is easy to conclude the proof of Theorem
\ref{sec:main-results-2}. Indeed, together with Theorem
\ref{sec:main-result-2} and Proposition~\ref{Lemmapetito}(i), it
leads to the convergence
\[
\frac{1}{\sqrt n} T'_n
\mathop{\longrightarrow}^{(\mathrm{d})}_{n\to\infty}
c_m\mathcal T_{1/2,\nu_{2} }
\]
for the Gromov--Hausdorff--Prokhorov topology, where $c_m=\sqrt
2/\break(\sqrt\pi\kappa\wt{\psi}^{(m)}(1/\rho))$. Then, by Proposition
\ref{Lemmapetito}(iii) and since $(\mathcal M_{\mathrm w},
d_{\mathrm{GHP}})$ is a complete separable space, we can apply a
Slutsky-type theorem to get
\[
\frac{1}{\sqrt n} T_n \mathop{\longrightarrow}^{(\mathrm{d})}_{n\to\infty}
c_m \mathcal T_{1/2,\nu_{2} }.
\]

The rest of this section is devoted to the proof of Proposition
\ref{PropH2}.

%s6.2.1 ###
\subsubsection{Negligible terms}\label{prelim}

We show in this section that the set
of partitions $\lambda\in\mathcal P_n$ where either $\lambda_1\geq
n(1-\varepsilon)$ or
$\lambda_{\mathrm r} \geq n\varepsilon$ plays a negligible role
in the limit of Proposition~\ref{PropH2} when we first let $n
\rightarrow\infty$ and then
$\varepsilon\rightarrow0$.
\begin{lemma}
\label{lem5}
There exists $C \in(0,\infty)$ such that, for all
$0<\varepsilon<1$,
\[
\limsup_{n \rightarrow\infty} \frac{\sqrt n}{\mathbf T_{n+1}^{(m)}}
\sum_{\lambda\in\mathcal P_n} \mathbf1_{\{\lambda_1 \geq
n(1-\varepsilon)\}}\biggl(1-\frac{\lambda_1}{n}\biggr) \mathbf
S^{(\lambda)}_{n+1}\leq\frac{C\sqrt\varepsilon}{(1-\varepsilon)^{3/2}}.
\]
\end{lemma}
\begin{pf}
Using~(\ref{bound1}) and then (\ref
{eq:20}), we get
\begin{eqnarray*}
&&
\sum_{\lambda\in\mathcal P_n, \lambda_1 \geq n(1-\varepsilon)}
\biggl(1-\frac{\lambda_1}{n}\biggr) \mathbf S^{(\lambda)}_{n+1}\\
&&\qquad\leq
\sum_{\lambda_1 = \lceil n(1-\varepsilon) \rceil}^n
\biggl(1-\frac{\lambda_1}{n}\biggr) \mathbf T_{\lambda_1}^{(m)}
\sum_{\mu\in\mathcal P_{n-\lambda_1}}
\mathbf S_{n+1-\lambda_1}^{(\mu)} \\
&&\qquad\leq \sum_{\lambda_1=
\lceil
n(1-\varepsilon) \rceil}^n \biggl(1-\frac{\lambda_1}{n}\biggr)
\mathbf T^{(m)}_{\lambda_1} \mathbf T^{(m)}_{n+1-\lambda_1} \\
&&\qquad \leq
K^2
\rho^{n+1} \sum_{\lambda_1= \lceil n(1 -\varepsilon) \rceil}^{n-1}
\frac{1-{\lambda_1}/{n}}{\lambda_1^{3/2}(n+1-\lambda_1)^{3/2}}
\\
&&\qquad\leq \frac{K^2 \rho^{n+1}}{(n(1-\varepsilon))^{3/2}}\times
\frac{1}{n^{3/2}} \times\sum_{\lambda_1= \lceil n(1 -\varepsilon)
\rceil}^{n-1} \frac{1}{(1-\lambda_1/n)^{1/2}}.
\end{eqnarray*}
We conclude with the fact that the sum $ \sum_{\lambda_1= \lceil n(1
-\varepsilon)
\rceil}^{n-1} (1-\lambda_1/n)^{-1/2}$ is smaller than the integral $
\int_{n(1-\varepsilon)}^n (1-x/n)^{-1/2} \,\mathrm dx = 2n \sqrt
\varepsilon$ and then use the lower bound of~(\ref{eq:20}) for
$\mathbf T_{n+1}^{(m)}$.
\end{pf}

To deal with the partitions where $\lambda_{\mathrm r} \geq
n\varepsilon$, we need the following lemma when $m=\infty$. We denote
by $\mathbf T^{(\infty,a-)}_{n}$ the number of trees of
$\bT^{(\infty)}_{n}$ whose subtrees born from the root have sizes at
most $a$, $a \geq1$.
\begin{lemma}
\label{lemmadeg}
Let $m=\infty$. There exists
$A,B>0$ such that
\[
\mathbf T^{(\infty,a-)}_{k+1} \leq
A\rho^{k}\exp(-B k/a )\qquad \forall k \in\mathbb N \mbox{ and
} a \geq1.
\]
\end{lemma}
\begin{pf}
Recall that $\bT$ denotes the set of all
(rooted, unordered, unlabeled) trees and rewrite the power series
$\psi=\psi^{(\infty)}$ as $ \psi(x)=\sum_{\mathsf t \in\bT} x^{\#
\mathsf t}$. According to~\cite{FlSe09}, Section VII.5, its radius
of convergence is $1/\rho<1$ and $\psi(1/\rho)=1$. Note also that
$\psi(0)=0$. Now, we consider a random tree $T$ in $\bT$ with
distribution defined by
\[
\mathbb P(T=\mathsf t)=\rho^{-\# \mathsf t}.
\]
If $c_{\varnothing}(\mathsf t)$ denotes the degree of the
root of $\mathsf t$, we just have to show that
%
%e38 ###
\begin{equation}
\label{eqdeg}\quad
\mathbb P\bigl(c_\varnothing(T)=r\bigr) \leq A' \exp(-B'r)\qquad
\mbox{for some } A',B'>0 \mbox{ and all }r \geq1.
\end{equation}
Indeed, each tree with $k+1$ vertices and a decomposition in subtrees
with sizes at most $a$ has a root degree larger or equal to $k/a$.
So,
if the above inequality holds, we will have
\begin{eqnarray*}
\mathbf T^{(\infty,a-)}_{k+1} &\leq&\rho^{k+1} \mathbb
P\bigl(c_\varnothing(T)\geq k/a , T \in\bT^{(\infty)}_{k+1} \bigr)\\
&\leq&\rho^{k+1} A' B'^{-1} \exp\bigl(-B'(k/a-1)\bigr)
\end{eqnarray*}
as required. To get
(\ref{eqdeg}), note that
\begin{eqnarray*}
\mathbb P\bigl(c_\varnothing(T)=r\bigr)&=&\sum_{\mathsf t \in\bT,
c_\varnothing(\mathsf t)=r} \rho^{-\# \mathsf t} \\
&=& \sum_{k=1}^r
\frac{1}{ k !} \mathop{\sum_{\mathsf{t}_1,\ldots,\mathsf{t}_k \in
\bT}}_{\mathrm{pairwise}\ \mathrm{distinct} }\mathop{\sum_{m_1+\cdots+m_k=r}}_{m_i
\geq1} \rho^{-1-\sum_{ 1 \leq i\leq k} m_i \# t_i},
\end{eqnarray*}
which is obtained by considering the multiset of $r$ subtrees of a
tree $\mathsf t$, made of $k$ distinct trees with multiplicities
$m_1,\ldots,m_k$. Hence,
\begin{eqnarray*}
\mathbb P\bigl(c_\varnothing(T)=r\bigr) &\leq& \rho^{-1} \sum_{k=1}^r
\frac{1}{ k !} \mathop{\sum_{m_1+\cdots+m_k=r}}_{m_i \geq1}
\prod_{i=1}^k\psi( \rho^{-m_i}) \\ & =& \rho^{-1} \sum
_{k=1}^{\lfloor
cr \rfloor} \frac{1}{ k !} \mathop{\sum_{m_1+\cdots+m_k=r}}_{m_i
\geq
1} \prod_{i=1}^k\psi( \rho^{-m_i})\\
&&{} + \rho^{-1} \sum_{k={\lfloor
cr \rfloor+1}}^r \frac{1}{ k !} \mathop{\sum_{m_1+\cdots+m_k=r}}_
{m_i \geq1} \prod_{i=1}^k\psi( \rho^{-m_i}),
\end{eqnarray*}
where the $c \in\ ]0,1[$ chosen for this split will be specified below.

We first bound from above the second term. Using that
$\psi(\rho^{-m_i}) \leq\psi(\rho^{-1})=1$ for $m_i \geq1$ and that
$\sum_{m_1+\cdots+m_k=r , m_i \geq1} 1 ={r-1\choose k-1}$, we obtain
\begin{eqnarray*}
\sum_{k={\lfloor cr \rfloor+1}}^r \frac{1}{ k !}
\mathop{\sum_{m_1+\cdots+m_k=r}}_{m_i \geq1} \prod_{i=1}^k\psi(
\rho^{-m_i}) &\leq&\frac{1}{\lfloor cr \rfloor!} \sum_{k=1}^r
\pmatrix{r-1\cr k-1}\\
&\leq&
\frac{2^{r-1}}{\lfloor cr \rfloor!},
\end{eqnarray*}
which decays exponentially fast as $r \rightarrow\infty$, for
every $ c \in\ ]0,1[$.

Now\vspace*{1pt} we will check that the sum $ \sum_{k=1}^{\lfloor cr \rfloor}
\frac{1}{ k !} \sum_{m_1+\cdots+m_k=r , m_i \geq1} \prod
_{i=1}^k\psi(
\rho^{-m_i}) $ also decays exponentially in $r$, provided that $c
\in\ ]0,1[$ is chosen sufficiently small.\vadjust{\goodbreak} Since $\psi(0)=0$, we
have that $\psi(x) \leq Cx$ for some $C<\infty$ and all $x
\in[0,\rho^{-1}]$. Hence
\begin{eqnarray*}
&&\sum_{k=1}^{\lfloor cr \rfloor} \frac{1}{ k !}
\mathop{\sum_{m_1+\cdots+m_k=r}}_{m_i \geq1} \prod_{i=1}^k\psi(
\rho^{-m_i}) \\
&&\qquad\leq\sum_{k=1}^{\lfloor cr \rfloor} \frac{C^k}{ k !}
\mathop{\sum_{m_1+\cdots+m_k=r}}_{m_i \geq1} \prod_{i=1}^k
\rho^{-m_i} \\
&&\qquad\leq \exp(C) \sum_{k=1}^{\lfloor cr \rfloor}
\rho^{-r} \pmatrix{r-1\cr k-1} \\
&&\hspace*{-15.6pt}\qquad\displaystyle\mathop{\leq}_{\mathrm{for}\ \mathrm{all}\
\lambda>0}
\exp(C) \rho^{-r}
\sum_{k=r-\lfloor cr \rfloor}^{r-1}
\pmatrix{r-1\cr k} \exp\bigl(\lambda k- \lambda(r-\lfloor cr
\rfloor)\bigr) \\
&&\qquad \leq \exp(C) \bigl(\rho^{-1} \exp\bigl(-\lambda
(1-c)\bigr) \bigl(\exp(\lambda)+1\bigr)\bigr)^r.
\end{eqnarray*}
When $c \rightarrow0$, $\rho^{-1} \exp(-\lambda(1-c)) (\exp
(\lambda
)+1) \rightarrow\rho^{-1} (1+\exp(-\lambda))$, which is strictly
smaller than $1$ for $\lambda$ large enough. Hence, fix such a large
$\lambda$ and then take $c>0$ sufficiently small so that $\rho^{-1}
\exp(-\lambda(1-c)) (\exp(\lambda)+1)<1$. This ends the proof.
\end{pf}
\begin{lemma}
\label{lem6}
For all $\varepsilon>0$,
\[
\frac{\sqrt n}{\mathbf T_{n+1}^{(m)}} \sum_{\lambda\in\mathcal P_n}
\mathbf1_{\{\lambda_{\mathrm r} \geq n \varepsilon\}}
\mathbf S^{(\lambda)}_{n+1} \mathop{\longrightarrow}_{n\to\infty}0.
\]
\end{lemma}
\begin{pf}
$\bullet$ If $m=2$, $\lambda_{\mathrm r}
=0$ for all
$\lambda\in\mathcal P_n$ and the assertion is obvious.

$\bullet$ Assume now that $3 \leq m<\infty$, and note that
when $\lambda\in\mathcal P_n$ with $p(\lambda) \leq m$, one has
that $\lambda_{\mathrm r} \geq
n\varepsilon$ implies $(m-2) \lambda_3 \geq n \varepsilon$, in
particular $\lambda_1 \geq\lambda_2 \geq n\varepsilon/(m-2)$. Hence,
\[
\sum_{\lambda\in\mathcal P_n} \mathbf1_{\{\lambda_{\mathrm r} \geq n
\varepsilon\}} \mathbf S^{(\lambda)}_{n+1}
\leq\sum_{\lambda_{\mathrm r} = \lceil n\varepsilon\rceil}^{n-2}
\sum_{\lambda_1=\lceil
n\varepsilon/(m-2 \rceil)}^{\lfloor n-\lambda_{\mathrm r}
-n\varepsilon/(m-2) \rfloor^+} \mathbf T^{(m)}_{\lambda_1} \mathbf
T^{(m)}_{n-\lambda_{\mathrm r} -\lambda_1}
\mathbf T^{(m)}_{\lambda_{\mathrm r} +1}.
\]
Then for $C$ a generic constant, using~(\ref{eq:20}), the latter term
multiplied by $\sqrt n / \mathbf T_{n+1}^{(m)}$ is bounded from above by
\begin{eqnarray*}
&& Cn^{1/2} (n+1)^{3/2} \sum_{\lambda_{\mathrm r} =\lceil
n\varepsilon\rceil}^{n-2} \sum_{\lambda_1= \lceil
n\varepsilon/(m-2)\rceil}^{\lfloor n-\lambda_{\mathrm r}
-n\varepsilon/(m-2) \rfloor^+} \frac{1}{
\lambda_1^{3/2}(n-\lambda_{\mathrm r} -\lambda_1)^{3/2}(\lambda
_{\mathrm r} +1)^{3/2}}
\\
&&\qquad\leq C
\frac{n^2(n-2)(\lfloor n-n\varepsilon/(m-2)\rfloor^+)}{n^{3*3/2}}
=O\biggl( \frac{1}{\sqrt n}\biggr).
\end{eqnarray*}

$\bullet$ Next we turn to the case where $m=\infty$. Let
$\gamma\in(5/6,1)$. On the one hand,
by the same token as for the $m<\infty$ cases,
\begin{eqnarray*}
&&
\frac{\sqrt n}{\mathbf T^{(\infty)}_{n+1}}\sum_{\lambda\in\mathcal
P_n} \mathbf1_{\{\lambda_{\mathrm r} \geq n
\varepsilon\}} \mathbf1_{\{\lambda_2 \geq n^{\gamma} \}}
\mathbf S^{(\lambda)}_{n+1} \\
&&\qquad\leq C n^2 \sum_{\lambda_{\mathrm r}
= \lceil
n\varepsilon\rceil}^{n-2} \sum_{\lambda_1= \lceil n^{\gamma}
\rceil}^{\lfloor n-\lambda_{\mathrm r} -n^{\gamma} \rfloor^+} \frac{1}{
\lambda_1^{3/2}(n-\lambda_{\mathrm r} -\lambda_1)^{3/2}(\lambda
_{\mathrm r} +1)^{3/2}}
\\
&&\qquad\leq C
\frac{n^4}{n^{3\gamma+ 3/2}} =O(
n^{5/2-3\gamma})=o(1),
\end{eqnarray*}
since $5/2-3\gamma<0$ when $\gamma>5/6$. On the other hand, since
$\lambda_2 < n^{\gamma}$ implies that $\lambda_i <n^{\gamma}$ for
all $i \geq3$,
we get by using Lemma~\ref{lemmadeg} that
\begin{eqnarray*}
\frac{\sqrt n}{\mathbf T^{(\infty)}_{n+1}}\!\sum_{\lambda\in\mathcal
P_n}\! \mathbf1_{\{\lambda_{\mathrm r} \geq n
\varepsilon\}} \mathbf1_{\{\lambda_2< n^{\gamma} \}}
\mathbf S^{(\lambda)}_{n+1} &\leq&\frac{\sqrt n}{\mathbf T^{(\infty
)}_{n+1}}\! \sum_{\lambda_{\mathrm r} = \lceil n\varepsilon\rceil}^{n-2}
\!\sum_{\lambda_1= 1}^{n-\lambda_{\mathrm r} -1} \!\mathbf T_{\lambda
_1}^{(\infty)} \mathbf T^{(\infty)}_{n-\lambda_{\mathrm r} -\lambda
_1} \mathbf T^{(\infty,n^{\gamma}-)}_{\lambda_{\mathrm r} +1}
\\ &\leq& Cn^4
\exp(-Bn^{1-\gamma}\varepsilon)= o( 1 ).
\end{eqnarray*}
\upqed
\end{pf}

%s6.2.2 ###
\subsubsection{\texorpdfstring{Proof of Proposition \protect\ref{PropH2}}
{Proof of Proposition 51}} \label{preuveProp}

We rely on the following lemma. Let $\mathcal P_n^{\mathrm{bin}}$ be
the subset of $\mathcal P_n$ of partitions of $n$ with exactly two
parts.
\begin{lemma}
\label{lem4}
Let $f\dvtx\mathcal S^{\downarrow} \rightarrow\mathbb R$ be continuous.

\begin{longlist}
\item
For all $a \in\mathbb Z_+$ and all $\varepsilon
\in
(0,1)$, as $n \rightarrow\infty$,
\begin{eqnarray*}
&&\frac{\sqrt n}{ \mathbf T_{n+1}^{(m)}} \mathop{\sum_{\lambda\in
\mathcal
P^{\mathrm{bin}}_{n-a}}}_{\lambda_1 \leq n(1-\varepsilon)}
f\biggl(\frac{\lambda_1}{n},\frac{\lambda_2+a}{n},0,\ldots\biggr)
\mathbf S_{n+1-a}^{(\lambda)} \\
&&\qquad\longrightarrow\frac{\kappa}{\rho^{1+a}}
\int_{1/2}^{1-\varepsilon} \frac{
f(x,1-x,0,\ldots)}{x^{3/2}(1-x)^{3/2}}\,
\mathrm dx.
\end{eqnarray*}

\item Moreover, there exists $C_{\varepsilon} \in
(0,\infty)$ such that, for all $n \geq1$, all $0 \leq a \leq
n\varepsilon/2$ and all nonincreasing nonnegative sequences
$(a_i,i\geq1)$ with $\sum_{i \geq1} a_i \leq a/n $,
\begin{eqnarray*}
&&\biggl|\frac{\sqrt n}{\mathbf T_{n+1}^{(m)}}\mathop{\sum_{\lambda
\in\mathcal P^{\mathrm{bin}}_{n-a}}}_{\lambda_1 \leq
n(1-\varepsilon)}
f\biggl(\frac{\lambda_1}{n},\frac{\lambda_2}{n}+a_1,a_2,a_3,\ldots
\biggr)\mathbf S_{n+1-a}^{(\lambda)} \mathbf T^{(m)}_{a+1}
\biggr|\\
&&\qquad\leq\frac{C_{\varepsilon}}{(a+1)^{3/2}}.
\end{eqnarray*}
\end{longlist}
\end{lemma}
\begin{pf}
(i)
For large enough $n$,
\begin{eqnarray*}
&&\mathop{\sum_{\lambda\in\mathcal P^{\mathrm{bin}}_{n-a}}}_{
\lambda_1 \leq n(1-\varepsilon)}
f\biggl(\frac{\lambda_1}{n},\frac{\lambda_2+a}{n},0,\ldots\biggr)
\mathbf S_{n+1-a}^{(\lambda)}\\
&&\qquad=
f\biggl(\frac{1}{2}-\frac{a}{2n},\frac{1}{2}+\frac{a}{2n},0,\ldots
\biggr)\bFF_{(n-a)/2}(2) \mathbf1_{\{ n-a \ \mathrm{is}\
\mathrm{even}\}}\\
&&\qquad\quad{}+\sum_{\lambda_1=\lfloor(n-a)/2 \rfloor+1}^{\lfloor
n(1-\varepsilon) \rfloor}
f\biggl(\frac{\lambda_1}{n},1-\frac{\lambda_1}{n},0,\ldots\biggr)
\mathbf T^{(m)}_{\lambda_1} \mathbf T^{(m)}_{n-a-\lambda_1}.
\end{eqnarray*}
On the one hand, by Otter's approximation result for $\mathbf T_{(n-a)/2}^{(m)}$
\begin{eqnarray*}
\bFF_{(n-a)/2}(2) &=& \mathbf T_{(n-a)/2}^{(m)}\bigl(\mathbf
T_{(n-a)/2}^{(m)}+1\bigr) /2 \sim\kappa^2 \rho^{n-a}\bigl((n-a)/2\bigr)^{-3}/2 \\
&=&o
\bigl(\mathbf T_{n+1}^{(m)} /\sqrt n \bigr).
\end{eqnarray*}
On the other hand, still using
Otter's result, we get that for all $\eta>0$, provided that $n$ is
large enough,
\begin{eqnarray*}
&&\frac{\sqrt n}{ \mathbf T_{n+1}^{(m)}} \sum_{\lambda_1=\lfloor(n-a)/2
\rfloor+1}^{\lfloor n(1-\varepsilon) \rfloor}
f\biggl(\frac{\lambda_1}{n},1-\frac{\lambda_1}{n},0,\ldots\biggr)
\mathbf T^{(m)}_{\lambda_1} \mathbf T^{(m)}_{n-a-\lambda_1} \\
&&\qquad\leq
\frac{(\kappa+\eta)^2
}{(\kappa-\eta)\rho^{1+a}}\frac{1}{n}\sum_{\lambda_1=\lfloor(n-a)/2
\rfloor+1}^{\lfloor n(1-{\varepsilon}) \rfloor}
f\biggl(\frac{\lambda_1}{n},1-\frac{\lambda_1}{n},0,\ldots\biggr)
\frac{(n+1)^{3/2}}{\lambda_1^{3/2}}\\
&&\qquad\quad\hphantom{\frac{(\kappa+\eta)^2
}{(\kappa-\eta)\rho^{1+a}}\frac{1}{n}\sum_{\lambda_1=\lfloor(n-a)/2
\rfloor+1}^{\lfloor n(1-{\varepsilon}) \rfloor}}
{}\times\frac{n^{3/2}}
{(n-a-\lambda_1)^{3/2}} \\
&&\qquad\mathop{\longrightarrow}_{n\to\infty}
\frac{(\kappa+\eta)^2
}{(\kappa-\eta)\rho^{1+a}} \int_{1/2}^{1-\varepsilon}
\frac{f(x,1-x,0,\ldots) }{x^{3/2}(1-x)^{3/2}}\, \mathrm dx.
\end{eqnarray*}
Letting $\eta\rightarrow0$, this gives
\begin{eqnarray*}
&&\limsup_{n \rightarrow\infty}\frac{\sqrt n}{ \mathbf T_{n+1}^{(m)}}
\sum_{\lambda_1=\lfloor(n-a)/2 \rfloor+1}^{\lfloor n(1-\varepsilon)
\rfloor}
f\biggl(\frac{\lambda_1}{n},1-\frac{\lambda_1}{n},0,\ldots\biggr)
\mathbf T^{(m)}_{\lambda_1} \mathbf T^{(m)}_{n-a-\lambda_1} \\
&&\qquad\leq
\frac{\kappa}{
\rho^{1+a}}\int_{1/2}^{1-\varepsilon}\frac{ f(x,1-x,0,\ldots)
}{x^{3/2}(1-x)^{3/2}} \,\mathrm dx.
\end{eqnarray*}
We obtain the liminf similarly, hence (i).

(ii) We will use that $\mathbf S_{n+1-a}^{(\lambda)}\leq\mathbf
T_{\lambda_1}^{(m)} \mathbf T_{n-a-\lambda_1}^{(m)}$ for all $\lambda
\in\mathcal P^{\mathrm{bin}}_{n-a}$.
Recall that $f$ is continuous, hence bounded, on the compact set
$\mathcal S^{\downarrow}$. There exits then a\vadjust{\goodbreak}
generic constant $C$ independent of $n$ and $a \leq n \varepsilon/2$
such that
\begin{eqnarray*}
&& \biggl|\frac{\sqrt n}{\mathbf T_{n+1}^{(m)}}\sum_{\lambda\in
\mathcal
P^{\mathrm{bin}}_{n-a}, \lambda_1 \leq n(1-\varepsilon)}
f\biggl(\frac{\lambda_1}{n},\frac{\lambda_2}{n}+a_1,a_2,a_3,\ldots
\biggr)\mathbf S_{n+1-a}^{(\lambda)}
\mathbf T^{(m)}_{a+1} \biggr| \\
&&\qquad\leq\frac{C\sqrt n}{ \mathbf T^{(m)}_{n+1}}
\sum_{\lambda_1=\lceil(n-a)/2 \rceil}^{\lfloor n(1-\varepsilon)
\rfloor} \mathbf
T^{(m)}_{\lambda_1} \mathbf T^{(m)}_{n-a-\lambda_1} \mathbf T^{(m)}_{a+1}
\\
&&\qquad\leq \frac{C}{ (a+1)^{3/2}
}\frac{1}{n}\sum_{\lambda_1=\lceil(n-a)/2 \rceil}^{\lfloor
n(1-\varepsilon) \rfloor} \frac{(n+1)^{3/2}}{\lambda_1^{3/2}}\times
\frac{n^{3/2}}{(n-a-\lambda_1)^{3/2}} \\
&&\qquad\leq \frac{C}{ (a+1)^{3/2}
}\frac{1}{n}\sum_{\lambda_1 = \lceil n(1-\varepsilon/2)/2
\rceil}^{\lfloor n(1-\varepsilon) \rfloor}
\frac{n^{3/2}}{\lambda_1^{3/2}} \times
\frac{n^{3/2}}{(n-\lambda_1)^{3/2}},
\end{eqnarray*}
where we have used for the last inequality that $n-a \geq
n(1-\varepsilon/2)$ since $a \leq n \varepsilon/2$ and that
$n-a-\lambda_1
\geq(n-\lambda_1)/2$ since $a \leq n \varepsilon/2$ and $\lambda_1
\leq
n(1-\varepsilon)$. This upper bound is of the form $Cu_n/(a+1)^{3/2}$
where $(u_n,n \geq1)$ is a sequence independent of $a$ converging to
a finite limit as $n \rightarrow\infty$. Hence the result.
\end{pf}
\begin{pf*}{Proof of Proposition~\ref{PropH2}} By Lemmas \ref
{lem5} and~\ref{lem6}, the set
of partitions where either $\lambda_1\geq n(1-\varepsilon)$ or
$\lambda_{\mathrm r} \geq n\varepsilon/3$ will play a negligible role
in the limit when we first let $n \rightarrow\infty$ and then
$\varepsilon\rightarrow0$. Hence we concentrate on the
following sums [for $\varepsilon\in(0,1)$], where we use that for all
$\lambda\in\mathcal P_n$, $\lambda_1 \leq
n(1-\varepsilon)$ and $\lambda_{\mathrm r} \leq n \varepsilon/3$ implies
$\lambda_2>\lambda_3$:
%
%e39 ###
\begin{eqnarray}
\label{eq1}
&&
\mathop{\sum_{\lambda\in\mathcal P_n}}_{\lambda_1 \leq
n(1-\varepsilon), \lambda_{\mathrm r} \leq n\varepsilon/3}
f\biggl(\frac{\lambda}{n}\biggr)
\mathbf S^{(\lambda)}_{n+1} \nonumber\\
&&\qquad= \sum_{k=0}^{\lfloor n\varepsilon/3
\rfloor}
\mathop{\sum_{\mu\in\mathcal P_{k}}}_{p(\mu) \leq m-2}
\mathop{\sum_{\lambda\in\mathcal
P^{\mathrm{bin}}_{n-k}}}_{\lambda_1 \leq n(1-\varepsilon)}
f\biggl(\frac{\lambda_1}{n},\frac{\lambda_2+k}{n},0,\ldots\biggr)
\mathbf S^{(\lambda)}_{n-k+1} \mathbf{ S}^{(\mu)}_{k+1} \nonumber\\
&&\qquad\quad{} + \sum_{k=0}^{\lfloor n\varepsilon/3
\rfloor} \mathop{\sum_{\mu\in\mathcal P_{k}}}_{p(\mu) \leq m-2}
\mathop{\sum_{\lambda\in\mathcal
P^{\mathrm{bin}}_{n-k}}}_{\lambda_1 \leq n(1-\varepsilon)} \biggl(
f\biggl(\frac{\lambda_1}{n},\frac{\lambda_2}{n},\frac{\mu
_1}{n},\ldots\biggr)\\
&&\qquad\quad\hspace*{19.8pt}\hspace*{106.2pt}{}-f
\biggl(\frac{\lambda_1}{n},\frac{\lambda_2+k}{n},0,\ldots
\biggr)\biggr)\nonumber\\
&&\qquad\quad\hspace*{119.5pt}{}\times
\mathbf S^{(\lambda)}_{n-k+1}\mathbf{S}^{(\mu)}_{k+1}.\nonumber
\end{eqnarray}
The first sum in the
right-hand side of~(\ref{eq1}) is equal to
%
%e40 ###
\begin{equation}
\label{eq4}
\sum_{k=0}^{\lfloor n\varepsilon/3 \rfloor} \sum_{\lambda\in
\mathcal
P^{\mathrm{bin}}_{n-k} , \lambda_1 \leq n(1-\varepsilon)}
f\biggl(\frac{\lambda_1}{n},\frac{\lambda_2+k}{n},0,\ldots\biggr)
\mathbf S^{(\lambda)}_{n-k+1} \wt{\mathbf{T}}{}^{(m)} _{k+1},
\end{equation}
which, multiplied by $\sqrt n / \mathbf T_{n+1}^{(m)}$, according to
Lemma~\ref{lem4}(i) and (ii) [(ii) implies dominated convergence],
converges to
%
%e41 ###
\begin{equation}
\label{eq2}
\sum_{k=0}^{\infty} \wt{\mathbf{T}}^{(m)} _{k+1} \frac{\kappa
}{\rho^{1+k}} \int_{1/2}^{1-\varepsilon}
\frac{f(x,1-x,0,\ldots)}{x^{3/2}(1-x)^{3/2}} \,\mathrm dx.
\end{equation}
Next, let $\delta>0$. Since $f$ is continuous (hence uniformly
continuous) on the
compact set $\mathcal S^{\downarrow}$, we can choose $\eps$ small
enough so that
the absolute value of the
second sum in the right-hand side of~(\ref{eq1}) is bounded from above
by
%
%e42 ###
\begin{equation}
\label{eq6}
2 \sum_{k=0}^{\lfloor n\varepsilon/3 \rfloor} \sum_{\lambda\in
\mathcal P^{\mathrm{bin}}_{n-k}, \lambda_1 \leq n(1-\varepsilon)}
\biggl(\delta\wedge\biggl(1-\frac{\lambda_1}{n}
\biggr)\biggr) \mathbf S^{(\lambda)}_{n-k+1} \wt{\mathbf{
T}}^{(m)} _{k+1} .
\end{equation}
Similarly as above, when multiplied by $\sqrt n / \mathbf
T^{(m)}_{n+1}$, this quantity converges to
%
%e43 ###
\begin{eqnarray}
\label{eq3}
&&
2 \sum_{k=0}^{\infty}\wt{\mathbf{T}}^{(m)} _{k+1} \frac{\kappa
}{\rho^{1+k}} \int_{1/2}^{1-\varepsilon}\frac{
\delta\wedge(1-x)}{x^{3/2}(1-x)^{3/2}} \,\mathrm
dx \nonumber\\[-8pt]\\[-8pt]
&&\qquad\leq2 \sum_{k=0}^{\infty}\wt{\mathbf{T}}^{(m)} _{k+1} \frac
{\kappa}{\rho^{1+k}} \int_{1/2}^{1}\frac{
\delta\wedge(1-x)}{x^{3/2}(1-x)^{3/2}} \,\mathrm
dx\nonumber
\end{eqnarray}
by Lemma~\ref{lem4}(i) and (ii).

Now let $\eta>0$ be fixed. For $\delta$ and $\varepsilon$
sufficiently small, the
terms~(\ref{eq3}) and the limsup of Lemma~\ref{lem5} are smaller than
$\eta$, and the term~(\ref{eq2}) is in a neighborhood of radius $\eta$
of the intended limit
%
%e44 ###
\begin{equation}
\label{eq7}
\kappa\sum_{k=0}^{\infty}\frac{ \wt{\mathbf{ T}}^{(m)} _{k+1} }{
\rho^{k+1}} \int_{1/2}^1
\frac{ f(x,1-x,0,\ldots)}{x^{3/2}(1-x)^{3/2}}\, \mathrm dx.
\end{equation}
Next, such small $\delta$ and $\varepsilon$ being fixed, letting $n
\rightarrow
\infty$, and using Lem\-ma~\ref{lem6} and the convergences of
(\ref{eq4}) to~(\ref{eq2}) and of~(\ref{eq6}) to~(\ref{eq3}), we get
that $ \sqrt n \sum_{\lambda\in\mathcal P_n}
f(\frac{\lambda}{n} )
\mathbf S^{(\lambda)}_{n+1}/\mathbf T_{n+1}^{(m)} $ is indeed in a
neighborhood of radius $7\eta$ of~(\ref{eq7}) for all $n$ large
enough.
\end{pf*}

\section*{Acknowledgments}
We are deeply indebted to Jean-Fran\c{c}ois
Marckert, who suggested that the methods of~\cite{HMPW06} could be
relevant to tackle the problem of scaling limits of uniform unordered
trees. This provided the initial spark for the present work and
\cite{HaMi09}.\vadjust{\goodbreak} Thanks are also due to the referee for useful comments,
and to Shen Lin for pointing at a mistake in an earlier version of the
paper.

%suskaldyti doi

% imsref loaded by lrinkeviciute, 2011-06-27 16:21:14
% imsref loaded by lrinkeviciute, 2011-06-27 16:24:29
% imsref loaded by lrinkeviciute, 2011-07-01 15:38:54

%
\printaddresses

\end{document}